%% file: main.tex
\documentclass[]{elsarticle}
\makeatletter
\def\ps@pprintTitle{%
 \let\@oddhead\@empty
 \let\@evenhead\@empty
 \def\@oddfoot{}%
 \let\@evenfoot\@oddfoot}
\makeatother

\usepackage[english]{babel}
\usepackage[utf8x]{inputenc}
\usepackage[T1]{fontenc}
\usepackage{subfiles}

\usepackage{graphicx}
\usepackage[colorinlistoftodos]{todonotes}
\usepackage{amsfonts,amsmath,amsthm,mathtools}
\usepackage{subfiles}
\usepackage{geometry}
\usepackage{todonotes}
\usepackage{tikz}
\usepackage{caption}
\usepackage{subcaption}
\newtheorem{theorem}{Theorem}
\newtheorem{remark}{Remark}
\newtheorem{definition}{Definition}
\newtheorem{lemma}{Lemma}

\newcommand{\ve}{\varepsilon}

\newcommand{\R}{\mathbb{R}}
\renewcommand{\S}{\mathcal{S}}

\newcommand{\C}{\mathcal{C}}
\newcommand{\Q}{\mathcal{Q}}
\usepackage{float}

\newcommand{\txta}{\textnormal{a}}
\newcommand{\txtr}{\textnormal{r}}

\newcommand{\hjk}[1]{{\color{black}{#1}}}
\newcommand{\e}{\color{black}}

\begin{document}
%\maketitle

%%% \footnote{Corresponding author.}

\begin{frontmatter}
\title{Geometric analysis of Oscillations in the Frzilator model}

\author[1]{Hadi Taghvafard\fnref{HT}\corref{htf}}
\ead{taghvafard@gmail.com}
\cortext[htf]{Corresponding author.}
\author[2]{Hildeberto Jard\'on-Kojakhmetov\fnref{hjk}}
\author[3]{Peter Szmolyan\fnref{ps}}
\author[1]{Ming Cao\fnref{MC}}

\fntext[HT]{H. Taghvafard was supported by the European Research Council (ERC-StG-307207).}
\fntext[hjk]{H. Jard\'on-Kojakhmetov was supported by the Alexander von Humboldt Foundation.}
\fntext[ps]{The research of P. Szmolyan was funded by WWTF under the grant MA14-049.}
\fntext[MC]{The work of Cao was supported in part by the European Research Council (ERC-CoG-771687) and the Netherlands Organization for Scientific Research (NWO-vidi-14134).}

\address[1]{
Engineering and Technology Institute, 
Faculty of Science and Engineering,
University of Groningen, 
9747 AG Groningen, The Netherlands}
\address[2]{Faculty of Mathematics,
Technical University of Munich,
85748 München, Germany}
\address[3]{Institute for Analysis and Scientific Computing,
Technische Universität Wien,
1040 Wien, Austria}

\begin{abstract}
A biochemical oscillator model, describing developmental stage of myxobacteria, is analyzed mathematically.
Observations from numerical simulations show that in a certain range of parameters, the corresponding system of ordinary differential equations displays stable and robust oscillations.
In this work, {\e{we use geometric singular perturbation theory and blow-up method to prove the existence of a strongly attracting limit cycle.
This cycle}} corresponds to a relaxation oscillation of an auxiliary system, whose singular perturbation nature originates from the small Michaelis-Menten constants of the biochemical model. 
In addition, we give a detailed description of the {\e{structure}} of the limit cycle, and the timescales along it.
\end{abstract}


\end{frontmatter}

%{\bf keywords}
%slow-fast system, relaxation oscillations, blow-up method, geometric singular perturbation theory, myxobacteria

%%%%%%%%%%%%%%%%%%%%%%%% the body of paper %%%%%%%%%%%%%%%%%%%%%%%%%%%%
\section{Introduction}
\subfile{subfiles/Intro.tex} 
%
\section{Detailed model and preliminary analysis}\label{sec.model}

In this section we provide a preliminary analysis of the biochemical oscillator proposed in \cite{igoshin2004biochemical}. We start in Subsection \ref{mod.des} by presenting a detailed description of the model under study. Furthermore, we describe the behavior of the trajectories and the role of parameters, and propose a unification of them. Afterwards, in Subsection \ref{2para.bif}, we present a two-parameter bifurcation analysis where we clarify the nature and the role of two distinct parameters of the system. Finally, in Subsection \ref{sec:sfs} we provide a brief introduction to slow-fast systems and the main techniques for their analysis.

\subsection{Model description}\label{mod.des}
\subfile{subfiles/ModelDesc.tex}

%
%%%%%----------------------------------------
\subsection{Two-parameter bifurcation analysis}\label{2para.bif}
\subfile{subfiles/App2Para.tex}

%%%%%----------------------------------------
\subsection{Preliminaries on slow-fast systems}\label{sec:sfs}
\subfile{subfiles/SlowFast.tex}

%
% \subsection{Outline of the paper}
% \subfile{subfiles/Outline.tex}

%\todo[inline]{Do not change anything from here on}
\section{Geometric singular perturbation analysis}\label{sec:geomdes}
The goal of this section is to give the detailed analysis of the slow-fast structure of the auxiliary system \eqref{repara}.
\subsection{Layer problem and the critical manifold}\label{critical.manifold}
\subfile{subfiles/GSPA1.tex}
\subsection{Reduced problem, slow manifolds, and slow dynamics}
\subfile{subfiles/ReducedProblem.tex}
\subsection{Singular cycle}\label{singularcycle}
\subfile{subfiles/SC.tex}
\subsection{Main result} 
\subfile{subfiles/MainResult.tex}

%\subsection{Global Uniqueness}
%\subfile{subfiles/GlobalUniqueness.tex}

%%%%%----------------------------------------
\section{Blow-up analysis}\label{blowup}
\subfile{subfiles/BlowUp.tex}
\subsection{Blow-up of the non-hyperbolic line $\ell_1\times\{0\}$}\label{bue}
\subfile{subfiles/BlowUpEaxis.tex}
\subsection{Blow-up of the non-hyperbolic line $\ell_2\times\{0\}$ and a sketch of the proof of Lemma \ref{entry2}}\label{buf}
\subfile{subfiles/BlowUpFaxis.tex}
\section{Conclusions}\label{sec.dis}
\subfile{subfiles/discuss.tex}

%%%%%----------------------------------------
\appendix
\section{Range of parameter $\gamma$ in Theorem \ref{main.result}}\label{relax.gamma}
\subfile{subfiles/RelaxGamma.tex}
\section*{References}
\bibliographystyle{plain}
\bibliography{References}

\end{document}

%% file: subfiles/Intro.tex
Oscillators are ubiquitous in different fields of science such as biology \cite{winfree1967biological}, biochemistry \cite{goldbeter1991minimal,goldbeter1996biochemical}, neuroscience \cite{hodgkin1952quantitative}, medicine \cite{goodwin1965oscillatory}, and engineering \cite{van1928lxxii}.
In particular, biochemical oscillations often occur in several contexts including signaling, metabolism, development, and regulation of important physiological cell functions \cite{novak2008design}.  
In this paper, we study a biochemical oscillator model that describes the developmental stage of \emph{myxcobacteria}.  Myxcobacteria are multicellular organisms that are common in the topsoil \cite{kaiser1989multicellular}.
During vegetation growth, i.e. when food is ample, myxobacteria constitute small swarms by a mechanism called ``gliding'' \cite{igoshin2004biochemical}.
In contrast, under a starvation condition, they aggregate and initiate a complex developmental cycle during which small swarms are transformed into a multicellular single body known as ``fruiting body'', whose role is to produce spores for next generation of bacteria \cite{kaiser1989multicellular}.
During the aforementioned transition, myxobactria pass through a developmental stage called the ``ripple phase'' \cite{igoshin2004biochemical,kaiser1989multicellular}, characterized by complex patterns of waves that propagate within the whole colony.

Two genetically distinct molecular motors are concentrated at the cell poles of myxobacteria, allowing them to glide on surfaces;
these two motors are called Adventurous (A-motility) and Social (S-motility) motors, respectively \cite{igoshin2004biochemical}.
The role of the former is to push the cells forward, while the role of the latter is to pull them together.
So, in order for a cell to reverse its direction, it has to alternatively activate its A-motility (push) and S-motility (pull) motors at opposite cell poles \cite{igoshin2004biochemical}.
As a result, by forward and backward motion of myxobacteria, complex spatial wave patterns are created.
In particular, wave patterns are produced by the coordination of motion of individual cells through a direct end-to-end contact signal, the ``C-signal''.
During the ripple phase of development, the C-signaling induces reversals, while suppresses them during the aggregation stage of development.
Observations from experiments resulted in proposing a biochemical oscillator in \cite{igoshin2004biochemical}, known as the Frzilator, which acts as a ``clock'' to control reversals.

The Frzilator is detailed in Section \ref{mod.des}.
From our numerical simulations, it appears that this biochemical oscillator is robust under small variation of parameters.
More importantly, it seems that (almost) all solutions converge to a ``unique'' limit cycle.
Regarding the previous property, in \cite{taghvafard2017parameter} it has been shown that within a certain range of parameter values, (almost) all trajectories are oscillatory, the system has a finite number of isolated periodic orbits, at least one of which is asymptotically stable.
Although some biological systems may produce more than one stable periodic solution for a certain range of parameters \cite{decroly1982birhythmicity}, the coexistence between multiple stable solutions has not yet been observed \emph{experimentally} \cite{goldbeter2017dissipative}.
%On the other hand, mathematically it is of great interest to prove the convergence of (almost) all solutions to a unique limit cycle.
%Therefore this motivates us to study the behavior of solutions and their convergence.

The main contribution of this paper is to prove that, within a certain range of parameter values, {\e{there exists a \emph{strongly attracting} periodic orbit for the Frzilator.
Moreover, the detailed description of the structure of such periodic orbit is given.}} 
The methodology used to prove the aforementioned result consists first on an appropriate rescaling of the original model, which leads to a slow-fast (or two timescales) system; next, we take advantage of the two timescales of the rescaled system to develop a geometric analysis via techniques of multi-timescale dynamical systems. From the multi-timescale nature of the problem, it turns out that the limit cycle is in fact a relaxation oscillator, meaning that there are several timescales along the orbit of the oscillator. 
From an analytical point of view, the main {\e{difficulty of this analysis}} is the detailed description of a transition along two non-hyperbolic lines (see details in Section \ref{sec:geomdes}). Our analysis is based on the approach developed in \cite{kosiuk2012relaxation,kosiuk2016geometric} where similar mechanisms, leading to an attracting limit cycle in the Goldbeter minimal model \cite{goldbeter1991minimal}, have been studied.

%todo[inline]{Outline here}
The rest of this paper is organized as follows.
In Section \ref{sec.model} we introduce the model, perform some preliminary analysis on
the model, and briefly introduce the tools which we are going to use in the paper.
In Section \ref{sec:geomdes} we give the slow-fast analysis of an auxiliary system, corresponding to the original system.
More precisely, we discuss the behavior of the dynamics when $\ve\to 0$.
In Section \ref{blowup} we present the blow-up analysis of the non-hyperbolic parts.
We conclude the paper with a discussion and outlook in Section \ref{sec.dis}.

%%%%%%%%%%%%-------------------------------------------
%Thus, a biochemical oscillator model that acts as a ``clock'' has been proposed in \cite{igoshin2004biochemical} to control reversals.
%In particular, the wave patterns are produced through motion movements by the C-signal, which controls reversals of gliding direction.
%Thus, a biochemical oscillator model that acts as a ``clock'' has been proposed in \cite{igoshin2004biochemical} to control reversals. ... {\textbf{@Hadi: I do not understand this part. It seems to me that there is some important piece of information missing. Could you please try to reformulate and/or elaborate on this? Otherwise we can discuss tomorrow }}}}

%% file: subfiles/ModelDesc.tex
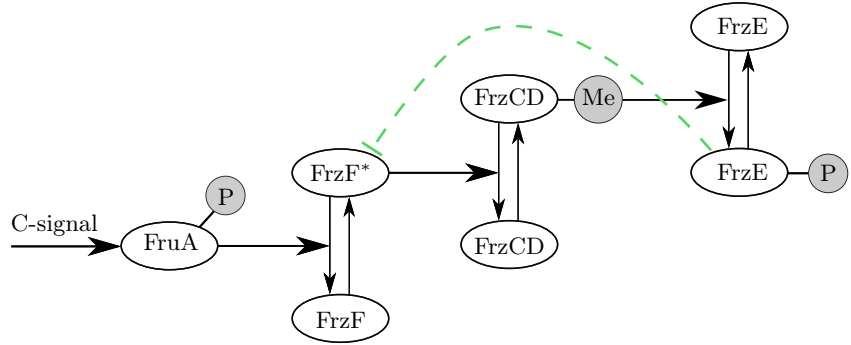
\begin{figure}[t]
\centering
\begin{tikzpicture}
	\pgftext{\includegraphics[scale=1]{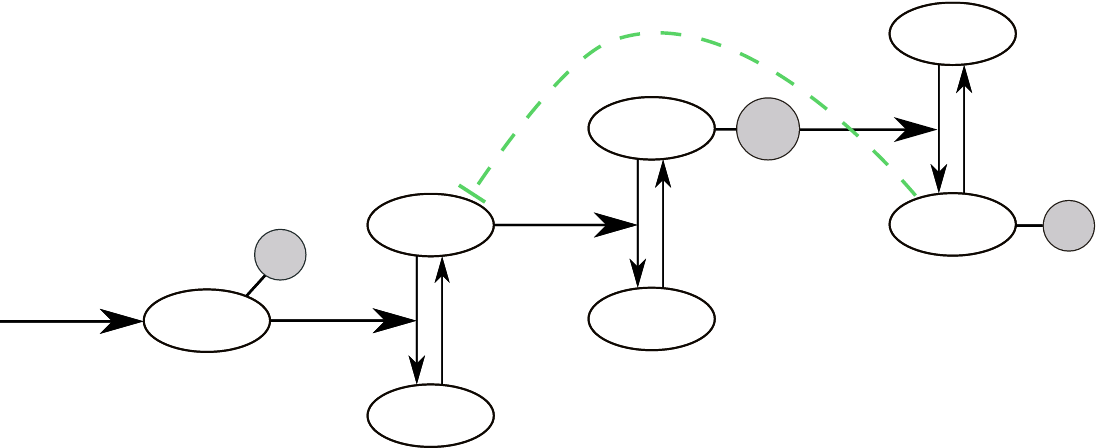}}
    \node at (-5,-0.7){\small  C-signal};
    \node at (-2.72,-0.33){\small P};
    \node at (-3.45,-0.93){\small FruA};
    \node at (-1.19,0){\small FrzF$^*$};
    \node at (-1.19,-1.97){\small FrzF};
    \node at (1.05,1){\small FrzCD};
    \node at (1.05,-0.98){\small FrzCD};
    \node at (2.23,1){\small Me};
    \node at (4.15, 1.95){\small FrzE};
    \node at (4.15, 0){\small FrzE};
    \node at (5.29, -0.01){\small P};
\end{tikzpicture}
\caption{Essential components of the Frzilator.}
\label{schematic.diag}
\end{figure}

We study a biochemical oscillator model which describes the social-behavior transition
phase of myxobacteria \cite{igoshin2004biochemical}.
This model, which is known as the Frzilator (or simply ``Frz'') model, is based on a negative feedback loop.
In the Frz model, there are three proteins, namely, a methyltransferase (FrzF), the cytoplasmic methyl-accepting protein (FrzCD), and a protein kinase (FrzE).
A direct and end-to-end collision of two myxobacteria results in producing a signal, so-called ``C-signal'', under which a protein called FruA is phosphorylated.
The signal from phosphorylated FruA (FruA-P) activates the Frz proteins as follows \cite{igoshin2004biochemical}: 
(i) the methyltransferase FrzF (FrzF$^*$) is activated by the protein FruA-P; 
(ii) in response to FrzF$^*$, the protein FrzCD is methylated (FrzCD-M); 
(iii) the phosphorylation of FrzE (FrzE-P) is activated by the methylated form of FrzCD; 
(iv) FrzF$^*$ is inhibited by the phosphorylated form of FrzE.
Figure \ref{schematic.diag} shows a schematic representation of interactions between proteins of the Frz system.
For a more detailed explanation of the model and its
biological background, see \cite{igoshin2004biochemical}.
Denote $f, c$ and $e$ respectively as the fraction of activated FrzF, methylated FrzCD, and phosphorylated FrzE. These fractions are given by \cite{igoshin2004biochemical}
\begin{equation*}
f= \frac{\textnormal{[FrzF$^*$]}}{\textnormal{[FrzF$^*$]} + \textnormal{[FrzF]}},\qquad\qquad
c= \frac{\textnormal{[FrzCD-M]}}{\textnormal{[FrzCD]} + \textnormal{[FrzCD-M]}},\qquad\qquad
e= \frac{\textnormal{[FrzE-P]}}{\textnormal{[FrzE]} + \textnormal{[FrzE-P]}}.
\end{equation*}
The interaction between the Frz proteins is modeled by Michaelis-Menten kinetics and hence leads to the dynamical system
\begin{equation}\label{model}
\begin{aligned}
	\frac{df}{d\tau} &= k_\textnormal{a}(1-f)-k_\textnormal{d}fe,\\
	\frac{dc}{d\tau} &= k_\textnormal{m}(1-c)f-k_{\textnormal{dm}}c,\\
	\frac{de}{d\tau} &= k_\textnormal{p}(1-e)c - k_{\textnormal{dp}}e,
\end{aligned}
\end{equation}
where
\begin{equation}\label{MichMen.constants}
\begin{split}
	& k_\textnormal{a}=\frac{k_\textnormal{a}^{\textnormal{max}}}{K_{\textnormal{a}} + (1 - f)}, \qquad\ \ k_\textnormal{d}=\frac{k_\textnormal{d}^{\textnormal{max}}}{K_\textnormal{d} + f},\\
	& k_\textnormal{m}=\frac{k_\textnormal{m}^{\textnormal{max}}}{K_{\textnormal{m}} + (1 - c)}, \qquad k_\textnormal{dm}=\frac{k_\textnormal{dm}^{\textnormal{max}}}{K_\textnormal{dm} + c},\\
    & k_\textnormal{p}=\frac{k_\textnormal{p}^{\textnormal{max}}}{K_{\textnormal{p}} + (1 - e)}, \qquad\ \  k_\textnormal{dp}=\frac{k_\textnormal{dp}^{\textnormal{max}}}{K_\textnormal{dp} + e}.
\end{split}
\end{equation} 

\begin{figure}[t]
\centering
\begin{tikzpicture}
	\pgftext{\includegraphics[scale=0.18]{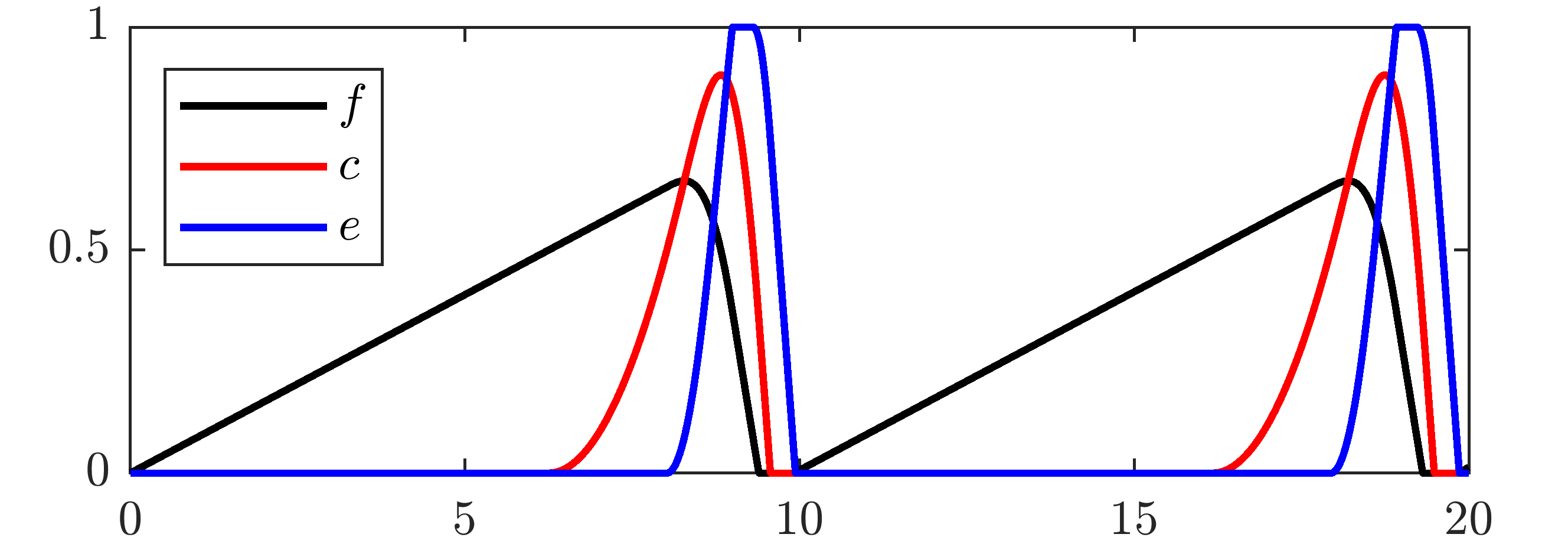}}
\end{tikzpicture}
\caption{Numerically computed time evolution of system \eqref{model} for $\ve=10^{-3}$ and $\gamma=0.08$. }
\label{TimeEvolution}
\end{figure}

\begin{remark}
Due to the fact that $f, c$ and $e$ represent fractions of active protein concentrations, their values are restricted to $[0, 1]$.
So the fraction of inactive protein concentrations are given by $(1-f)$, $(1-c)$ and $(1-e)$.
Therefore, hereafter, our analysis is restricted to the unit cube
\begin{equation}\label{cube:C}
\hjk{\Q}=\left\{ (f,c,e)\in\R^3 \, | \, f\in[0,1],\,  c\in[0,1],\,  e\in[0,1] \right\}.
\end{equation}
\end{remark}
As mentioned in \cite{igoshin2004biochemical}, the Frz system has the well-known property of ``zero-order ultrasensitivity'' which requires that the Michaelis-Menten constants $K_{\textnormal{a}}, K_{\textnormal{d}}, K_{\textnormal{m}}, K_{\textnormal{dm}}, K_{\textnormal{p}}$ and $K_{\textnormal{dp}}$ have to be \emph{small} \cite{goldbeter1984ultrasensitivity}.
It is observed \emph{numerically} in \cite{igoshin2004biochemical} that for the parameter values 
$K_{\textnormal{a}}=10^{-2}$, $K_{\textnormal{d}}=K_{\textnormal{m}}= K_{\textnormal{dm}}=K_{\textnormal{p}}=K_{\textnormal{dp}}=5\times 10^{-3}$, $k_\textnormal{d}^{\textnormal{max}}=1$ min$^{-1}$, $k_\textnormal{m}^{\textnormal{max}}=k_\textnormal{p}^{\textnormal{max}}=4$ min$^{-1}$,
$k_\textnormal{dm}^{\textnormal{max}}=k_\textnormal{dp}^{\textnormal{max}}=2$ min$^{-1}$, and $k_\textnormal{a}^{\textnormal{max}}=0.08$ min$^{-1}$, system \eqref{model} has an attracting periodic solution.
For simplicity, we unify all the dimensionless Michaelis-Menten constants by $K_{\textnormal{a}}=2K_{\textnormal{d}}=2K_{\textnormal{m}}= 2K_{\textnormal{dm}}=2K_{\textnormal{p}}=2K_{\textnormal{dp}}=\ve\ll 1$.
After unifying all Michaelis-Menten constants by $\ve$, denoting $\gamma:=k_\textnormal{a}^{\textnormal{max}}$, and substituting \eqref{MichMen.constants} in \eqref{model}, we obtain the following dynamical system 
\begin{align}\label{syseps}
	&\frac{df}{d\tau} = \frac{\gamma(1-f)}{\ve + (1-f)} - \frac{2fe}{\ve+2f},\nonumber\\
	&\frac{dc}{d\tau} = \frac{8(1-c)f}{\ve + 2(1-c)} - \frac{4c}{\ve+2c},\\
	&\frac{de}{d\tau} = \frac{8(1-e)c}{\ve + 2(1-e)} - \frac{4e}{\ve+2e}.\nonumber
\end{align}
{\e{Figures \ref{cycle} and \ref{TimeEvolution} show numerically computed attracting limit cycle as well as time evolution of system \eqref{syseps} for $\ve=10^{-3}$ and $\gamma=0.08$}}. 
\begin{remark}

%%----------------------------------
\begin{figure}[t]
\centering
\begin{tikzpicture}
	\pgftext{\includegraphics[scale=1]{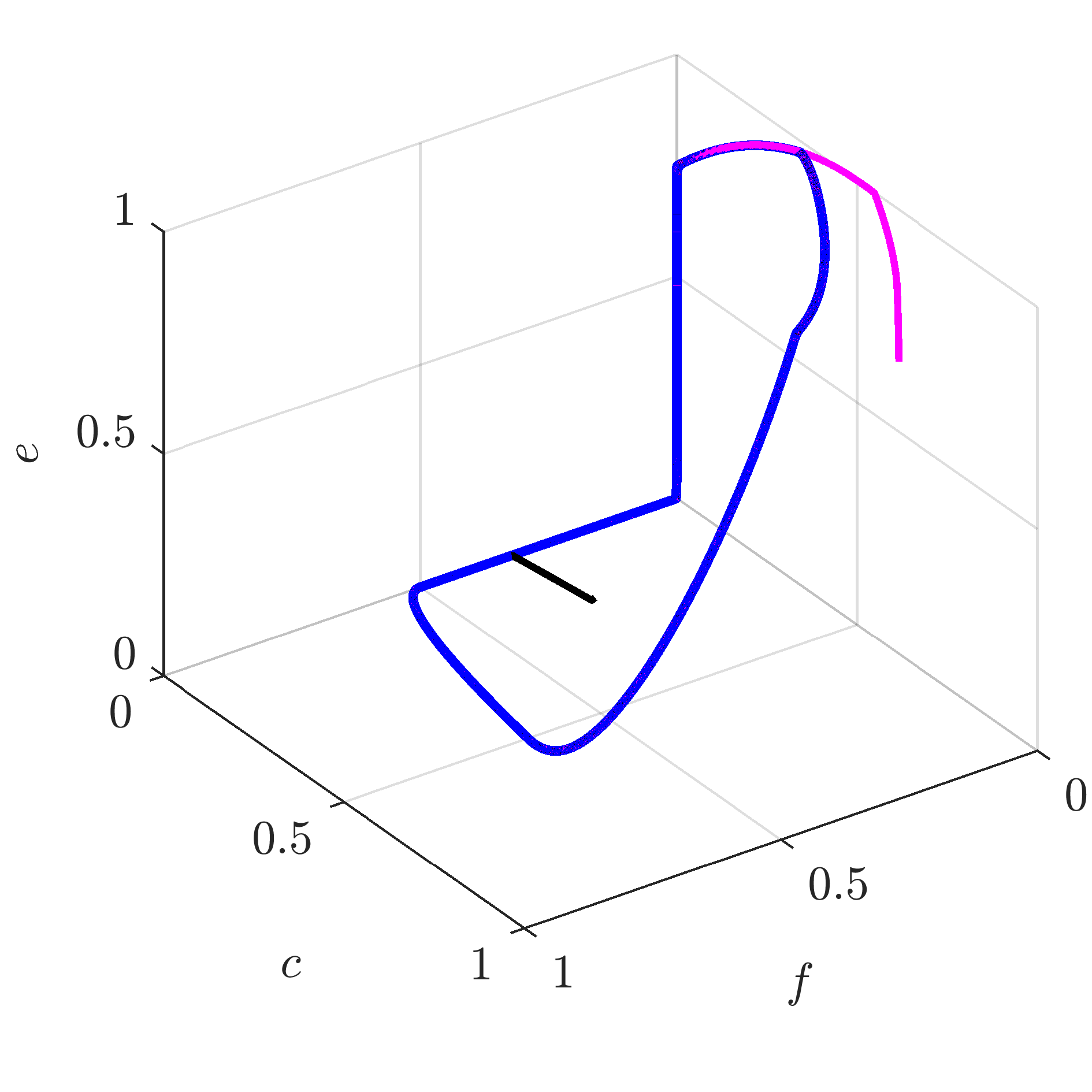}}
\end{tikzpicture}
\caption{Numerically computed attracting limit cycle of system \eqref{model} for $\ve=10^{-3}$ and $\gamma=0.08$.}
\label{cycle}
\end{figure}
%%%---------------------------------
For our analysis in this paper, we fix $\gamma=0.08$, while later we show that this parameter can be relaxed to some extent, see Remark \ref{gamma.relax} and Appendix \ref{relax.gamma}.
\end{remark}
{\e{The dynamics along the limit cycle, shown in Fig. \ref{cycle}, can be summarized as follows.}}
Initially, all protein ratios $f, c$ and $e$ are close to zero, under the dynamics \eqref{syseps}, the variable $f$ increases (due to the action of the $C$-signal), while $c$ and $e$ stay close to zero.
Once the variable $f$ passes the activation threshold $f^*:=0.5$, the variable $c$ increases very fast. 
Next, once the variable $c$ passes the threshold $c^*:=0.5$, the variable $e$ is activated and also increases very fast until it reaches its maximum value, i.e., $e=1$.
Due to the fact that there is a negative feedback from $e$ to $f$, the increase in $e$ results in the degradation of variable $f$.
Once $f$ reaches the threshold $f^*$, variable $c$ decreases, and once $c$ reaches the threshold $c^*$, the variable $e$ decreases vary fast.
As a result, the variables $f$ and $c$ reach their lowest values (i.e. very close to zero), but the variable $e$ reaches the threshold $e^*:=\gamma$.
Once the variable $e$ drops below the threshold $e^*$, the variable $f$ is activated and increases.
This behavior is repeated in a periodic manner and a limit cycle is formed (see Figure \ref{cycle}).
%%%%%%%%%%

%%%%%--------------

%%%%%%%----------

For system \eqref{syseps}, a parameter-robustness analysis with respect to $\ve$ and $\gamma=0.08$ is presented in \cite{taghvafard2017parameter}.
More precisely, using bifurcation analysis, it is shown that system \eqref{syseps} is robust under the variation of $\ve$ for $\ve\in(0, \ve^*)$ with $\ve^*:=0.05517665$.  
Moreover, it is proven that for $\ve\in(0,\ve^*)$, almost all trajectories converge to a \emph{finite} number of periodic solutions, one of which is orbitally asymptotically stable. 
%From numerical simulations, as shown in Figure \ref{cycle}, it appears that the stable limit cycle is unique, and that there are no other attractors. 
In this article, we prove the existence of a strongly attracting limit cycle which explains the numerically computed periodic orbit, for sufficiently small $\ve>0$.
%Similar mechanisms, leading to an attracting limit cycle, occur in the Goldbeter minimal model~\cite{goldbeter1991minimal}, which have been studied in~\cite{kosiuk2016geometric}.}

%\textbf{Existence and uniqueness of a stable limit cycle.} For the system given by \eqref{syseps}, and for $\ve>0$ sufficiently small, (almost all) trajectories converge to a unique and smooth limit cycle. See Theorem \ref{main.result}.

%\noindent\emph{Idea of the proof:} The convergence of almost all trajectories to at least one asymptotically stable limit cycle is proven in \cite{taghvafard2017parameter}. In this article we use techniques from singular perturbation theory to prove uniqueness of such a stable limit cycle. Moreover, we provide a precise description of the topological properties of the limit cycle, and the timescales along it.

%% file: subfiles/App2Para.tex
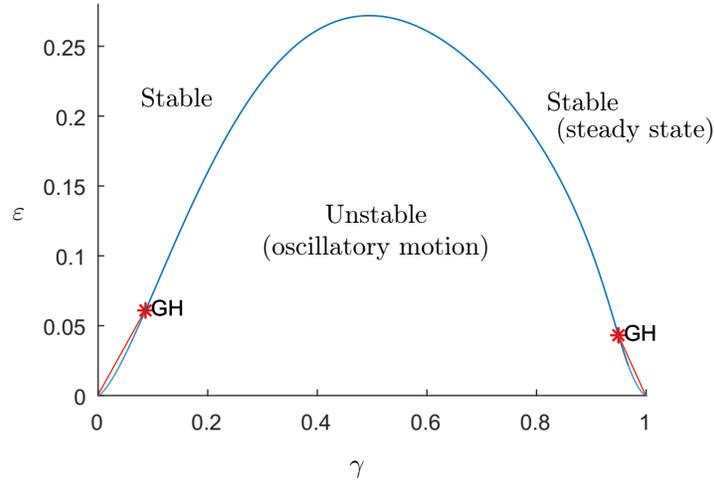
\begin{figure}[t]
\centering
\begin{tikzpicture}
\pgftext{\includegraphics[scale=0.3]{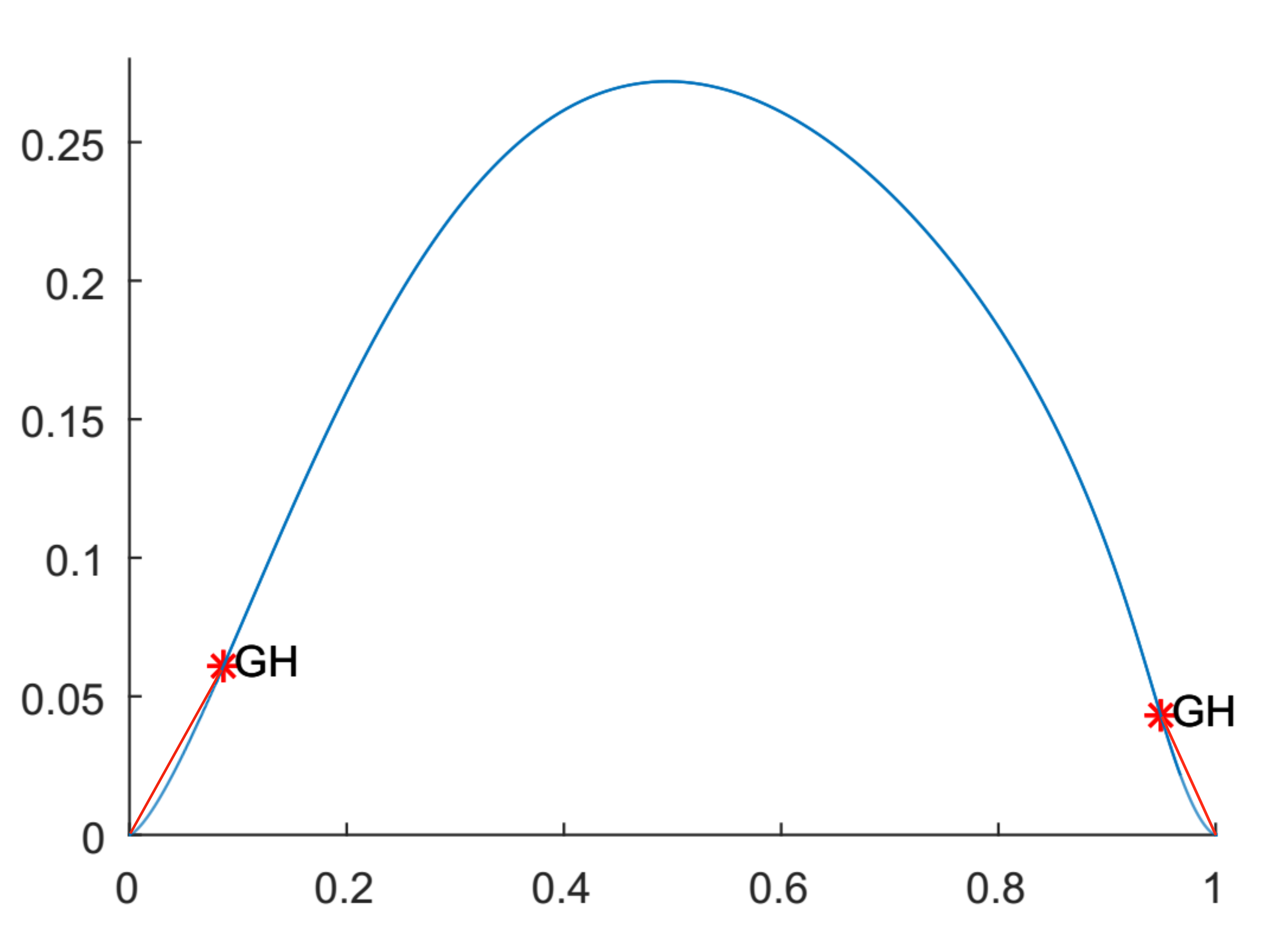}}
\node at (-4.5, 0){$\ve$};
\node at (0, -3.4){$\gamma$};
\node (u) at (0.25, 0){\textnormal{Unstable}};
\node[below=of u, yshift=1.1cm]{\textnormal{(oscillatory motion)}};
\node at (-2.4, 1.55){\textnormal{Stable}};
\node (s) at (3, 1.5){\textnormal{Stable}};
\node[below=of s, xshift=0.7cm, yshift=1.15cm]{\textnormal{(steady state)}};
\end{tikzpicture}
\caption{Two-parameter bifurcation analysis of \eqref{syseps} with respect to $(\ve, \gamma)$}
\label{bif.2para}
\end{figure}

This section is devoted to the two-parameter bifurcation analysis of \eqref{syseps}. 
In particular, we are interested in understanding the behavior of system \eqref{syseps} under the variation of parameters $(\ve, \gamma)$. 
To this end, let us represent \eqref{syseps} by
\begin{equation}\label{gen.2para}
	\dot{x}=G(x;\ve,\gamma),
\end{equation}
where $x=
\begin{bmatrix}
f & c & e
\end{bmatrix}^\top$, and $G(x;\ve,\gamma)$ denotes the right-hand side of \eqref{syseps}.
We have used the numerical continuation software \textsc{Matcont} \cite{dhooge2003matcont} to compute the two-parameter bifurcation diagram of \eqref{gen.2para} with respect to $(\ve,\gamma)$, presented in Fig. \ref{bif.2para}, where the vertical and the horizontal axes show, respectively, the behavior of $G(x;\ve,\gamma)$ with respect to $ \ve$ and $\gamma$. 
The blue curve indicates that for any $0<\gamma < 1$ and any $\ve$ below the curve, the system has unstable equilibria and hence exhibits oscillatory behavior. 
For those values of $\ve$ which are above the blue curve, the system is not oscillatory anymore, i.e. the equilibrium point is stable.
In fact, the blue curve is a curve of Hopf bifurcations where the equilibria of the system switches from being stable to unstable:
with fixed $0<\gamma<1$, as $\ve$ passes through the curve from above to below, a limit cycle is generated.

\begin{figure}[t]
\centering
\begin{tikzpicture}
	\pgftext{\includegraphics[scale=1]{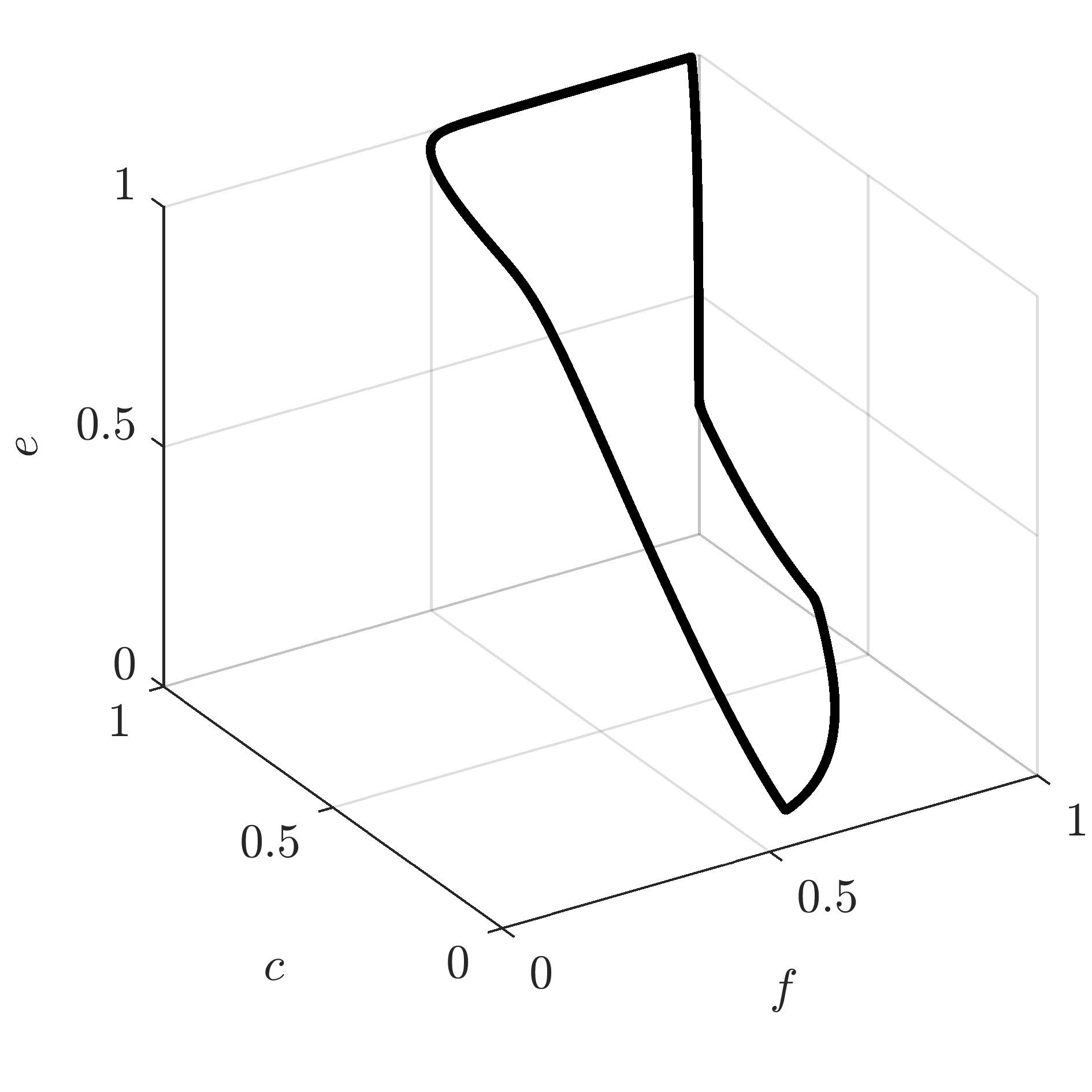}}
\end{tikzpicture}
\caption{Numerically computed attracting limit cycle of system \eqref{model} for $\ve=10^{-3}$ and $\gamma=0.9$.}
\label{reverse.cycle}
\end{figure}

As shown in Fig. \ref{bif.2para}, there are two points, denoted by ``GH'', which are generalized Hopf (or Bautin) Bifurcation points. 
At these points, the equilibria of \eqref{gen.2para} have a pair of purely imaginary eigenvalues at which the first Lyapunov exponent coefficient of the Hopf bifurcation vanishes \cite{kuznetsov2013elements}.  
Computed by $\textsc{Matcont}$, the values of $(\ve,\gamma)$ at ``GH'' points are as follows:
\begin{equation}
	(\ve_1, \gamma_1)=(0.060907128, 0.086423772), \qquad\qquad
    (\ve_2, \gamma_2)=(0.043172692, 0.949470320).
\end{equation}
In Fig. \ref{bif.2para}, the red curves are the {\e{curves of}} ``limit points” (or saddle-node bifurcation) of cycles.
For parameter values $(\ve,\gamma)$ between the blue and red curves in Fig. \ref{bif.2para}, at least two limit cycles exist simultaneously, i.e., for $\gamma$ close to 0 or $\gamma$ close to 1, with a suitable $0<\ve\ll 1$, at least one stable and one unstable limit cycle coexist.
%As it is clear from Fig. \ref{bif.2para}, these two limit cycles collide with each other and disappear at the ``GH'' points.
\begin{remark}\label{gamma.relax}
As we mentioned in Section \ref{mod.des}, due to the property of ``zero-order ultrasensitivity'', the Michaelis-Menten constants and hence $\ve$ have to be small.
Our observation from numerical simulations shows that, for sufficiently small $\ve$, system \eqref{syseps} has similar qualitative behaviors when $\gamma$ belongs to certain bounds which are close to 0 and 1. 
In this regard, we emphasize that although the position of the limit cycle changes when $\gamma$ is close to 1 (see, for instance, Fig. \ref{reverse.cycle}), the geometric analysis of the dynamics is the same as the case that $\gamma$ is close to 0, for sufficiently small $\ve$.
\end{remark}

\begin{remark}
{\e{In Section \ref{mod.des}, we have unified all the Michaelis-Menten constants of system \eqref{model} by $\ve$, resulted in system \eqref{syseps}.
Although $\gamma$ has similar size as the Michaelis-Menten constants, we have not unified it with them.
One reason is that the unit of $\gamma$ is ``$\textnormal{min}^{-1}$'', while the Michaelis-Menten constants are unitless.
The other reason is that the simultaneous limit $(\ve,\gamma)\to(0,0)$ is very singular because a certain point $(0,0,\gamma)$, playing crucial role in our analysis, approaches $(0,0,0)$ which is the intersection of three critical manifolds $f=0$, $c=0$, and $e=0$.
It would be interesting to study this limit further, which could explain the coalescence of the Hopf curve and the saddle-node curve at $(0,0)$, see Fig. \ref{bif.2para}.
Similar remark holds as $(\ve,\gamma)\to(0,1)$.}}
\end{remark}
%%%%%--------------
%}

%% file: subfiles/SlowFast.tex
Our goal is to understand the dynamics of \eqref{syseps} for small $\ve$ in the limit $\ve\to0$. However, as it is seen in \eqref{syseps}, when the variables $f, c$ and $e$ are very close to the boundary of $\Q$, the limiting behavior is different from the case that they are away from the boundary.
To resolve the aforementioned problem, one possibility is to consider an auxiliary system which is smoothly equivalent to \eqref{syseps}.
To this end, let us define 
\begin{equation*}
	H^\ve(f,c,e) := H_1^\ve(f) H_2^\ve(c) H_3^\ve(e),
\end{equation*}
where 
\begin{equation}\label{H123eps}
\begin{split}
	 & H_1^\ve(f) := (\ve + 1 - f)(\ve + 2f), \\
	 & H_2^\ve(c) := (\ve + 2 - 2c)(\ve + 2c), \\
	 & H_3^\ve(e) := (\ve + 2 - 2e)(\ve + 2e).
\end{split}
\end{equation}
Note that $H^\ve(f,c,e)>0$ for any $\ve>0$ and any $(f,c,e)\in\Q$. Therefore, we can reparametrize time of system \eqref{syseps} by multiplying both sides of \eqref{syseps} in $H^\ve(f,c,e)$, which leads to the following dynamical system
\begin{align}\label{dot}
	&\frac{df}{d\tau} = \left(\frac{\gamma(1-f)}{\ve + (1-f)} - \frac{2fe}{\ve+2f}\right)H^\ve(f,c,e),\nonumber\\
	&\frac{dc}{d\tau} = \left(\frac{8(1-c)f}{\ve + 2(1-c)} - \frac{4c}{\ve+2c}\right)H^\ve(f,c,e),\\
	&\frac{de}{d\tau} = \left(\frac{8(1-e)c}{\ve + 2(1-e)} - \frac{4e}{\ve+2e}\right)H^\ve(f,c,e),\nonumber	
\end{align}%
where, for simplicity, we recycle $\tau$ to denote the reparametrized time. 
One can rewrite \eqref{dot} as follows
\begin{align}\label{repara}
X_\ve:
\begin{cases}
	\frac{df}{d\tau} = \left[\gamma(1 - f)(\ve + 2 f) - 2fe(\ve + 1 - f)\right] H_2^\ve(c) H_3^\ve(e),\\
	\frac{dc}{d\tau} = \left[8(1 - c)f(\ve + 2c) - 4c(\ve + 2 - 2c)\right] H_1^\ve(f) H_3^\ve(e),\\
	\frac{de}{d\tau} = \left[8(1-e)c(\ve + 2e)- 4e(\ve + 2 - 2e)\right] H_1^\ve(f) H_2^\ve(c).	
\end{cases}
\end{align}
%%%%%
The vector field \eqref{repara} is smoothly equivalent to \eqref{syseps} for $\ve>0~$\cite{chicone2006}, which from now on is the object of study.
The main reason to rewrite system \eqref{syseps} into the form of system \eqref{repara} is that the latter is a singularly perturbed ODE which allows us to analyze the system using geometric methods.
Moreover, note that in contrast to \eqref{syseps}, system \eqref{repara} is polynomial, which is another of its  advantages. 
%-----------------------------------------------------------
\subsubsection{Slow-fast systems}
A Slow-Fast System (SFS) is a singularly perturbed ordinary differential equation with two timescales often presented as
\begin{equation}\label{eq:sfs1}
\begin{split}\dot x &= F(x,y,\ve),\\
\ve\dot y &= G(x,y,\ve),
\end{split}
\end{equation}
where the ``dot'' denotes derivative with respect to the \emph{slow} time $t$, $F$ and $G$ are assumed to be smooth, $x\in\R^{n_s}$, $y\in\R^{n_f}$, and $0<\ve\ll 1$ is a small parameter that describes the timescale separation between $x$ and $y$. 
The SFS presented in \eqref{eq:sfs1}, where the timescale separation is explicitly given, is said to be in the \emph{standard form}. To study standard SFSs we usually define a new \emph{fast} time $\tau=\frac{t}{\ve}$ with which system \eqref{eq:sfs1} can be rewritten as
\begin{equation}\label{eq:sfs2}
\begin{split}
x' &= \ve F(x,y,\ve),\\
y' &= G(x,y,\ve),
\end{split}
\end{equation}%
where now the ``prime'' denotes $\frac{d}{d\tau}$.  Since $\ve$ is a small parameter, we would like to draw conclusions on the overall behavior of the trajectories of a SFS from limiting systems obtained by taking the limit $\ve\to0$. In such a limit \eqref{eq:sfs1} becomes a Differential Algebraic Equation\footnote{Also known as Constrained Differential Equation \cite{takens1976constrained}.}(DAE) of the form
\begin{equation}\label{eq:cde}
\begin{split}
\dot x &= F(x,y,0),\\
0 &= G(x,y,0),
\end{split}
\end{equation}
{\e{which is called the \emph{reduced problem}.}}
On the other hand \eqref{eq:sfs2} becomes {\e{the \emph{layer problem}}}
\begin{equation}\label{eq:layer}
\begin{split}
x' &= 0,\\
y' &= G(x,y,0).
\end{split}
\end{equation}%
\begin{remark}
The SFS \eqref{eq:sfs2} is a \emph{regular} $\ve$-perturbation problem. Therefore, its solutions can be expressed as $O(\ve)$ perturbations of solutions of \eqref{eq:layer}. However, the previous is valid \emph{only} for time $\tau$ of order $O(1)$, or equivalently for time $t$ of order $O(\ve)$ in \eqref{eq:sfs1}. To describe trajectories for longer time, techniques outside the scope of regular perturbation theory are needed.
\end{remark}

Systems \eqref{eq:cde} and \eqref{eq:layer} are not equivalent. However, the critical manifold provides a relationship between the two.
\begin{definition}[Critical manifold] The critical manifold is defined as
\begin{equation}
\hjk{\C_0}=\left\{ (x,y)\in\R^{n_s}\times\R^{n_f}\, | \, G(x,y,0)=0\right\}.
\end{equation}
\end{definition}

Note that the critical manifold $\hjk{\C_0}$ serves as the phase space of the DAE and as the set of equilibrium points of the layer problem. An important property that a critical manifold may posses is normal hyperbolicity.

\begin{definition}[Normal Hyperbolicity] Consider a SFS \eqref{eq:sfs1} and its associated critical manifold $\hjk{\C_0}$. A point $\hjk{p}\in \hjk{\C_0}$ is said to be hyperbolic if the matrix \hjk{$\textnormal{D}_yG(p)$, where $D_y$ denotes the total derivative with respect to $y$,} has all its eigenvalues with non-zero real part. The critical manifold $\hjk{\C_0}$ is said to be normally hyperbolic (NH) if every point $\hjk{p}\in \hjk{\C_0}$ is hyperbolic.
\end{definition}

Fenichel theory \cite{kuehn2015multiple} describes the dynamics of a SFS with a normally hyperbolic critical manifold.

\begin{theorem}[Fenichel] Let $\hjk{\S_0\subseteq\C_0}$ be a compact {\e{and normally hyperbolic critical}} manifold of an SFS. Then, for $\ve>0$ sufficiently small, the followings hold:
\begin{itemize}
	\item There exists a locally invariant manifold $\hjk{\S_\ve}$ which is diffeomorphic to $\hjk{\S_0}$ and lies within distance of order $O(\ve)$ from $\hjk{\S_0}$.
	\item {\e{The vector field $X_\ve$ restricted to $\hjk{\S_\ve}$ is a smooth perturbation of the reduced problem.}}
	\item $\hjk{\S_\ve}$ has the same stability properties as $\hjk{\S_0}$.
\end{itemize}
\end{theorem}

In words, Fenichel theory says that if a SFS has a compact and normally hyperbolic critical manifold $\hjk{\S_0}$, the dynamics of the slow-fast systems can be inferred from the reduced flow along $\hjk{\S_0}$ and the flow of the layer equation, which provide the stability properties of $\hjk{\S_0}$.

{\e{Often slow-fast systems have critical manifolds which lose normal hyperbolicity at certain points.}}
In fact, like the system studied in this article, many interesting phenomena in several timescales such as relaxation oscillations and canards, are associated to the loss of normal hyperbolicity \cite{kuehn2015multiscale,kosiuk2016geometric,brons2015mixed,roberts2016mixed}. 

\begin{remark}
In a more general context, SFSs do not have to be given in standard form as in \eqref{eq:sfs2}. That is, SFSs can be defined by and ODE of the form $z' = H(z,\ve)$. In such a case the corresponding critical manifold $\hjk{\S_0}$ is defined by $\hjk{\S_0}=\left\{ z\in\R^{n_s+n_f} \,|\, H(z,0)=0\right\}$, while the layer equation reads as $z'=H(z,0)$. Under normal hyperbolicity of the critical manifold, all the Fenichel theory results hold for the aforementioned general case \cite{kuehn2015multiple}. When the critical manifold has non-hyperbolic points, a careful combination of Fenichel theory and the blow-up method can be employed for a detailed analysis of the dynamics of the SFS.  In the following subsection we briefly describe the blow-up method. We later show that \eqref{repara} is indeed a general SFS, and provide a detailed geometric analysis of \eqref{repara} by means of the blow-up method and Fenichel theory.
\end{remark}

\subsubsection{The blow-up method}\label{sec:blowup}

The blow-up method was introduced to describe the dynamics of SFSs near non-hyperbolic points, and is the main mathematical technique used in forthcoming section of this article. Here we just provide a brief description of the method, for more details the interested reader is referred to \cite{dumortier1996canard,krupa2001extending,kuehn2015multiple,jardon2019survey}.

First of all, note that a SFS written in the fast-time scale is an $\ve$-parameter family of vector fields. Thus, it is convenient to lift such family up and instead consider a single vector field of the form
\begin{align}\label{GeneralSFS}
X:
\begin{cases}
x' = \ve F(x,y,\ve),\\
y' = G(x,y,\ve),\\
\ve'= 0.
\end{cases}
\end{align}

\begin{definition} Consider a generalized polar coordinate transformation
\begin{eqnarray}\label{eq:polar}
	&& \Phi:\s^{n_s+n_f}\times I\to\R^{{n_s+n_f}+1}\nonumber\\
	&& \Phi(\bx,\by,\be,\br)\mapsto(\br^{\alpha}\bx,\br^{\beta}\by,\br^\gamma\be)=(x,y,\ve),
\end{eqnarray}%
where $ \sum_{i=1}^{n_s} \bx_i^2 +  \sum_{j=1}^{n_f} \by_j^2 + \be^2=1$ and $\br\in I$ where $I$ is a (possibly infinite) interval containing $0\in\R$. The corresponding \emph{(quasi-homogeneous)\footnote{A homogeneous blow-up (or simply blow-up) refers to all the exponents $\alpha$, $\beta$, $\gamma$ set to $1$.  } blow-up } is defined by $(\bx,\by,\be,\br)=\Phi^{-1}(x,y,\ve)$. 
The map $\Phi$ is called \emph{blow down}\footnote{Note that the blow-up maps the origin $0\in\R^{{n_s+n_f}+1}$ to the sphere $\s^{n_s+n_f}\times\{0\}$ while the blow down does the opposite, hence the names.}.
\end{definition}

For the purposes of this article, it is sufficient to let $\br\in[0,\rho)$, with $\rho>0$. The main idea of the blow-up method is to construct a new, but equivalent, vector field to $X$, which is defined in a higher dimensional manifold, but whose singularities are simpler compared to those of $X$.

\begin{definition} The blown up vector field $\bar X$ is induced by the blow-up map as $\bar X = \textnormal{D}\Phi^{-1}\circ X\circ\Phi$, where $\textnormal{D}\Phi$ denotes the derivative of $\Phi$. 
{\e{If $\tilde X$ vanishes on $\s^{n_s+n_f}\times\{0\}$ with order $m\in\mathbb{N}$, we define \hjk{\emph{the desingularized vector field}} $\tilde X=\frac{1}{\br^m}\bar X$.}}
%If $\tilde X$ is degenerate at $\s^{n_s+n_f}\times\{0\}$, then we define the desingularized vector field as $\bar X=\frac{1}{\br^m}\tilde X$, where $m\in\mathbb N$ is as large as possible such that $\bar X$ is well-defined and non-degenerate at $\s^{n_s+n_f}\times\{0\}$.
\end{definition}

Note that the vector fields $\bar X$ and $\tilde X$ are equivalent on $\s^{n_s+n_f}\times \left\{ \br>0\right\}$. Moreover, if the weights $(\alpha,\beta,\gamma)$ are well chosen, the singularities of $\tilde X|_{\br=0}$ are partially hyperbolic or even hyperbolic, making the analysis of $\tilde X$ simpler than that of $X$. Due to the equivalence between $X$ and $\tilde X$, one obtains all the local information of $X$ around $0\in\R^{{n_s+n_f}+1}$ from the analysis of $\tilde X$ around $\s^{n_s+n_f}\times \left\{\br\geq 0\right\}$.

While doing computations, it is more convenient to study the vector field $\tilde X$ in charts. A chart is a parametrization of a hemisphere of $\s^{n_s+n_f}\times I$ and is obtained by setting one of the coordinates $(\bx,\by,\be)\in\s^{n_s+n_f}$ to $\pm 1$ in the definition of $\Phi$. For example, one of the most important charts in the blow-up method is the central chart defined by $K_{\bar\ve}=\left\{ \be=1 \right\}$. After we study the dynamics in the relevant charts, we connect the flow together via \emph{transition maps}, allowing us a complete description of the flow of $\tilde X$ near $\s^{n_s+n_f}\times\{0\}$. In turn,  and as mentioned above, the flow of $\tilde X$ is equivalent to the flow of $X$ for $\ve>0$ sufficiently small. For more details see Section \ref{blowup} and \cite{kuehn2015multiple}. 
{\e{\begin{remark}
It is also possible to blow-up \emph{only some} of the variables in the system \eqref{GeneralSFS}, and keep the others unchanged.
In this paper, we blow-up a non-hyperbolic line of equilibria to a cylinder, see Section \ref{blowup}.
\end{remark}}}

%% file: subfiles/GSPA1.tex
\begin{figure}[t]
\centering
\begin{tikzpicture}
	\pgftext{\includegraphics[scale=1]{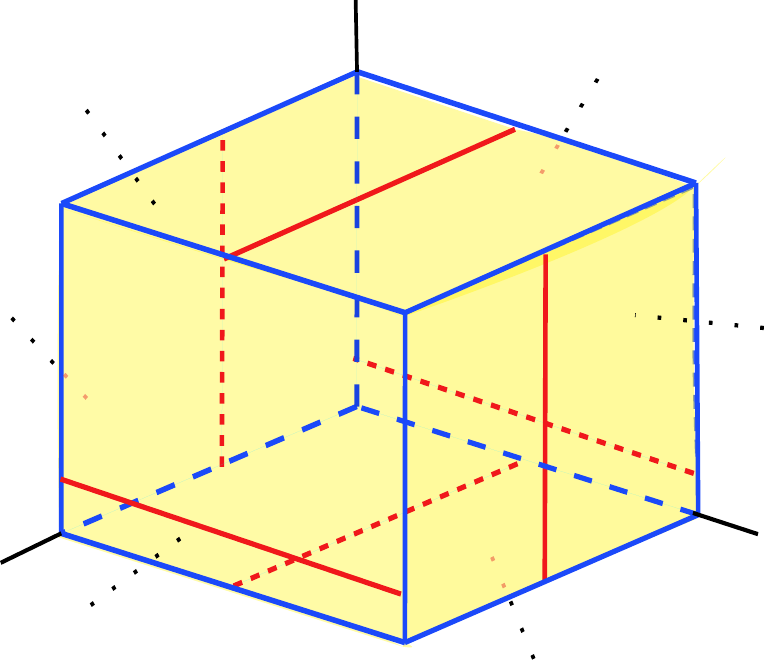}}
     %%%% S_i
            \node at (2.4,2.8){$\hjk{\S_{0,1}}$}; 
            \node at (-4.1,0.4){$\hjk{\S_{0,2}}$};
            \node at (1.5,-3.65){$\hjk{\S_{0,3}}$};
            \node at (-3.3,-3.1){$\hjk{\S_{0,4}}$};
            \node at (4.2, 0.15){$\hjk{\S_{0,5}}$};
            \node at (-3.1,2.55){$\hjk{\S_{0,6}}$};
            %%%% x
    		\node at (-4.3,-2.5){$f$};
    		\node at (4.15,-2.2){$c$};
    		\node at (-0.3,3.7){$e$};
            %%%%%% \ell lines
            \node at (-1.9, -0.1){$\ell_f$};
            \node at (0.8, -0.4){$\ell_e$}; 
            \node at (-0.5, -1.9){$\ell_c$};
            \node at (2, 0.3){$\ell^f$};
            \node at (-0.8, 1.3){$\ell^c$};
            \node at (-1.15, -1.95){$\ell^e$};
            \node at (0.02, 0.85){$\ell_1$};
            \node at (-1, -0.8){$\ell_2$};
            %%%%%% equilibrium
            \node at (-1.65, 2.3){$f^*$};
            \node at (-1.6, -2.8){$c^*$};
            \node at (3.45, -1.4){$e^*$};
\end{tikzpicture}
\caption{The critical manifold $\hjk{\S_{0}}=\bigcup_{i=1}^6 \hjk{\S_{0,i}}$, non-hyperbolic lines $\ell_f, \ell_c, \ell_e, \ell^f, \ell^c, \ell^e$ in red, all 12 non-hyperbolic edges in blue, and in particular, \hjk{ the two non-hyperbolic edges $\ell_1$ and $\ell_2$ shall play an important role in our analysis.}}
\label{unit.cube}
\end{figure}
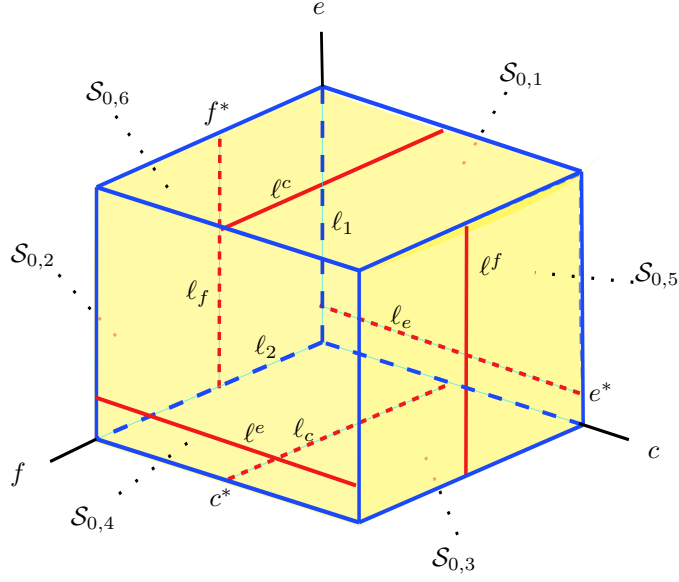

Setting $\ve=0$ in \eqref{repara} results in the layer problem
\begin{equation}\label{layer}
	\begin{split}
		\frac{df}{d\tau} &= \left(\gamma - e\right)H^0(f,c,e),\\
		\frac{dc}{d\tau} &= 2\left(2f - 1\right)H^0(f,c,e),\\
		\frac{de}{d\tau} &= 2(2c - 1)H^0(f,c,e),		
	\end{split}
\end{equation}
with 
\begin{equation*}
 	H^0(f,c,e) = 32fce(1-f)(1-c)(1-e).
\end{equation*}
%%%
Apart form the isolated equilibrium point $P:=(0.5, 0.5, \gamma)$, which is inside the cube $\hjk{\Q}$, the boundary of $\hjk{\Q}$, which consists of six planes, is the equilibria set of the layer problem \eqref{layer}.
We denote each plane of equilibria by $\hjk{\S_{0,i}}$ $(i=1,2,...,6)$ as follows:
\begin{equation}\label{Si}
\begin{split}
\hjk{\S_{0,1}}&:=\left\{ (f,c,e)\in\R^3 \, | \, f=0, \, c\in[0,1],\, e\in[0,1] \right\},\\
\hjk{\S_{0,2}}&:=\left\{ (f,c,e)\in\R^3 \, | \, f\in[0,1], \, c=0,\, e\in[0,1] \right\},\\
\hjk{\S_{0,3}}&:=\left\{ (f,c,e)\in\R^3 \, | \, f\in[0,1], \, c\in[0,1],\, e=0 \right\},\\
\hjk{\S_{0,4}}&:=\left\{ (f,c,e)\in\R^3 \, | \, f=1, \, c\in[0,1],\, e\in[0,1] \right\},\\
\hjk{\S_{0,5}}&:=\left\{ (f,c,e)\in\R^3 \, | \, f\in[0,1], \, c=1,\, e\in[0,1] \right\},\\
\hjk{\S_{0,6}}&:=\left\{ (f,c,e)\in\R^3 \, | \, f\in[0,1], \, c\in[0,1],\, e=1 \right\}.
\end{split}
\end{equation}
Therefore $\hjk{\S_{0}}:=\bigcup_{i=1}^6 \hjk{\S_{0,i}}$ is the critical manifold.
The stability of system \eqref{repara} changes at lines 
$\ell_f\in \hjk{\S_{0,2}} ,\, \ell^f\in \hjk{\S_{0,5}}$ (given by $f=f^*$); 
$\ell_c\in \hjk{\S_{0,3}} ,\, \ell^c\in \hjk{\S_{0,6}}$ (given by $c=c^*$); and
$\ell_e\in \hjk{\S_{0,1}} ,\, \ell^e\in \hjk{\S_{0,4}}$ (given by $e=e^*$).
Moreover, {\e{the 12 edges of the unit cube, where the 6 planes $\hjk{\S_{0,i}}$ intersect, are non-hyperbolic lines as well.}}
However, for our analysis, only the lines $\ell_1=\hjk{\S_{0,1}}\cap \hjk{\S_{0,2}}$ and $\ell_2=\hjk{\S_{0,2}}\cap \hjk{\S_{0,3}}$ are crucial (see Figure \ref{unit.cube}).
%%%---------------------
The stability of points in $\hjk{\S_{0}}$ is summarized in the following lemma.
\begin{lemma}
The critical manifold $\hjk{\S_{0}}$ of the layer problem  \eqref{layer} has the following properties:
\begin{itemize}
\item $\hjk{\S_{0,1}}$ is attracting for $e>e^*$ and repelling for $e<e^*$.  
\item $\hjk{\S_{0,2}}$ is attracting for $f<f^*$ and repelling for $f>f^*$.  
\item $\hjk{\S_{0,3}}$ is attracting for $c<c^*$ and repelling for $c>c^*$.  
\item $\hjk{\S_{0,4}}$ is attracting for $e<e^*$ and repelling for $e>e^*$.  
\item $\hjk{\S_{0,5}}$ is attracting for $f>f^*$ and repelling for $f<f^*$. 
\item $\hjk{\S_{0,6}}$ is attracting for $c>c^*$ and repelling for $c<c^*$. 
\item The equilibrium $P:=(0.5, 0.5, \gamma)$ is a saddle-focus point.
\item The lines $\ell_f\in \hjk{\S_{0,1}}, \ell_c\in \hjk{\S_{0,2}}, \ell_e\in \hjk{\S_{0,3}}, \ell^f\in \hjk{\S_{0,4}}, \ell^c\in \hjk{\S_{0,5}}, \ell^e\in \hjk{\S_{0,6}}$, {\e{all 12 edges of the unit cube, and in particular,}} the edges $\ell_1=\hjk{\S_{0,1}}\cap \hjk{\S_{0,2}}$ and $\ell_2=\hjk{\S_{0,2}}\cap \hjk{\S_{0,3}}$ are non-hyperbolic.
\end{itemize}
\end{lemma}
\begin{proof}
The eigenvalues of the linearization of system \eqref{layer} at points, e.g., in the plane $\hjk{\S_{0,1}}$ are given by 
\begin{equation*}
\lambda_1=\lambda_2=0, \qquad \lambda_3 = -32ce(c-1)(e-1)(e-\gamma).
\end{equation*}
It is clear that $\lambda_3$ is zero at the boundary of $\hjk{\S_{0,1}}$, and also along the line $l_e$ given by $e=e^*$.
Therefore, $\hjk{\S_{0,1}}$ is attracting for $e>e^*$ and it is repelling for $e<e^*$.
% The attracting (resp. repelling) part of $S_1$ is denoted by $S^1_a$ (resp. $S^1_r$), and defined by 
% \begin{equation*}
% \begin{split}
% S^1_a &=\left\{ (c,e) \in\mathbb{R}^2 \, | \, (c,e)\in\left[\delta,\ 1-\delta\right]\times\left[\gamma + \delta, \ 1-\delta\right] \right\}, \\
% S^1_r &=\left\{ (c,e) \in\mathbb{R}^2 \, | \, (c,e)\in\left[\delta,\ 1-\delta\right]\times\left[\delta, \gamma-\delta\right] \right\}.
% \end{split}
% \end{equation*}
The proof of the other cases is performed analogously. 
\end{proof}
%%%%-------------- 
%%%%%--------------------
We denote the interior of the cube $\hjk{\Q}$ by $\hjk{\mathring\Q}$.
Note that when $(f,c,e)\in\hjk{\mathring\Q}$, the layer problem \eqref{layer} can be divided by the positive term $H^0(f,c,e) = 32fce(1-f)(1-c)(1-e)$.
Therefore away from the critical manifold $S$, all the variables evolve on the fast time scale $\tau$ and the orbits of the layer problem \eqref{layer} are identical to the orbits of the linear system 
\begin{equation}\label{identical}
	\begin{split}
		\frac{df}{d\tau} &= \gamma - e,\\
		\frac{dc}{d\tau} &= 2(2f - 1),\\
		\frac{de}{d\tau} &= 2(2c - 1).	
	\end{split}
\end{equation}
\begin{remark}
System \eqref{identical} is the limit of \eqref{syseps} when $\ve\to 0$ and $(f,c,e)\in\hjk{\mathring\Q}$. 
\end{remark}

%% file: subfiles/ReducedProblem.tex
From Subsection \ref{critical.manifold}, we know that \hjk{the boundary of $\Q$} is the critical manifold $\hjk{\S_0}$.
Any compact subset of $\hjk{\S_0}$ that does not contain any non-hyperbolic point is normally hyperbolic, and hence Fenichel theory \cite{fenichel1979geometric} is applicable.
In other words, this theory implies that the normally hyperbolic parts of $\hjk{\S_0}$ perturb to slow manifolds, which lie within a distance of order $O(\ve)$ of the critical manifold $\hjk{\S_0}$.
In the following, we compute the slow manifolds and analyze the reduced flows in the planes $\hjk{\S_{0,1}},\, \hjk{\S_{0,2}},\, \hjk{\S_{0,3}}$ and $\hjk{\S_{0,6}}$ which are essential for our analysis.
%%%%%%%%%%%%%%%%%%%%%%%%%%%%% Analysis of S^1 %%%%%%%%%%%%%%%%%%%%%%%%%%%%%%%%%%%%
\begin{lemma}\label{slow1}
For sufficiently small $\delta>0$, there exist $\ve_0>0$ and a smooth function $\hjk{h_{\ve,1}}(c,e)$ defined on $\hjk{I_1^\txta}=\left[\delta,\ 1-\delta\right]\times\left[\gamma + \delta,\ 1-\delta\right]$ such that the manifold 
\begin{equation}
\hjk{\S_{\ve,1}^\txta}=\left\{ (f,c,e)\in\hjk{\Q}\, | \, f=\hjk{h_{\ve,1}}(c,e), \, (c,e)\in \hjk{I_1^\txta} \right\},
\end{equation}
is a locally invariant attracting manifold of \eqref{repara} for $\ve\in(0,\ve_0]$. 
The function $\hjk{h_{\ve,1}}(c,e)$ has the expansion
\begin{equation}
	\hjk{h_{\ve,1}}(c,e) = \frac{\gamma}{2(e -\gamma)}\ve + O(\ve^2).
\end{equation}
\end{lemma}
\begin{proof}
Since the set $\hjk{I_1^\txta}$ is hyperbolic, Fenichel theory implies that there exists a sufficiently small $\ve_0>0$ such that the function $\hjk{h_{\ve,1}}(c,e)$ has the expansion $\hjk{h_{\ve,1}}(c,e) = \eta(c,e)\ve+O(\ve^2)$ for all $\ve\in(0,\ve_0]$. 
Due to invariance, we can substitute $\hjk{h_{\ve,1}}(c,e)$ into the equation of $\frac{df}{d\tau}$ in \eqref{repara} and identify coefficients of $\ve$. 
By doing so, we obtain
\begin{equation}\label{eta2}
\eta(c,e)= \frac{\gamma}{2(e -\gamma)}.
\end{equation}
Note that \eqref{eta2} reflects the fact that the manifold $\hjk{\S_{\ve,1}^\txta}$ is not well-defined when $e=\gamma$. 
Thus, the invariant manifold $\hjk{\S_{\ve,1}^\txta}$ is given as stated in the lemma, which completes the proof.
\end{proof}

\hjk{
For the sake of brevity, we summarize the analysis in the planes $\hjk{\S_{0,2}},\, \hjk{\S_{0,3}}$ and $\hjk{\S_{0,6}}$ in Table \ref{tab:s}, which is shown by following the same line of reasoning as the one of Lemma \ref{slow1}.
For more details, the interested reader is referred to~\cite{taghvafard2018modeling}.
\renewcommand{\arraystretch}{1.5}
\begin{table}[htbp]\centering
    \begin{tabular}{l|l|l}
     $i$ & \multicolumn{1}{c|}{$I_i^\txta$}  & \multicolumn{1}{c}{$\S_{\ve,i}^\txta$}   \\
     \hline
     $2$ & $(f,e)\in\left[\delta, \ \frac{1}{2}-\delta\right]\times\left[\delta,\ 1-\delta\right]$ & $c=\frac{f}{1-2f}\ve + O(\ve^2)$\\
     $3$ & $(f,c)\in\left[\delta, \ 1-\delta\right]\times\left[\delta,\ \frac{1}{2}-\delta\right]$& $e=\frac{c}{1-2 c}\ve + O(\ve^2)$\\
     $6$ & $(f,c)\in\left[\delta, \ 1-\delta\right]\times\left[\frac{1}{2}+\delta, 1-\delta\right]$ & $e=1+\frac{1}{2(1-2 c)}\ve + O(\ve^2)$
\end{tabular}
\caption{For each row $i$ we show the interval of definition of $\S_{\ve,i}^\txta$ and the relation by which it is defined, all analogous to Lemma~\ref{slow1}.}
\label{tab:s}
\end{table}
\renewcommand{\arraystretch}{1.0}
}
\begin{remark}
{\e{Similar results can be obtained for the ``repelling'' parts $\hjk{\S_{\ve,i}^\txtr}$, $i=1,2,...,6$.
However, these are not needed in our analysis.}} \hjk{Nonetheless, we point out that the slow manifolds $\S_{\ve,i}^\txtr$ would be expressed by the same functions $h_{\ve,i}$ and appropriate intervals $I_i^\txtr$.} %Furtheremore we shall define $\S_{\ve,i}$ by the union $\S_{\ve,i}=\S_{\ve,i}^\txta\cup\S_{\ve,i}^\txtr$.}
\end{remark}
\begin{remark}
{\e{The expansions of the functions $\hjk{h_{\ve,i}}(\cdot, \cdot)$, $i=1, 2, 3, 6$, also explain why it is necessary to restrict the domain of definition of the slow manifolds to $I_i^\txta$ to exclude their singularities.}}
\end{remark}
%%%%%%%%%%---------------------------------
We now turn to the analysis of the reduced flows in the planes $\hjk{\S_{0,1}},\, \hjk{\S_{0,2}},\, \hjk{\S_{0,3}}$ and $\hjk{\S_{0,6}}$ which, respectively, means the planes $f=0, c=0, e=0$ and $e=1$.
We know that system \eqref{repara} has the fast time scale $\tau$.
By substituting the functions $\hjk{h_{\ve,i}}$, $i=1, 2, 3, 6$ into \eqref{repara}, transforming the fast time variable to the slow one by $t=\ve\tau$, and setting $\ve=0$, the equations governing the slow dynamics on the critical manifold $\hjk{\S_{0,i}}$ are computed.
In the following, we give the analysis in the plane $\hjk{\S_{0,1}}$.

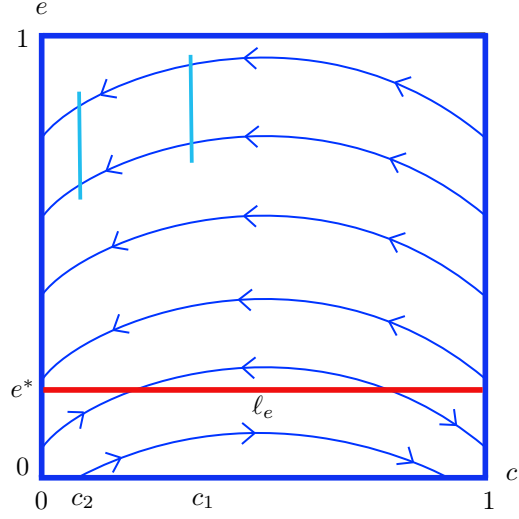
\begin{figure}[htbp]
\centering
\begin{tikzpicture}
\pgftext{\includegraphics[scale=1]{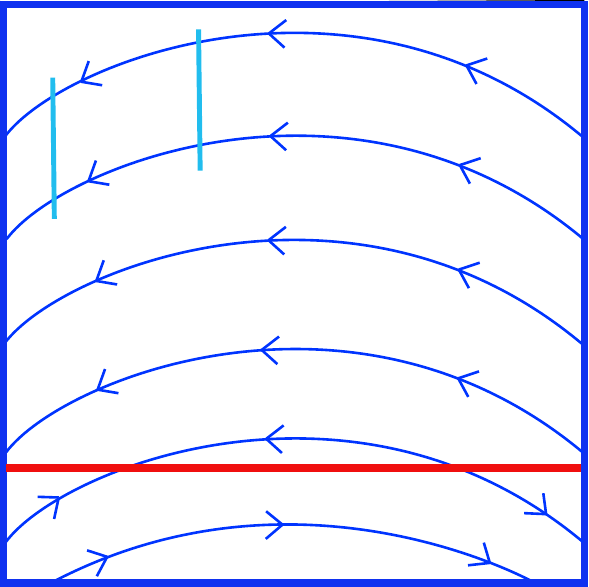}}
\node at (3.3, -2.9){$c$};
\node at (-2.95, -3.25){0};
\node at (-2.4, -3.25){$c_2$};
\node at (-0.8, -3.25){$c_1$};
%\node at (-0.2, -3.55){$c$};
\node at (3, -3.25){1};
%%%%%
\node at (-2.95, 3.3){$e$};
\node at (-3.2, -1.75){$e^*$};
\node at (0, -2){$\ell_e$};
\node at (-3.2, -2.8){0};
\node at (-3.2, 2.9){1};
\end{tikzpicture}
\caption{Flow of the slow vector field in $\hjk{\S_{0,1}}$, non-hyperbolic line $\ell_e$ in red, and sections $c_1, c_2$ in cyan.}
\label{redVF1}
\end{figure}

After substituting $\hjk{h_{\ve,1}}$ into system \eqref{repara}, the dynamics of the reduced system in $\hjk{\S_{0,1}}$, i.e., on the plane $f=0$, is governed by 
\begin{equation}\label{Red.S11}
\begin{split}
c' & = \frac{-32 ce^2(c-1)(e-1)}{e-\gamma}\ve + O(\ve^2),\\
e' & =  \frac{32 ce^2(c-1)(e-1)(2c-1)}{e-\gamma}\ve + O(\ve^2),
\end{split}
\end{equation}
where $'$ denotes the differentiation with respect to $\tau$. 
Now by dividing out a factor of $\ve$, which corresponds to switching from the fast time variable to the slow one, we have
\begin{equation}\label{Red1_a}
\begin{split}
\dot{c} & = \frac{-32 ce^2(c-1)(e-1)}{e-\gamma} + O(\ve),\\
\dot{e} & =  \frac{32 ce^2(c-1)(e-1)(2c-1)}{e-\gamma} + O(\ve),
\end{split}
\end{equation}
where the overdot represents differentiation with respect to $t=\ve\tau$.
Now, by setting $\ve=0$ in \eqref{Red1_a}, the reduced flow on $\hjk{\S_{0,1}}$ is given by
\begin{equation}\label{Red.S12}
\begin{split}
\dot{c} & = \frac{-32 ce^2(c-1)(e-1)}{e-\gamma},\\
\dot{e} & =  \frac{32 ce^2(c-1)(e-1)(2c-1)}{e-\gamma}.
\end{split}
\end{equation}
As it is clear, the vector field \eqref{Red.S12} is singular at the line $\ell_e$, given by $e=e^*$. 
In other words, the flow \eqref{Red.S12} is not defined on the line $\ell_e$. 
{\e{The lines $c=0$, $e=0$, $c=1$, and $e=1$, shown in Figure \ref{redVF1}, are lines of equilibria.}}
The line $c=0$ is attracting for $e>e^*$ and it is repelling for $e<e^*$, while the line $c=1$ is attracting for $e<e^*$ and repelling for $e>e^*$.

By dividing out the factor $\frac{32 ce^2(c-1)(e-1)}{e-\gamma}$ in \eqref{Red.S12}, the orbits of the reduced flow can be derived from the desingularized system 
\begin{equation}\label{Red.S13}
\begin{split}
	\dot{c} & = -1,\\
	\dot{e} & =  2c-1,
\end{split}
\end{equation} 
which can be integrated explicitly.
\begin{remark}
For $e>e^*$, systems \eqref{Red.S12} and \eqref{Red.S13} have qualitatively the same dynamics when $c,e\in(0,1)$. 
In particular, the vector field \eqref{Red.S13} is $C^\infty$-equivalent but not $C^\infty$-conjugate to the vector field \eqref{Red.S12}.
For the case that $e<e^*$, the direction of the vector field \eqref{Red.S12} is not preserved in the vector field \eqref{Red.S13}. 
However, for our analysis, it suffices to study the flow of system \eqref{Red.S12} when $e>e^*$, \hjk{or equivalently on $\S_{0,1}^\txta$}.
\end{remark}
{\e{\begin{lemma}\label{contract.s1}
For $e>e^*$, the reduced flow \eqref{Red.S12} on $\hjk{\S_{0,1}}$ and hence the slow flow \eqref{Red1_a} on $\hjk{\S_{\ve,1}^\txta}$ maps section $\left\{c=c_1\right\}$ to $\left\{c=c_2\right\}$, where $0<c_2<c_1<\frac{1}{2}$;
this map is well-defined and its first derivative with respect to $e$ is equal to one.
\end{lemma}
}}
{
\hjk{
\begin{proof}
It suffices to consider \eqref{Red.S13}. Let $\Pi(e)$ denote the map from $\left\{c=c_1\right\}$ to $\left\{c=c_2\right\}$ induced by the flow of \eqref{Red.S13}. Then, it is straightforward to get $\Pi(e)=e+c_2-c_2^2-c_1+c_1^2$, from which the statement follows. 
\end{proof}
}
}
% \begin{lemma}\label{contract.s1}
% For $e>e^*$, the reduced flow \eqref{Red.S12} on $S^1$, and hence the slow flow \eqref{Red1_a} on $S_{a,\ve}^1$ are contracting in $e$, i.e., the induced map between sections $c=c_1$ and $c=c_2$ with $0<c_2<c_1<\frac{1}{2}$ contracts the variable $e$ (see Figure \ref{redVF1}).
% \end{lemma}
%%%%%%%%%%%%%%%---------------------------------------------

\begin{figure}
    \captionsetup{singlelinecheck=off}

    \begin{subfigure}[h]{0.45\linewidth}
    \centering
    \begin{tikzpicture}
    \pgftext{\includegraphics[scale=1]{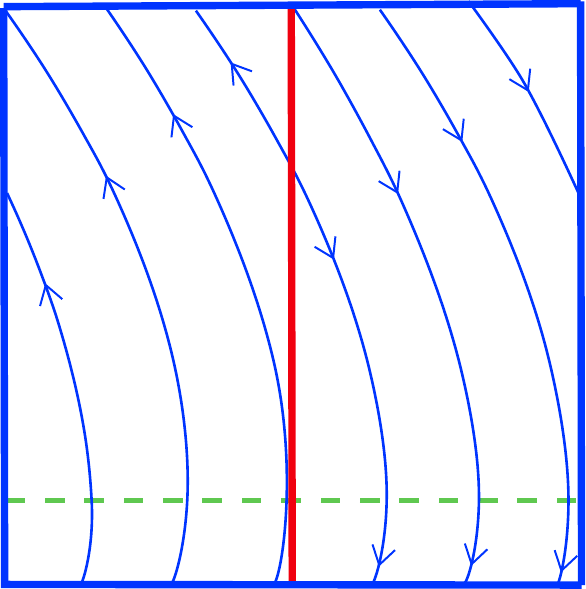}}
    \node at (-2.9, -3.2){1};
    \node at (3, -3.25){0};
    \node at (3.2, -2.83){0};
    \node at (3.15, 2.85){1};
    \node at (-3.3, -2.9){$f$};
    %%%%%
    \node at (0, -3.25){$f^*$};
    \node at (-0.25, 0.55){$\ell_f$};
    \node at (3.2, -2.15){$\gamma$};
    \node at (2.9, 3.3){$e$};
    \end{tikzpicture}
    \caption[]{The desingularized slow flow along $\hjk{\S_{0,2}}$ is
    \begin{equation}\label{Red.S23}
        \begin{split}
            &\dot f=\gamma-e\\
            &\dot e=-2.
        \end{split}
    \end{equation}}
    \label{redVF2}
    \end{subfigure}
    \hfill
    \begin{subfigure}[h]{0.45\linewidth}
    \centering
    \begin{tikzpicture}
    \pgftext{\includegraphics[scale=1]{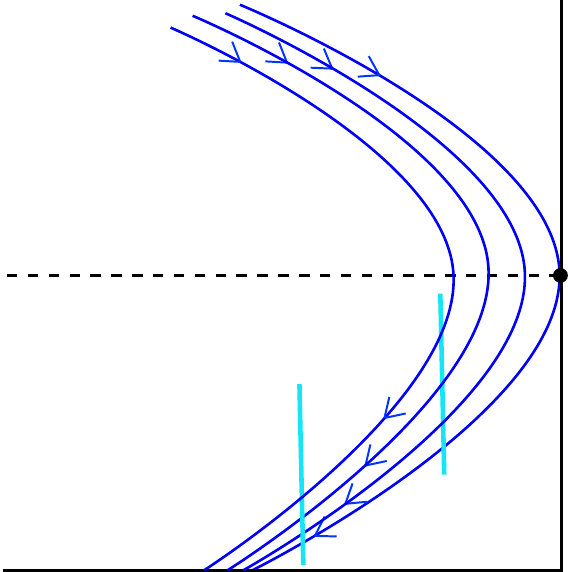}}
    \node at (2.8, 3.25){$e$};
    \node at (2.8, -3.25){0};
    \node at (-3.2, -2.9){$f$};
    \node at (3.1, 0.1){$\gamma$};
    \node at (1.7, -3.25){$f_1$};
    \node at (0.2, -3.25){$f_2$};
    \end{tikzpicture}
    \caption{The slow flow of vector field~\eqref{Red.S23} along $\hjk{\S_{0,2}}$ around $e=\gamma$ as well as the sections $f_1,f_2$ in cyan close to zero.
        Note that the variable $f$ is tangent to the line $f=0$ at $e=\gamma$.}
    \label{redVF2zoomin}
    \end{subfigure}
    %%%%%%-------------------
    \vfill
    \begin{subfigure}[h]{0.45\linewidth}
    \centering
    \begin{tikzpicture}
    \pgftext{\includegraphics[scale=1]{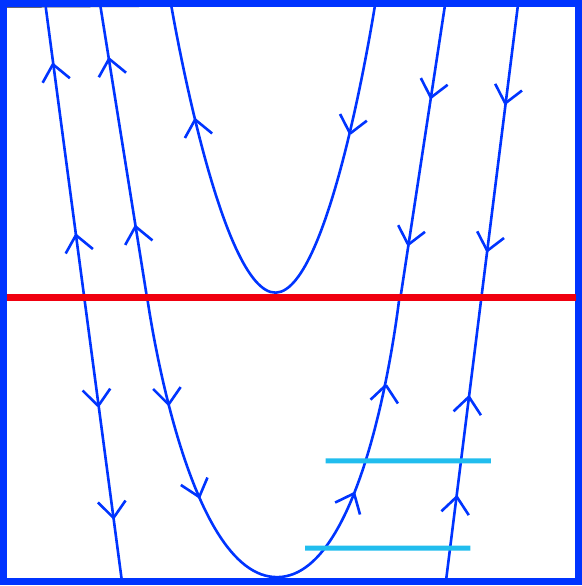}}
    \node at (-0.1, -3.25){$f^*$};
    \node at (-2.7, -3.25){0};
    \node at (2.8, -3.25){1};
    \node at (3.3, -2.9){$f$};
    %%%%%
    \node at (-3.2, 0){$c^*$};
    \node at (0.2, -0.3){$\ell_c$};
    \node at (-3.2, -1.65){$c_2$};
    \node at (-3.2, -2.6){$c_1$};
    \node at (-3.2, -3){0};
    \node at (-3.2, 2.9){1};
    \node at (-2.9, 3.3){$c$};
    \end{tikzpicture}
    \caption[]{The desingularized slow flow along $\hjk{\S_{0,3}}$ is
    \begin{equation}\label{Red.S33}
        \begin{split}
            \dot f &= \gamma, \\
            \dot c &=  2(2f-1).
        \end{split}
    \end{equation}
    For $f>f^*$, the reduced flow~\eqref{Red.S33} contracts the variable $f$ between sections $c=c_1$ and $c=c_2$ with $0<c_1<c<c_2<c^*$.}
    \end{subfigure}
    \hfill
    \begin{subfigure}[h]{0.45\linewidth}
    \centering
    \begin{tikzpicture}
    \pgftext{\includegraphics[scale=1]{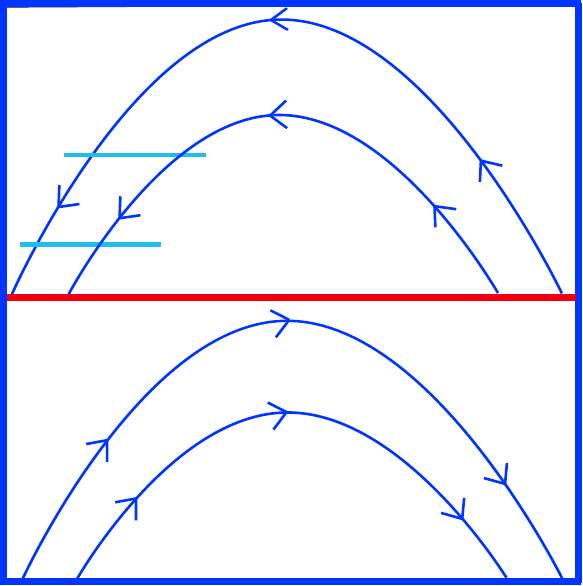}}
    \node at (-0.1, -3.25){$f^*$};
    \node at (-2.7, -3.25){0};
    \node at (3, -3.25){1};
    %%%%%
    \node at (-3.2, 0){$c^*$};
    \node at (0.5, 0.2){$\ell^c$};
    \node at (-3.2, 0.5){$c_2$};
    \node at (-3.2, 1.3){$c_1$};
    \node at (-3.2, -2.8){0};
    \node at (-3.2, 2.9){1};
    \node at (-2.9, 3.2){$c$};
    \end{tikzpicture}
    \caption[]{The desingularized slow flow along $\hjk{\S_{0,6}}$ is
    \begin{equation}\label{Red.S66}
        \begin{split}
            \dot f &= \gamma-1, \\
            \dot c &=  2f-1.
        \end{split}
    \end{equation}     
    For $f<f^*$, the reduced flow~\eqref{Red.S66} contracts the variable $f$ between sections $c=c_1$ and $c=c_2$ with $c^*<c_1<c<c_2<1$.}
    \end{subfigure}
    \caption{The reduced flows along $\hjk{\S_{0,2}^\txta},\, \hjk{\S_{0,3}^\txta}$ and $\hjk{\S_{0,6}^\txta}$.
    The vector fields are singular at red lines.
    The thick blue lines are lines of equilibria.
    The direction of vector fields shows which (part of) line of equilibria is either attracting or repelling.}
    \label{compact.figures}
\end{figure}

In order to obtain the equations governing the slow flow along $\hjk{\S_{\ve,2}^\txta},\, \hjk{\S_{\ve,3}^\txta}$ and $\hjk{\S_{\ve,6}^\txta}$, a similar analysis can be done by inserting the functions $\hjk{h_{\ve,2}},\, \hjk{h_{\ve,3}}$ and $\hjk{h_{\ve,6}}$ into \eqref{repara} and dividing out a factor of $\ve$, which corresponds to switching to the slow time scale $t = \ve\tau$.     
Next, by setting $\ve=0$ one obtains the reduced flow on the critical manifolds $\hjk{\S_{0,2}},\, \hjk{\S_{0,3}}$ and $\hjk{\S_{0,6}}$.
For the sake of brevity, we have summarized the slow flows along $\hjk{\S_{0,2}},\, \hjk{\S_{0,3}}$ and $\hjk{\S_{0,6}}$ in Fig.~\ref{compact.figures}.
For more details, the interested reader is referred to~\cite{taghvafard2018modeling}.
%For the sake of brevity, we summarize the slow flows on $\hjk{\S_{0,2}},\, \hjk{\S_{0,3}}$ and $\hjk{\S_{0,6}}$ in Lemmas \ref{red_s2}, \ref{red_s3} and \ref{red_s6}, which are crucial for our analysis.}

%%%%%%-----------------------------------------
%%%%%%%%------ compact presentation of figures by Hadi -------------

%% file: subfiles/SC.tex
%%%%%%% Singular cycle section
In this section, we present the overall behavior of the singular cycle, which is a closed curve consisting of alternating parts of the layer problem, and the critical manifold $\S_0$.
However, by the information that we have so far from the critical manifold and the layer problem, we cannot fully describe the singular cycle close to the non-hyperbolic lines $\ell_1$ and $\ell_2$.
A full description of the singular cycle for those parts that cannot be derived from the critical manifold and the layer problem is presented in Section \ref{blowup} by the blow-up method.

\begin{figure}[htbp]
\centering
\begin{tikzpicture}
\pgftext{\includegraphics[scale=1]{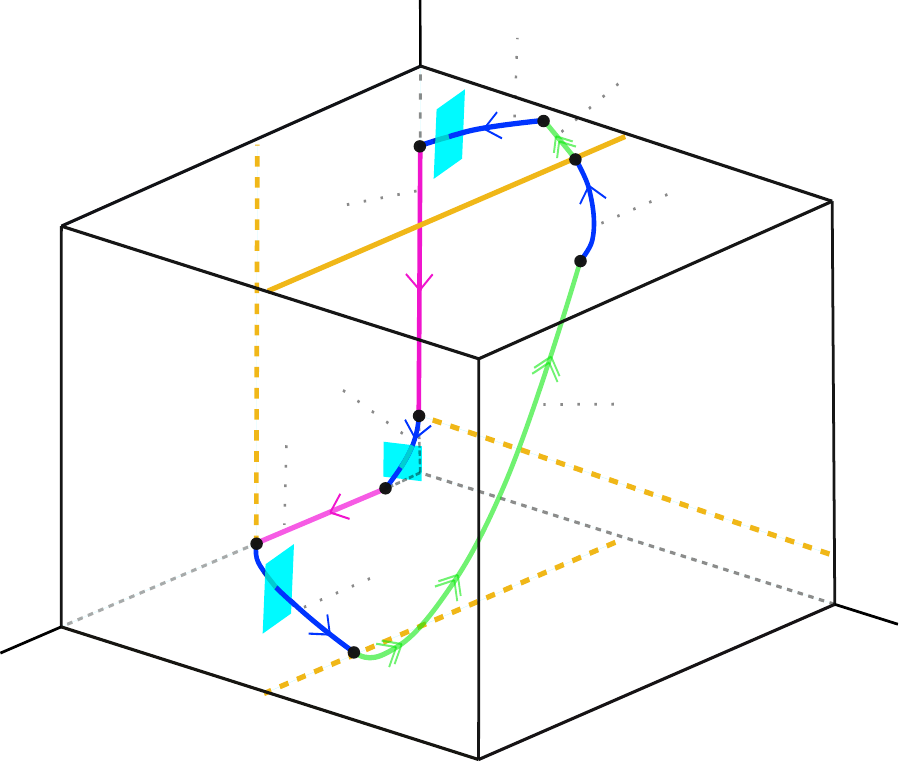}}
\node at (-0.9, -3.1){$p_1$};
\node at (0.8, 2.35){$q_e$};
\node at (1.15, 1.9){$q^e$};
\node at (-0.55, 2.3){$p^e$};
\node at (-0.05, -0.1){$p_e$};
\node at (-0.5, -1.4){$p_f$};
\node at (-2.2, -1.95){$p^f$};
\node at (1.65, 1.2){$p_2$};
\node at (0.3, -4.2){1};
\node at (3.9, -2.6){1};
%%%%%
\node at (-1.8, -3.55){$0.5$};
\node at (-4.2, -2.25){1};
\node at (-4.8, -2.85){$f$};
\node at (4.8, -2.6){$c$};
\node at (-0.3, 3.3){1};
\node at (-0.3, 4.2){$e$};
\node at (4.1, -1.8){$\gamma$};
\node at (-2.2, 2.7){$0.5$};
%%%%--
\node at (0.1, 2){$\Sigma_1$};
\node at (-0.9, -0.75){$\Sigma_2$};
\node at (-2.1, -2.6){$\Sigma_3$};
%%%%%-----
\node at (-0.45, -2){$\omega_1$};
\node at (2.1, -0.25){$\omega_2$};
\node at (2.6, 1.9){$\omega_3$};
\node at (1.9, 3.2){$\omega_4$};
\node at (0.7, 3.7){$\omega_5$};
\node at (-1.35, 1.75){$\omega_6$};
\node at (-1.05, 0.15){$\omega_7$};
\node at (-1.6, -0.3){$\omega_8$};
\end{tikzpicture}
\caption{Schematic diagram of the singular cycle $\Gamma_0$.}
\label{SinCyc}
\end{figure}
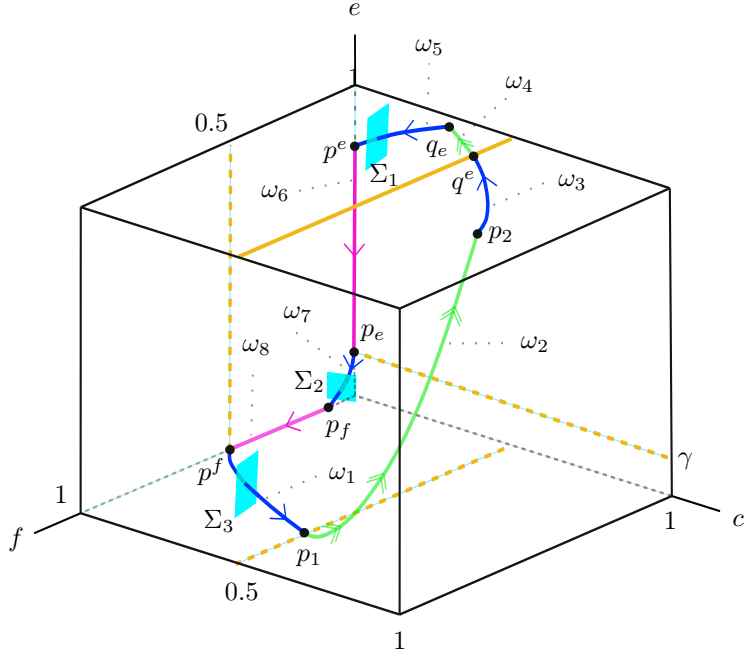

{\e{The construction of the singular cycle $\Gamma_0$ starts at the point $p^f:=(0.5, 0,0)$.
This point is connected to the point $p_1:=(\frac{1+\sqrt{\gamma}}{2},0.5,0)\in\ell_c$ through the orbit $\omega_1$ of the reduced flow \eqref{Red.S33}. 
Starting at $p_1$, the layer problem \eqref{identical} intersects the attracting part of the plane $\hjk{\S_{0,6}^\txta}$ in a point, denoted by $p_2$.
This point is connected to a point, denoted by $q^e\in\ell^c$, through the orbit $\omega_3$ of the reduced flow \eqref{Red.S66}.
Starting at $q^e$, through the layer problem \eqref{identical}, the orbit $\omega_4$ intersects the plane $\hjk{\S_{0,1}^\txta}$ at a point, denoted by $q_e$.
The orbit $\omega_5$ of the reduced flow \eqref{Red.S11} connects $q_e$ to a point, denoted by $p^e\in\ell_1$, which is the intersection of $\hjk{\S_{0,1}^\txta}$ and $\hjk{\S_{0,2}^\txta}$;
$p^e$ is connected to the point $p_e:=(0,0,\gamma)$ by a segment on the line $\ell_1$, denoted by $\omega_6$.
The orbit $\omega_7$ of the reduced flow \eqref{Red.S23} connects $p_e$ to the point $p_f:=(\frac{\gamma^2}{4}, 0, 0)$;
Finally, $p_f$ is connected to $p^f$ by a segment on the line $\ell_2$, denoted by $\omega_8$.}}
Hence, the singular cycle $\Gamma_0\in\mathbb{R}^3$ of system \eqref{repara} for $\ve=0$ is defined as follows (see Fig. \ref{SinCyc}):
\begin{equation}\label{singcycle}
\Gamma_0:=\omega_1\cup\omega_2\cup\omega_3\cup\omega_4\cup\omega_5\cup\omega_6\cup\omega_7\cup\omega_8.
\end{equation}

%%%%%
%%%%%%

{\e{
\begin{remark}\label{analytic.point}
All the orbits $\omega_j$ ($j=1,2, ..., 8$) are known analytically.
\end{remark}
Owing to the fact that the layer problem is linear, all the points that connect $\omega_j$ to $\omega_{j+1}$ are explicitly known. 
For the particular quantity $\gamma=0.08$, we have $p^f=(0.5, 0, 0)$, $p_1\approx (0.6414, 0.5, 0)$, $p_2\approx(0.3638, 0.8485, 1)$, $q^e\approx(0.0771, 0.5, 1)$, $q_e\approx(0, 0.3438, 0.9743)$, $p^e\approx(0, 0, 0.7487)$, $p_e=(0, 0, 0.08)$, and $p_f=(0.0016, 0, 0)$.

\begin{remark}
At the singular level, there is no visible flow on the segments $\omega_6$ and $\omega_8$.
The blow-up analysis, carried out in Section \ref{blowup}, will reveal a hidden flow for such segments.
\end{remark}}}

%% file: subfiles/MainResult.tex
In view of the singular cycle $\Gamma_0$, introduced in the previous subsection, we are now ready to present the main result. 
\begin{theorem}\label{main.result}
Assume that $\Gamma_0$ is the singular cycle described in Section \ref{singularcycle}.
Then for sufficiently small $\ve>0$, there exists a unique attracting periodic orbit $\Gamma_\ve$ of the auxiliary system \eqref{repara}, which tends to the singular cycle $\Gamma_0$ as $\ve\to 0$. 
\end{theorem}

In order to prove Theorem \ref{main.result}, we need to introduce the following sections
\begin{align}\label{entry123}
\begin{split}
	& \Sigma_1 := \{(f,c,e)\in\R^3 \ | \ (f,\ e)\in R_1, \ c=\delta_1\},\\
    & \Sigma_2 := \{(f,c,e)\in\R^3 \ | \ (c,\ e)\in R_2, \ f=\delta_2\},\\
    & \Sigma_3 := \{(f,c,e)\in\R^3 \ | \ (f,\ e)\in R_3, \ c=\delta_3\},
\end{split}
\end{align}
where $R_j$ $(j=1,2,3)$ are suitable small rectangles, and $\delta_j$ are chosen sufficiently small.
Note that $\Sigma_1$ is transversal to $\omega_4$, $\Sigma_2$ is transversal to $\omega_6$, and $\Sigma_3$ is transversal to $\omega_8$, see Fig. \ref{SinCyc}.

According to the definition of the sections $\Sigma_i$, introduced in \eqref{entry123}, we define the following Poincar\'e maps for the flow of the system \eqref{repara}
\begin{equation}\label{maps123}
	\begin{split}
		& \pi_{1}: \Sigma_1\to\Sigma_2, %\qquad  \textnormal{({\e{describing the passage from $\Sigma_1$ to $\Sigma_2$ along the non-hyperbolic line $\ell_1$}})}\\
        \\
        & \pi_{2}: \Sigma_2\to\Sigma_3,
        %\qquad  \textnormal{({\e{describing the passage from $\Sigma_2$ to $\Sigma_3$ along the non-hyperbolic line $\ell_2$}})}\\
        \\
        & \pi_{3}: \Sigma_3\to\Sigma_1,
        %\qquad\qquad  \textnormal{(describing flow form $S^3_{a,\ve}$ to $S^6_{a,\ve}$ and then to $S^1_{a,\ve}$)}\\
	\end{split}
\end{equation}
{\e{where the map $\pi_1$ describes the passage from $\Sigma_1$ to $\Sigma_2$ along the non-hyperbolic line $\ell_1$, the map $\pi_2$ describes the passage from $\Sigma_2$ to $\Sigma_3$ along the non-hyperbolic line $\ell_2$, and the map $\pi_{3}$ describes the passage from $\Sigma_3$ to $\Sigma_1$.
The map $\pi_3$ consists of slow flow along $\hjk{\S_{\ve,3}^\txta}$, followed by the fast dynamics from a neighborhood of $p_1$ to a neighborhood of $p_2$, followed by the slow flow along $\hjk{\S_{\ve,6}^\txta}$ to a neighborhood of $q^e$.
Through the fast dynamics, this neighborhood is mapped to a neighborhood of $q_e$, followed by the slow flow along $\hjk{\S_{\ve,1}^\txta}$ to $\Sigma_1$.}}

We summarize the properties of the above maps in the following lemmas.
%%%--------
\begin{lemma}\label{entry1}
	If the section $\Sigma_1$ is chosen sufficiently small, then there exists $\ve_0>0$ such that the map
    \begin{equation}
		\pi_{1}: \Sigma_1\to\Sigma_2, \qquad\qquad (f,e)\mapsto (\pi_{1}^{c}(f,e,\ve), \pi_{1}^{e}(f,e,\ve)),
    \end{equation}
    {\e{is well-defined for $\ve\in[0, \ve_0]$ and smooth for $\ve\in(0, \ve_0]$.}}
    The map $\pi_{1}$ is a strong contraction with contraction rate $\exp(-K/\ve)$ for some $K>0$.
    The image of $\Sigma_1$ is a two-dimensional domain of exponentially small size, which converges to the point $q_2:=\Sigma_2\cap\omega_7$ as $\ve\to 0$.
\end{lemma}
%%%--------
\begin{lemma}\label{entry2}
	If the section $\Sigma_2$ is chosen sufficiently small, then there exists $\ve_0>0$ such that the map
    \begin{equation} 
    		\pi_{2}: \Sigma_2\to\Sigma_3, \qquad\qquad (c,e)\mapsto (\pi_{2}^{f}(c,e,\ve), \pi_{2}^{e}(c,e,\ve)),
    \end{equation}
    {\e{is well-defined for $\ve\in[0, \ve_0]$ and smooth for $\ve\in(0, \ve_0]$.}}
    The map $\pi_{2}$ is a strong contraction with contraction rate $\exp(-K/\ve)$ for some $K>0$.
    The image of $\Sigma_2$ is a two-dimensional domain of exponentially small size, which converges to the point $q_3:=\Sigma_3\cap\omega_1$ as $\ve\to 0$.
\end{lemma}
The proofs of Lemmas \ref{entry1} and \ref{entry2} are based on the blow-up analysis of the lines $\ell_1$ and $\ell_2$, respectively, which will be presented in Subsections \ref{bue} and \ref{buf}.

%%%------------------------
%
{\e{
\begin{remark}\label{fold}
%The lines $\ell_c$, when $f^*<f<1$, and $\ell^c$, when $0<f<f^*$, are non-hyperbolic, see Figs. \ref{unit.cube}, \ref{redVF3} and \ref{redVF6}.
%On the other hand, we know that the flow along the critical manifold in $S^3$ and $S^6$ reach, respectively, the lines $\ell_c$ and $\ell^c$ in finite time (see Lemmas \ref{contraction:s3} and \ref{contraction:s6}).
The points on the line $\ell_c$ when $0.5<f<1$, and on the line $\ell^c$ when $0<f<0.5$ are jump points, i.e., the trajectory switches from the slow dynamics to the fast dynamics.
Further, it is shown that this behavior is very similar to the behavior of standard slow-fast systems with two slow variables and one fast variable near a generic ``fold'' line, studied in \cite{szmolyan2004relaxation} based on the blow-up method.
The critical manifolds $\hjk{\S_{0,3}}$ and $\hjk{\S_{0,6}}$ of system \eqref{repara} can be viewed as a standard folded critical manifold, which has been straightened out by a suitable diffeomorphism.
This leads to the curved fibers of the layer problem \eqref{layer}.
Therefore, we can use the results of \cite{szmolyan2004relaxation} to understand the behavior of \eqref{repara} close to the non-hyperbolic lines $\ell_c$ and $\ell^c$. 
\end{remark}}}
The following lemma describes the map from the section $\Sigma_3$ to the section $\Sigma_1$, defined in \eqref{maps123}.
\begin{lemma}\label{entry3}
	If the section $\Sigma_3$ is chosen sufficiently small, then there exists $\ve_0>0$ such that the map
    \begin{equation}
    	\pi_{3}: \Sigma_3\to\Sigma_1, \qquad\qquad (f,e)\mapsto (\pi_{3}^{f}(f,e,\ve), \pi_{3}^{e}(f,e,\ve)),
    \end{equation}
    is well-defined for $\ve\in[0, \ve_0]$ and smooth for $\ve\in(0, \ve_0]$.
    The image of $\Sigma_3$ is an exponentially thin strip {\e{lying exponentially close to}} $S_{a,\ve}^1\cap\Sigma_1$, i.e., its width in the $f$-direction is $O(\exp(-K/\ve))$ for some $K>0$.
Moreover, {\e{$\pi_{3}(\Sigma_3)$ converges to a segment of $\hjk{\S_{0,1}^\txta}\cap\Sigma_1$ as $\ve\to 0$.}}
\end{lemma}
{\e{
\begin{proof}
%The basic idea of the proof is to construct the map for $\ve=0$, i.e. $\pi_{3,0}$, and then treat $\pi_{3,\ve}$ as an $\ve$-perturbation of $\pi_{3,0}$.
The basic idea of the proof is based on the map that has been already described in Fig.~\ref{SinCyc} for $\ve=0$, denoted by $\pi_{3}^0$, and then treat $\pi_{3}$ as an $\ve$-perturbation of $\pi_{3}^0$.
If the section $\Sigma_3$ is chosen sufficiently small, then the trajectories starting in $\Sigma_3$ can be described by the slow flow along the manifold $\hjk{\S_{\ve,3}^\txta}$ combined with the exponential contraction towards the slow manifold until they reach a neighborhood of the jump points on the line $\ell_c$.
Applying \cite[Theorem 1]{szmolyan2004relaxation} close to the jump pints, the trajectories switch from the slow dynamics to the fast dynamics, and hence pass the non-hyperbolic line $\ell_c$; this transition is well-defined for $\ve\in[0,\ve_1]$, and smooth for $\ve\in(0,\ve_1]$ for some $\ve_1>0$.
Note that \cite[Theorem 1]{szmolyan2004relaxation} guarantees that the contraction of the solutions in the $e$-direction persists during the passage through the fold-line $\ell_c$, as it is at most algebraically expanding.
After that, the solutions follow the fast dynamics $\omega_2$ until they reach a neighborhood of the point $p_2$, see Fig. \ref{SinCyc}.
Next, the solutions follow the slow flow along the manifold $\hjk{\S_{\ve,6}^\txta}$ combined with the exponential contraction towards the slow manifold until they reach a neighborhood of the point $q^e$.
Again applying \cite[Theorem 1]{szmolyan2004relaxation} close to the jump points, the solutions which are very close to the non-hyperbolic line $\ell^c$ switch from the slow dynamics to the fast dynamics, and hence pass the non-hyperbolic line $\ell^c$, where the corresponding transitions are well-defined for $\ve\in[0,\ve_2]$, and smooth for $\ve\in(0,\ve_2]$ for some $\ve_2>0$, and then follow the fast dynamics ($\omega_4$) until they reach a neighborhood of the point $q_e$.
Finally, the solutions follow the slow flow along the manifold $\hjk{\S_{\ve,1}^\txta}$ combined with the exponential contraction towards the slow manifold until they reach the section $\Sigma_1$.

Theorem 1 of \cite{szmolyan2004relaxation} implies that the map $\pi_{3}$ is at most algebraically expanding in the direction of $e$ when $\Sigma_3$ is chosen sufficiently small.
On the other hand, the slow manifold $\hjk{\S_{\ve,1}^\txta}$ is exponentially contracting in the direction of $f$ (Fenichel theory).
Therefore, the image of $\Sigma_3$ is a thin strip lying exponentially close to $\hjk{\S_{\ve,1}^\txta}\cap\Sigma_1$.
Hence, the statements of the lemma follow.
\end{proof}}}
Now we are ready to give the proof of the main result.
\newline
\textit{Proof of Theorem \ref{main.result}.}
Let us define the map $\pi:\Sigma_3\to \Sigma_3$ as a combination of the maps $\pi_{j}$ ($j=1,2,3$), described in Lemmas \ref{entry1}, \ref{entry2} and \ref{entry3}.
More precisely, we define
\begin{equation*}
	\pi=\pi_{2}\circ\pi_{1}\circ\pi_{3}:\Sigma_3\to \Sigma_3.
\end{equation*}
If the section $\Sigma_3$ is chosen sufficiently small, Lemma \ref{entry3} implies that there exists $\ve_3>0$ such that the map $\pi_{3}$ is well-defined for $\ve\in[0,\ve_3]$ and smooth for $\ve\in(0,\ve_3]$, and the image of $\Sigma_3$ is a thin strip lying exponentially close to $\hjk{\S_{\ve,1}^\txta}\cap\Sigma_1$, i.e.,  $\pi_{3}(\Sigma_3)$ is exponentially contracting with rate $\exp(-K_3/\ve)$, for some $K_3>0$, in the $f$-direction while it is bounded in the $e$-direction.

Next, if the entry section $\Sigma_1$ is chosen such that $\Sigma_1\supset\pi_{3}(\Sigma_3)$, Lemma \ref{entry1} implies that there exists $\ve_1>0$ such that the map $\pi_{1}$ is well-defined for any $\ve\in[0, \ve_1]$ and smooth for $\ve\in(0, \ve_1]$, and $\pi_{1}$ is an exponential contraction with rate $\exp(-K_1/\ve)$ for some $K_1>0$.
Finally, if the entry section $\Sigma_2$ is chosen such that $\Sigma_2\supset\pi_{1}(\Sigma_1)$, Lemma \ref{entry2} implies that there exists $\ve_2>0$ such that the map $\pi_{2}$ is well-defined for any $\ve\in[0, \ve_2]$ and smooth for any $\ve\in(0, \ve_2]$, and further, $\pi_{2}$ is an exponential contraction with rate $\exp(-K_2/\ve)$, for some $K_2>0$, such that $\Sigma_3\supset\pi_{2}(\Sigma_2)$.

Denoting $\ve_0:=\min\{\ve_1,\ve_2,\ve_3\}$ and $K:=\min\{K_1,K_2, K_3\}$, the map $\pi:\Sigma_3\to\Sigma_3$ is well-defined for any $\ve\in[0, \ve_2]$, and smooth for $\ve\in(0,\ve_0]$.
Further, based on the contracting properties of the maps $\pi_{i}$, $i=1,2,3$, we conclude that $\pi(\Sigma_3)\subset\Sigma_3$ is contraction with rate $\exp(-K/\ve)$.
{\e{The Banach fixed-point theorem implies the existence of a unique fixed point for the map $\pi$, corresponding to the attracting periodic orbit of the system \eqref{repara}.}}
Moreover, due to the last assertion of Lemmas \ref{entry1}, \ref{entry2} and \ref{entry3}, the periodic orbit $\Gamma_\ve$ tends to the singular cycle $\Gamma_0$ as $\ve\to 0$.
This completes the proof.

%% file: subfiles/BlowUp.tex
{\e{The slow-fast analysis that we have done in Section \ref{sec:geomdes} does not explain the dynamics of system~\eqref{repara} close to the non-hyperbolic lines $\ell_1$ and $\ell_2$.
As the segments $\omega_5$ and $\omega_7$ lie on these lines (see Fig. \ref{SinCyc}), we need a detailed analysis close to the lines $\ell_1$ and $\ell_2$, which is carried out in this section via the blow-up method~\cite{kuehn2015multiple, gucwa2009geometric, krupa2001extending}.
To apply this, we extend system \eqref{repara} by adding $\ve$ as a trivial dynamic variable and obtain}}
% \vspace{2cm}
% As mentioned in Section \ref{singularcycle}, just by the information that we have from the layer problem and the slow manifold, it is not possible to give the detailed description of the singular cycle $\Gamma_0$ close to the non-hyperbolic lines $\ell_1$ and $\ell_2$.
% The goal of this section is to analyze the dynamics of the system \eqref{repara} close to the lines $\ell_1$ and $\ell_2$ by blow-up method \cite{kuehn2015multiple, gucwa2009geometric, krupa2001extending}.
% In order to get insights into the non-hyperbolic parts of the critical manifold by means of blow-up method, we extend system \eqref{repara} by adding $\ve$ as a trivial dynamic variable and obtain
\begin{equation}\label{busys}
\begin{split}
	&\frac{df}{d\tau} = \left[\gamma(1 - f)(\ve + 2 f) - 2fe(\ve + 1 - f)\right] H_2^\ve(c) H_3^\ve(e),\\
	&\frac{dc}{d\tau} = \left[8(1 - c)f(\ve + 2c) - 4c(\ve + 2 - 2c)\right] H_1^\ve(f) H_3^\ve(e),\\
	&\frac{de}{d\tau} = \left[8(1-e)c(\ve + 2e)- 4e(\ve + 2 - 2e)\right] H_1^\ve(f) H_2^\ve(c),\\
    &\frac{d\ve}{d\tau} = 0, 
\end{split}
\end{equation}
where $H_1^\ve(f), H_2^\ve(c)$ and $H_3^\ve(e)$ are defined in \eqref{H123eps}.
Note that for the extended system \eqref{busys}, the lines $\ell_1\times\{0\}$ and $\ell_2\times\{0\}$ are sets of equilibria.
Due to the fact that the linearization of \eqref{busys} around these lines has quadruple zero eigenvalues, system \eqref{busys} is very degenerate close to $\ell_1\times\{0\}$ and $\ell_2\times\{0\}$.
To resolve these degeneracies, we use the blow-up method, given in next subsections.

%% file: subfiles/BlowUpEaxis.tex
\begin{figure}[t]
\centering
\begin{tikzpicture}
\pgftext{\includegraphics[scale=1]{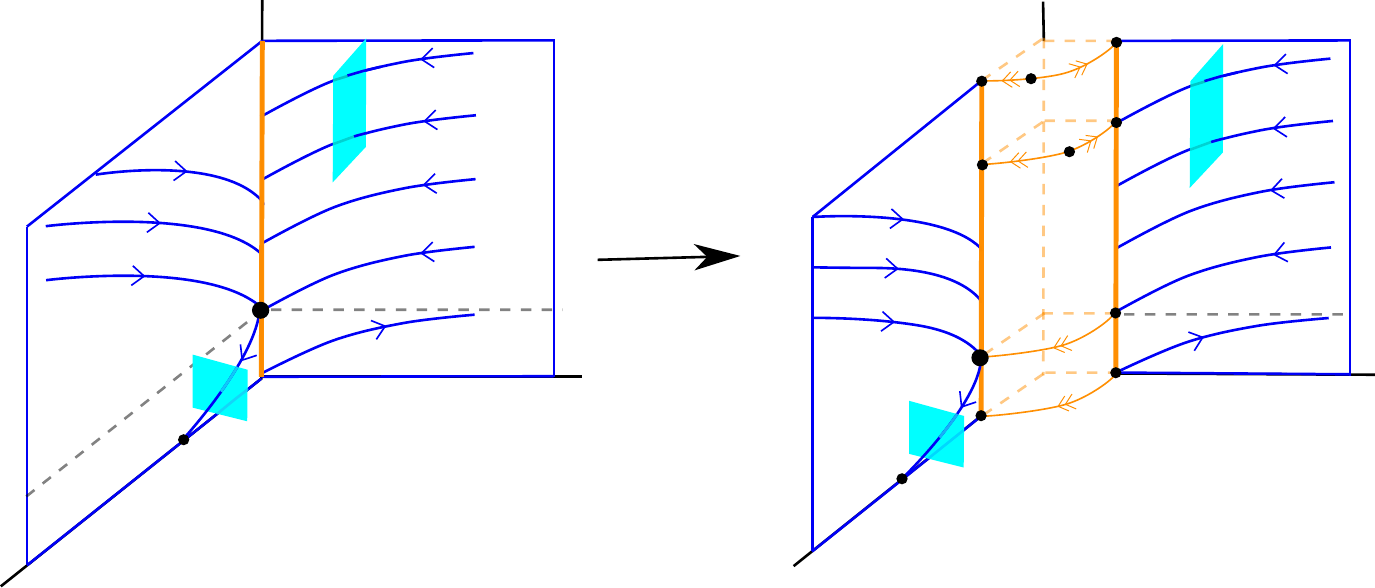}}
\node at (-7.2, -3.2){$f$};
\node at (0.8, -3.1){$\bar{f}$};
\node at (-0.9, -0.8){$c$};
\node at (7.2, -0.8){$\bar{c}$};
\node at (3.55, 3.2){$\bar{e}$};
\node at (-4.35, 3.2){$e$};
\node at (-3.05, 1.95){$\Sigma_1$};
\node at (5.65, 1.9){$\bar{\Sigma}_1$};
\node at (-4.5, 1.6){$\ell_1$};
\node at (-4.4, -1.5){$\Sigma_2$};
\node at (2.75, -2){$\bar{\Sigma}_2$};
\node at (-4.2, -0.4){$\gamma$};
\node at (2.7, -0.7){$\bar{p}_{e}$};
\node at (4.6, 1.52){$\bar{p}^{e}$};
\node at (6.4, 2.9){$\hjk{\bar\S_{0,1}^\txta}$};
\node at (1.5, 1.5){$\hjk{\bar\S_{0,2}^\txta}$};
\node at (-1.9, 2.9){$\hjk{\S_{0,1}^\txta}$};
\node at (-6.5, 1.4){$\hjk{\S_{0,2}^\txta}$};
\end{tikzpicture}
\caption{The left figure shows the dynamics close to the non-hyperbolic line $\ell_1$.
The right figure shows the corresponding dynamics in the blown-up space.}
\label{blowup.space}
\end{figure}
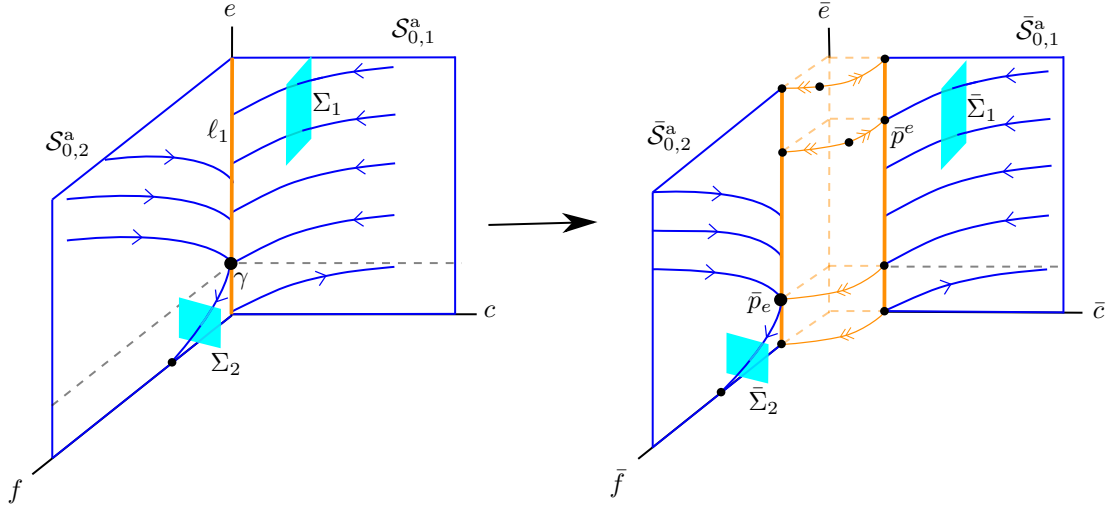

%%%%%% blow-up of the e-axis
The blow-up of the non-hyperbolic line $\ell_1\times\{0\}$ is presented in this subsection.
To this end, we transform the non-hyperbolic line of steady states $\ell_1\times\{0\}$ by
\begin{equation}\label{trans.e}
	f=r\bar{f}, \qquad\qquad c=r\bar{c}, \qquad\qquad \ve=r\bar{\ve},\qquad\qquad e=\bar{e},
\end{equation}
where $\bar{f}^2 + \bar{c}^2 + \bar{\ve}^2 = 1$ and $r\geq 0$.
Note that since $(f,c,e)\in\Q$, we may further assume that $\bar{f}, \bar{c}\geq 0$ and $\bar{e}\in[0, 1]$.
{\e{Since all weights are equal to 1 in \eqref{trans.e}, this is a homogeneous blow-up.}}
For fixed $\bar{e}$, each point $(0, 0, \bar{e})$ is blown-up to a sphere {\e{$\mathbb{S}^2$, and the line $\ell_1\times\{0\}$ is blown-up to a cylinder $\mathbb{S}^2\times[0,1]$, see Fig. \ref{blowup.space}.}}

For the analysis of system \eqref{busys} near the line $\ell_1\times\{0\}$, we define three charts $K_1, K_2$ and $K_3$ {\e{by setting $\bar{c}=1$, $\bar{\ve}=1$, and $\bar{f}=1$} in \eqref{trans.e}, respectively:}
\begin{align}
	& K_1: \qquad f=r_1 f_1, \qquad\qquad c=r_1, \qquad\qquad \ve=r_1\ve_1,\qquad\qquad e=e_1,\label{k1e}\\
	& K_2: \qquad f=r_2 f_2, \qquad\qquad c=r_2 c_2, \qquad\quad \ve=r_2,\qquad\qquad\quad  e=e_2,\label{k2e}\\
	& K_3: \qquad f=r_3, \qquad\qquad\quad  c=r_3 c_3, \qquad\quad \ve=r_3\ve_3,\qquad\qquad e=e_3,\label{k3e}
\end{align}
The changes of coordinates for the charts $K_1$ to $K_2$, and $K_2$ to $K_3$ in the blown-up space are given in the following lemma.
\begin{lemma}\label{lem:kappas}
The changes of coordinates $K_1$ to $K_2$, and $K_2$ to $K_3$ are given by
\begin{align}
	& \kappa_{12}:\qquad
    f_2 = \frac{f_1}{\ve_1}, \qquad\qquad c_2=\frac{1}{\ve_1},\qquad\qquad 
    \ve_2 = r_1\ve_1, \qquad\qquad e_2=e_1,\qquad \ve_1>0, \label{k1tok2}\\
    & \kappa_{23}:\qquad
    r_3 = r_2 f_2, \qquad\quad\ \   c_3= \frac{c_2}{f_2} ,\qquad\qquad 
    \ve_3 =\frac{1}{f_2}, \qquad\qquad \  e_3 = e_2, \qquad f_2>0.\label{k2tok3}
\end{align}
%We denote $\kappa_{ji}$ to be the inverse transformations from the chart $K_j$ to $K_i$ $(i,j=1,2,3,  \ j\neq i)$
\end{lemma}
%%%%%%%%%---------------------
The goal of this subsection is to construct the transition map $\pi_{1}:\Sigma_1\to\Sigma_2$, defined in \eqref{maps123}, and prove Lemma \ref{entry1}. 
Before going into the details, let us briefly describe our approach. 
%Subsection \ref{sec:blowup}, we shall use the blow-up method. By this we mean that 
We describe the transition map $\pi_{1}:\Sigma_1\to\Sigma_2$ via an equivalent one in the blown-up space. More specifically we define
{\e{
\begin{equation}\label{bu.map}
\pi_{1}:=\Phi\circ\bar{\pi}_1\circ\Phi^{-1},
\end{equation}%
where 
\begin{equation*}
\bar\pi_1:=\Pi_{3}\circ\kappa_{23}\circ\Pi_{2}\circ\kappa_{12}\circ\Pi_{1},
\end{equation*}
and $\Phi:\mathbb S^2\times[0,1]\times[0,r_0)\to\R^4$ is the cylindrical blow-up defined by \eqref{trans.e}, the maps $\Pi_i$ are local transitions induced by the blown-up vector fields which are detailed below, and $\kappa_{12}$ and $\kappa_{23}$ denote the changes of coordinates, given in Lemma \ref{lem:kappas}. 
$\bar\pi_1$ is the transition map}} in the blown-up space and due to the fact that $\Phi$ is a diffeomorphism, it is equivalent to $\pi_{1}$. 
A schematic of the problem at hand is shown in Fig.~\ref{blowup.space}.

{\e{The left picture in Fig. \ref{blowup.space} illustrates the critically manifolds $\hjk{\S_{0,1}^\txta}$ and $\hjk{\S_{0,2}^\txta}$, and the corresponding flows in blue.
The non-hyperbolic line $\ell_1$ is shown in orange.
For $e>\gamma$, the reduced flows on both critically manifolds approach the line $\ell_1$.
At the point on the line $\ell_1$ with $e=\gamma$, a transition from $\hjk{\S_{0,1}^\txta}$ to $\hjk{\S_{0,2}^\txta}$ is possible as indicated in the figure.
The right picture in Fig. \ref{blowup.space} schematically shows the configuration in the blown-up space.
The cylinder corresponding to $r=0$ is show in orange.
The part of the phase space corresponding to $\bar{\ve}=0$ and $r>0$ are shown outside of the the cylinder.
Here we recover the layer problem, the critically manifolds, and the reduced flows in $\hjk{\bar\S_{0,1}^\txta}$ and $\hjk{\bar\S_{0,2}^\txta}$.
In the blown-up space, the manifolds $\hjk{\S_{0,1}^\txta}$ and $\hjk{\S_{0,2}^\txta}$ are separated and hence gained hyperbolicity, \hjk{in particular they are} attractive, as indicated \hjk{below in Fig. \ref{fig:chartK11}}.
All these assertions will be proven in this  section.}}

{\e{Roughly speaking, in chart $K_1$ we continue the attracting slow manifold $\hjk{\bar\S_{0,1}^\txta}$ onto the cylinder.
Chart $K_2$ is used to track the flow across the cylinder.
The exit of the flow from the cylinder and its transition to $\hjk{\bar\S_{0,2}^\txta}$ is studied in chart $K_3$, see Figs. \ref{blowup.space} and \ref{CompleteblowupE}.
The detailed analysis of the maps $\Pi_i$ introduced in \eqref{bu.map}, is given in the forthcoming subsections.}}
% In few words, in chart $K_1$ we study the way trajectories approach the non-hyperbolic line $\ell_1$. Next, in chart $K_2$ we detail the behavior of the flow within a small neighborhood of $\ell_1$. Finally, in chart $K_3$ we describe the way trajectories leave $\ell_1$. The formal description of the maps $\Pi_i$ is provided in the forthcoming subsections.

%%%%%%%%%%%%%%%%%%---------------------------
\subsubsection{Analysis in chart $K_1$}
After substituting \eqref{k1e} into \eqref{busys} and dividing out all the equations by the common factor $r_1$, the equations governing the dynamics in chart $K_1$ are given by
\begin{equation}\label{bue.k11}
\begin{split}
	& f_1'= -4f_1 \Gamma_1 G_{11} + [\gamma(1-r_1f_1)(\ve_1+2f_1) - 2f_1e_1(r_1\ve_1+1-r_1f_1)]G_{12},\\
	& r_1' = 4 r_1 \Gamma_1 G_{11},\\
    & e_1' = 4r_1[2r_1(1-e_1)(r_1\ve_1+2e_1) - e_1(r_1\ve_1+2-2e_1)]G_{13},\\
	& \ve_1' = -4\ve_1[2r_1f_1(1-r_1)(\ve_1+2) - (r_1\ve_1+2-2r_1)]G_{11},
\end{split}
\end{equation}
where we denote
\begin{equation}
\begin{split}
	& \Gamma_1 := [2r_1f_1(1-r_1)(\ve_1+2) - (r_1\ve_1+2-2r_1)],\\
	& G_{11} := 	(r_1\ve_1 + 1 - r_1f_1)(r_1\ve_1 + 2 - 2e_1)(\ve + 2f_1)(r_1\ve_1 + 2e_1),\\
    & G_{12} := (r_1\ve_1+2-2r_1)(r_1\ve_1 + 2 - 2e_1)(\ve_1+2)(r_1\ve_1 + 2e_1),\\
	& G_{13} := (r_1\ve_1+1-r_1f_1)(r_1\ve_1+2-2r_1)(\ve_1+2f_1)(\ve_1+2).
\end{split}
\end{equation}
% %%%%%%%%%%%%%%%%%%--------------------------------------------
% %%%%% this is the extended version of the above equations and they are CORRECT
% \begin{equation}\label{bue.k11}
% \begin{split}
% 	& f_1'= -4f_1[2r_1f_1(1-r_1)(\ve_1+2) - (r_1\ve_1+2-2r_1)]G_{11}  + [\gamma(1-r_1f_1)(\ve_1+2f_1) - 2f_1e_1(r_1\ve_1+1-r_1f_1)]G_{12},\\
% 	& r_1' = 4 r_1[2 r_1 f_1(1-r_1)(\ve_1 + 2) - (r_1\ve_1+2-2r_1)]G_{11},\\
%     & e_1' = 4r_1[2r_1(1-e_1)(r_1\ve_1+2e_1) - e_1(r_1\ve_1+2-2e_1)]G_{13},\\
% 	& \ve_1' = -4\ve_1[2r_1f_1(1-r_1)(\ve_1+2) - (r_1\ve_1+2-2r_1)]G_{11},
% \end{split}
% \end{equation}
% where
% \begin{equation}
% \begin{split}
% 	& G_{11} = 	(r_1\ve_1 + 1 - r_1f_1)(r_1\ve_1 + 2 - 2e_1)(\ve + 2f_1)(r_1\ve_1 + 2e_1),\\
%     & G_{12} = (r_1\ve_1+2-2r_1)(r_1\ve_1 + 2 - 2e_1)(\ve_1+2)(r_1\ve_1 + 2e_1),\\
% 	& G_{13} = (r_1\ve_1+1-r_1f_1)(r_1\ve_1+2-2r_1)(\ve_1+2f_1)(\ve_1+2).
% \end{split}
% \end{equation}
% %%%%%%%%%%%%%%%%--------------------------------------
From \eqref{bue.k11} it is clear that the planes $r_1=0$ and $\ve_1=0$ are invariant.
Hence, we consider the following cases:
\begin{enumerate}
\item $r_1=\ve_1=0$: in this case, the dynamics \eqref{bue.k11} is simplified to
\begin{equation}\label{bue.k12}
\begin{split}
	& e_1'=0,\\
    & f_1'= 32f_1e_1(1-e_1)[2f_1+\gamma-e_1].
\end{split}
\end{equation}
For fixed $e$, the equilibria of the system \eqref{bue.k12} are the attracting point $p_1^a=(f_1,r_1,e_1,\ve_1)=(0,0,e_1,0)$, and the repelling point $p_1^r=(f_1,r_1,e_1,\ve_1)=(\frac{e_1-\gamma}{2},0,e_1,0)$.
Note that the two hyperbolic points $p_1^a$ and $p_1^r$ intersect at the non-hyperbolic point $(f_1,r_1,e_1,\ve_1)=(0,0,\gamma,0)$, see Fig. \ref{fig:chartK11}.
%%-----------------
\item $\ve_1=0$: in this case, the dynamics \eqref{bue.k11} is represented by
\begin{equation}\label{bue.k13}
\begin{split}
	& f_1' = 32f_1e_1(1-e_1)(1-r_1)(1-r_1f_1)[(\gamma-e_1) - 2f_1(2r_1f_1-1)],\\
	& r_1' = 64r_1f_1e_1(1-e_1)(1-r_1)(1-r_1f_1)[2r_1f_1-1],\\
    & e_1' = 64r_1f_1e_1(1-e_1)(1-r_1)(1-r_1f_1)[2r_1-1].
\end{split}
\end{equation}
From \eqref{bue.k13}, one concludes that the plane $f_1=0$ is the plane of equilibria which is denoted by $\hjk{\bar\S_{0,1}^\txta}$, see Fig. \ref{fig:chartK11}.
The non-zero eigenvalue along $\hjk{\bar\S_{0,1}^\txta}$ is given by $\lambda=32e_1(1-e_1)(1-r_1)(\gamma-e_1)$.
For $0\leq r_1<1$ and $e_1>\gamma$, the plane $\hjk{\bar\S_{0,1}^\txta}$ is attracting.
As the $e_1$-axis is a part of $\hjk{\bar\S_{0,1}^\txta}$, we denote that part of the $e_1$-axis that $\gamma\leq e_1\leq 1$ by $\ell_{e_1}$.
We also have another curve of equilibria which is defined by $r_1=0$, and $f_1=\frac{e_1-\gamma}{2}$, denoted by $\mathcal{M}_1^r$, see Fig. \ref{fig:chartK11}.
This curve of equilibria is of saddle-type with the eigenvalues $\lambda=\pm 32e_1(e_1-1)(e_1-\gamma)$.
Note that we have recovered the information of the previous case here. 
% \begin{figure}[H]
% \centering
% \begin{tikzpicture}
% \pgftext{\includegraphics[scale=1]{figures/chartK11}}
% \node at (3.3, -0.8){$r_1$};
% \node at (0, 0.9){$p_1^a$};
% \node at (-1.2, 0.5){$p_1^r$};
% \node at (3.1, 2.2){$\hjk{\bar\S_{0,1}^\txta}$};
% \node at (-0.55, 2.7){1};
% \node at (-0.15, -0.4){$\gamma$};
% \node at (-3.3, -3){$f_1$};
% \node at (-0.2, 3.15){$e_1$};
% \node at (-1, 1.2){$\ell_{e_1}$};
% \node at (-3, 0.05){$\mathcal{M}_1^r$};
% \end{tikzpicture}
% \caption{Dynamics of \eqref{bue.k11} in the invariant plane $\ve_1=0$}
% \label{fig:chartK11_old}
% \end{figure}
%%-----------------
\item $r_1=0$: in this case, the dynamics \eqref{bue.k11} is represented by
\begin{equation}\label{bue.k14}
\begin{split}
	& e_1' = 0,\\
    & f_1' = 8e_1(1-e_1)[(\ve_1+2)(\gamma(\ve_1+2f_1) - 2f_1e_1) + 4f_1(\ve_1+2f_1)],\\
    & \ve_1' = 32e_1\ve_1(1-e_1)(\ve_1+2f_1).
\end{split}
\end{equation}
By setting $\ve_1=0$, we again have the line $\ell_{e_1}$ and the curve $\mathcal{M}_1^r$.
The Jacobian matrix at a point in $\ell_{e_1}$ has two eigenvalues: one zero and the other one is $\lambda= 32e_1(1-e_1)(\gamma-e_1)$.
So the line $\ell_{e_1}$ is attracting when $e>\gamma$.
As in this case we have two zero eigenvalues, it implies that there exists a two-dimensional center manifold, namely, $\mathcal{C}_{a,1}$.
\begin{remark}
In chart $K_1$, the most important role is played by the two-dimensional center manifold $\mathcal{C}_{a,1}$, see Lemma \ref{centermanifold.k1e}.
In fact, this is the continuation of the critical manifold $\hjk{\bar\S_{0,1}^\txta}$. 
\end{remark}
% \textbf{\begin{figure}[H]
% \centering
% \begin{tikzpicture}
% \pgftext{\includegraphics[scale=1]{figures/chartK12}}
% \node at (3.3, -0.8){$\ve_1$};
% \node at (-0.6, 1){$p_1^a$};
% \node at (-1.25, 0.5){$p_1^r$};
% \node at (-0.3, 3.2){$e_1$};
% \node at (-0.15, -0.4){$\gamma$};
% \node at (-3.2, -3.05){$f_1$};
% \node at (-0.1, 2.65){1};
% \node at (0.7, 1.9){$\ell_{e_1}$};
% \node at (-3, 0){$\mathcal{M}_1^r$};
% \end{tikzpicture}
% \caption{Dynamics of system \eqref{bue.k11} in the invariant plane $r_1=0$}
% \label{chartK12_old}
% \end{figure}}
\end{enumerate}

\hjk{

\begin{figure}[htbp]
    \centering
    \begin{subfigure}[t]{0.45\textwidth}\centering
    \begin{tikzpicture}
    \pgftext{\includegraphics[scale=1]{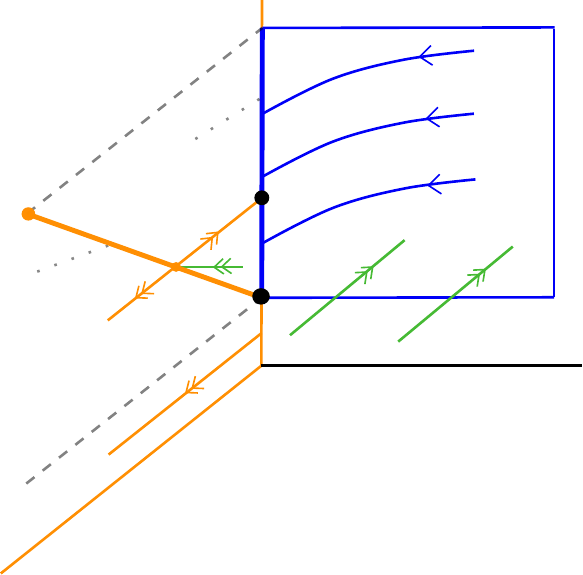}}
    \node at (3.3, -0.8){$r_1$};
    \node at (0, 0.9){$p_1^a$};
    \node at (-1.2, 0.5){$p_1^r$};
    \node at (3.1, 2.2){$\hjk{\bar\S_{0,1}^\txta}$};
    \node at (-0.55, 2.7){1};
    \node at (-0.15, -0.4){$\gamma$};
    \node at (-3.3, -3){$f_1$};
    \node at (-0.2, 3.15){$e_1$};
    \node at (-1, 1.2){$\ell_{e_1}$};
    \node at (-3, 0.05){$\mathcal{M}_1^r$};
\end{tikzpicture}
    \caption[]{Dynamics of \eqref{bue.k11} restricted to $\ve_1=0$.}
    \label{fig:chartK11}
    \end{subfigure}\hfill
    \begin{subfigure}[t]{0.45\textwidth}\centering
    \begin{tikzpicture}
    \pgftext{\includegraphics[scale=1]{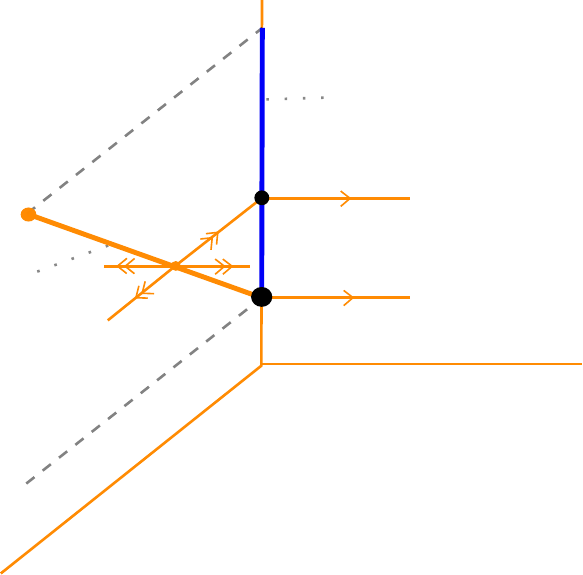}}
    \node at (3.3, -0.8){$\ve_1$};
    \node at (-0.6, 1){$p_1^a$};
    \node at (-1.25, 0.5){$p_1^r$};
    \node at (-0.3, 3.2){$e_1$};
    \node at (-0.15, -0.4){$\gamma$};
    \node at (-3.2, -3.05){$f_1$};
    \node at (-0.1, 2.65){1};
    \node at (0.7, 1.9){$\ell_{e_1}$};
    \node at (-3, 0){$\mathcal{M}_1^r$};
\end{tikzpicture}
    \caption[]{Dynamics of \eqref{bue.k11} restricted to $r_1=0$.}
    \label{fig:chartK12}
    \end{subfigure}
    \caption[]{Dynamics of \eqref{bue.k11} restricted to invariant subspaces.}
    \label{fig:my_label}
\end{figure}
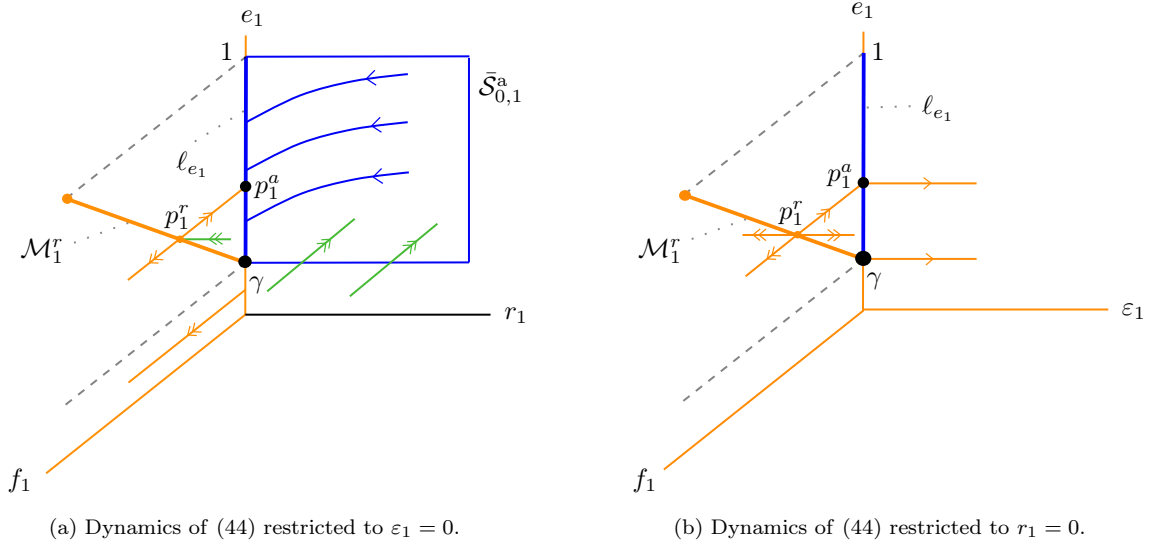
}

We summarize the analysis performed in this subsection in the following lemmas.
\begin{lemma}
System \eqref{bue.k11} has the following manifolds of equilibria: 
\begin{enumerate}
\item The plane $\hjk{\bar\S_{0,1}^\txta}$ which includes the line $\ell_{e_1}$,
\item $\mathcal{M}_1^r=\{(f_1,r_1,e_1,\ve_1) \ | \  f_1=\frac{e_1-\gamma}{2}   , \ r_1=0, \ e_1\in[\gamma, 1], \ \ve_1=0\}$
\end{enumerate}
\end{lemma}
\begin{lemma}\label{centermanifold.k1e}
The following properties hold for system \eqref{bue.k11}:
\begin{enumerate}
\item The linearization of \eqref{bue.k11} along $\hjk{\bar\S_{0,1}^\txta}$ has three zero eigenvalues, and  the nonzero eigenvalue $\lambda=32e_1(1-e_1)(1-r_1)(\gamma-e_1)$, which for $r_1=0$ corresponds to the flow in the invariant plane $(f_1,e_1)$. 
\item There exists a three-dimensional center manifold $\mathcal{W}_{a,1}^c$ of the line $\ell_{e_1}$ which contains the plane of equilibria $\hjk{\bar\S_{0,1}^\txta}$ and the two-dimensional center manifold $\mathcal{C}_{a,1}$.
The manifold $\mathcal{W}_{a,1}^c$ is attracting, and in the set $D_1$, defined by
\begin{equation*}
D_1:=\{(f_1,r_1,e_1,\ve_1) \ | \ 0\leq r_1\leq \delta_1, \ , \ e_1\in I_1, \ 0\leq\ve_1\leq\alpha_1\},
\end{equation*}
is given by the graph    
\begin{equation*}
	f_1=h_{a,1}(r_1, e_1, \ve_1),
\end{equation*}
where $I_1$ is a suitable interval, and $\alpha_1, \delta_1>0$ are sufficiently small.
For the particular point $p_{a,1}\in\ell_{e_1}$ where $e^0\in I_1$, the function $h_{a,1}(r_1, e^0, \ve_1)$ has the expansion
\begin{equation}\label{2dcenter.k1}
	h_{a,1}(r_1, e^0, \ve_1)= \frac{\gamma}{2(e^0-\gamma)}\ve_1 + O(\ve_1^2).
\end{equation}
\item There exists $K>0$ such that the orbits that are near the center manifold $\mathcal{W}_{a,1}^c$ are attracted to $\mathcal{W}_{a,1}^c$ by an exponential rate of order $O(\exp(-Kt_1))$.
\end{enumerate}
\end{lemma}
\begin{proof}
A straightforward calculation shows the first claim.
Due to the fact that the linearization of \eqref{bue.k11} along $\hjk{\bar\S_{0,1}^\txta}$ has three zero eigenvalues, there exists \cite{chow1994normal,guckenheimer2013nonlinear} an attracting three-dimensional center manifold $\mathcal{W}_{a,1}^c$ at the point $p_{a,1}$.
To derive equation \eqref{2dcenter.k1}, we first expand $f_1$ to the first order of variables $r_1, e_1$ and $\ve_1$, and then plug into \eqref{bue.k11}.
By comparing the coefficients of $r_1, e_1$ and $\ve_1$, equation \eqref{2dcenter.k1} is obtained.  
The last claim is proven by the center manifold theory %\cite{chow1994normal,guckenheimer2013nonlinear} 
applied at the point $p_{a,1}$.
\end{proof}
{\e{\begin{remark}
The attracting center manifold $\mathcal{W}_{a,1}^c$ recovers parts of the slow manifold {\hjk{$\S_{\ve,1}^\txta$}} away form the line $\ell_1$, and extends it into an $O(\ve)$ neighborhood of $\ell_1$.
The slow manifold {\hjk{$\S_{\ve,1}^\txta$}} is obtained as a section $\ve=constant$ of $\mathcal{W}_{a,1}^c$.
In chart $K_1$, this center manifold is given by the graph \eqref{2dcenter.k1}.
\end{remark}}}
Note that in chart $K_1$, our goal is to understand the dynamics \eqref{bue.k11} close to the center manifold $\mathcal{W}_{a,1}^c$, which corresponds to a sufficiently small neighborhood of the slow manifold $\hjk{\bar\S_{0,1}^\txta}$.
Assume that $\delta_1, \alpha_1, \beta_1>0$ are small constants.
Let us define the sections

\begin{equation}
\begin{split}
\Delta_1^{in} &:=\{(f_1, r_1, e_1, \ve_1) \ | \  (f_1, r_1, e_1, \ve_1)\in D_1, \  r_1=\delta_1\},\\
\Delta_1^{out} &:=\{(f_1, r_1, e_1, \ve_1) \ | \  (f_1, r_1, e_1, \ve_1)\in D_1, \ \ve_1=\alpha_1\},\\
R_1^{in} &:=\{(f_1, r_1, e_1, \ve_1) \ | \  (f_1, r_1, e_1, \ve_1)\in D_1, \  r_1=\delta_1, |f_1|\leq \beta_1\}.
\end{split}
\end{equation}
{\e{Note that by the way we have defined $\Delta_1^{in}$, we in fact have $\Delta_1^{in}=\bar{\Sigma}_1:=\Phi^{-1}(\Sigma_1\times\{[0,\rho_1]\})$ for some $\rho_1>0$, see Fig.~\ref{blowup.space}.}} 
Furthermore, the constants $\delta_1, \alpha_1, \beta_1$ are chosen such that $R_1^{in}\subset\Delta_1^{in}$, and the intersection of the center manifold $\mathcal{W}_{a,1}^c$ with $\Delta_1^{in}$ lies in $R_1^{in}$, i.e., $\mathcal{W}_{a,1}^c\cap\Delta_1^{in}\subset R_1^{in}$.

Let us denote $\Pi_{1}$ as the transition map from $\Delta_1^{in}$ to $\Delta_1^{out}$, induced by the flow of \eqref{bue.k11}.
In order to construct map $\Pi_{1}$, we reduce system \eqref{bue.k11} to the center manifold $\mathcal{W}_{a,1}^c$ and analyze the system based on the the dynamics on $\mathcal{W}_{a,1}^c$. 
To this end, by substituting \eqref{2dcenter.k1} into \eqref{bue.k11} and rescaling time, the flow of the center manifold is given by 
\begin{equation}\label{center.manifold11}
\begin{split}
& r_1' = -r_1,\\
& e_1' = -\frac{1}{2}[O(r_1)+O(r_1\ve_1)],\\
& \ve_1'= \ve_1,
\end{split}
\end{equation}
where the derivative is with respect to the new timescale, namely, $t_1$.
Now let us consider a solution of \eqref{center.manifold11}, namely, $(r_1(t_1), e_1(t_1), \ve_1(t_1))$  which satisfies the following conditions:
\begin{equation}
\begin{split}
& r_1(0)=\delta_1, \qquad\qquad r_1(T^{out})=r_1^{out},\\
& e_1(0)=e_1^{in}, \qquad\qquad e_1(T^{out})=e_1^{out},\\
& \ve_1(0)=\ve_1^{in}, \qquad\qquad \ve_1(T^{out})=\alpha_1.
\end{split}
\end{equation}
From equation $\ve_1'= \ve_1$ with the conditions $\ve_1(0)=\ve_1^{in}$ and $\ve_1(T^{out})=\alpha_1$, we can calculate the time that $(r_1(t_1), e_1(t_1), \ve_1(t_1))$ needs to travel from $\Delta_1^{in}$ to $\Delta_1^{out}$, which is given by
\begin{equation}\label{time.out1}
T^{out}=\ln\frac{\alpha_1}{\ve_1^{in}}.
\end{equation}
Since $e_1' = -\frac{1}{2}[O(r_1)+O(r_1\ve_1)]$ with $e_1(T^{in})=e_1^{in}$, we can estimate the time evolution of $e_1(t_1)$, which is given by 
\begin{equation}
e_1(t_1)=\frac{r_1^{in}}{2}\left[\exp(-t_1)-1-t_1\ve_1^{in}\right] + e_1^{in}, \qquad 0\leq t_1\leq T^{out}.
\end{equation}
Hence, in view of \eqref{time.out1}, one has
\begin{equation}\label{e1out}
e_1(T^{out})=e_1^{out}:=\frac{r_1^{in}}{2}\left(\frac{\ve_1^{in}}{\alpha_1}-1-\ve_1^{in}\ln\frac{\alpha_1}{\ve_1^{in}}\right) + e_1^{in}.
\end{equation}
We summarize the analysis performed for chart $K_1$ in the following theorem.
\begin{theorem}
For system \eqref{bue.k11} with sufficiently small $\delta_1, \alpha_1, \beta_1$ and $R_1^{in}\subset\Delta_1^{in}$, the transition map $\Pi_{1}:R_1^{in}\to\Delta_1^{out}$ is well-defined and has the following properties:
\begin{enumerate}
\item $\Pi_{1}(R_1^{in})\subset\Delta_1^{out}$ is a three-dimensional wedge-like region in $\Delta_1^{out}$.
\item The transition map $\Pi_{1}$ is given by
\begin{equation*}
\Pi_{1}
\begin{pmatrix}
f_1\\
\delta_1\\
e_1\\
\ve_1
\end{pmatrix}
=
\begin{pmatrix}
h_{a,1}(\frac{\delta_1}{\alpha_1}\ve_1, e_1^{out}, \alpha_1) + \Psi(\delta_1,e_1,\ve_1)\\
\frac{\delta_1}{\alpha_1}\ve_1\\
e_1^{out}\\
\alpha_1
\end{pmatrix},
\end{equation*}
where $e_1^{out}$ is given in \eqref{e1out}, $\Psi(\cdot)$ is an exponentially small function, and $h_{a,1}(\cdot)$ is of order $O(\ve_1)$, due to \eqref{2dcenter.k1}.
\end{enumerate}
\end{theorem}
%%%%%%%%%%%%%%%%%%%%%%%%%%%%%%%%%%%-----------------------------
\subsubsection{Analysis in chart $K_2$}
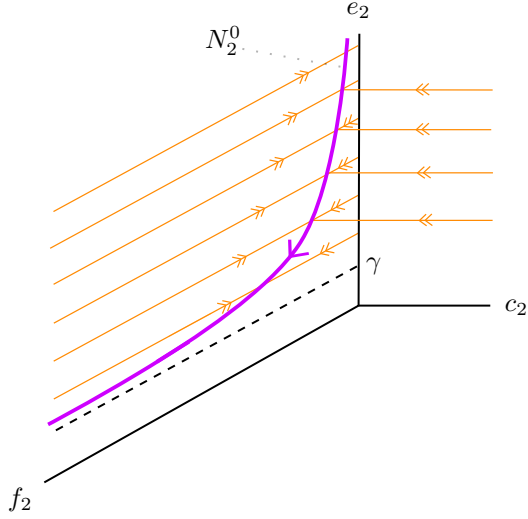
\begin{figure}[h]
\centering
\begin{tikzpicture}
\pgftext{\includegraphics[scale=1]{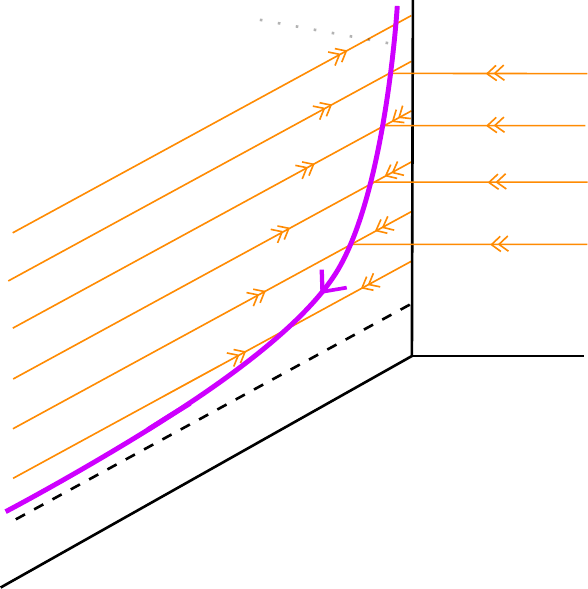}}
\node at (1.4, -0.1){$\gamma$};
\node at (3.3, -0.65){$c_2$};
\node at (1.2, 3.3){$e_2$};
\node at (-0.6, 2.9){$N_2^0$};
\node at (-3.3, -3.2){$f_2$};
%%%%%
\end{tikzpicture}
\caption{Fully attracting critical manifold $N_2^0$ in purple, and the slow and fast dynamics in chart $K_2$.}
\label{chartK2e}
\end{figure}

After substituting \eqref{k2e} into \eqref{busys} and dividing out all the equations by the common factor $r_2$, the equations governing the dynamics in chart $K_2$ are given by
\begin{equation}\label{bue.k21}
\begin{split}
	& f_2' = 8 e_2 \left[\gamma(1 + 2 f_2) - 2 f_2 e_2\right](1 - e_2)(1 + 2c_2) + O(\ve),\\
	& c_2' = -32 c_2 e_2 (1 - e_2)(1 + 2f_2) + O(\ve),\\
    & e_2' = -16\ve e_2 (1-2e_2)(1+2f_2)(1+2c_2) + O(\ve^2),\\
    & \ve'=0.
\end{split}
\end{equation}
Due to the fact that $r_2=\ve$ in chart $K_2$, we have presented \eqref{bue.k21} in terms of $\ve$.
Note that since $r_2'=\ve'=0$, system \eqref{bue.k21} is a family of three-dimensional vector fields which are parametrized by $\ve$.
Moreover, system \eqref{bue.k21} is a slow-fast system in the standard form, i.e., $e_2$ is the slow variable, and $f_2$ and $c_2$ are the fast variables.
Since the differentiation $'$ in \eqref{bue.k21} is with respect to the fast time variable, namely $\tau_2$, by transforming it to the slow time variable we have $t_2=\ve\tau_2$, and hence
% %%%%------------------
% \textcolor{red}{\begin{equation}\label{bue.k22}
% \begin{split}
% 	& \ve \dot{f}_2 = \left[\gamma\left(1- \ve f_2\right)\left(1+2 f_2\right) - 2 f_2 e_2\left(\ve + 1 - \ve f_2\right)\right]G_{21}(c_2, e_2, \ve),\\
% 	& \ve \dot{c}_2 = \left[8 \ve f_2\left(1 - \ve c_2\right)\left(1 + 2 c_2\right) - 4 c_2\left(\ve + 2 - 2 \ve c_2\right)\right]G_{22}(f_2, e_2, \ve),\\
% 	& \dot{e}_2 = \left[8 \ve c_2 \left(1 - e_2\right)\left(\ve + 2 e_2\right) - 4 e_2\left(\ve + 2 - 2 e_2\right)\right]G_{23}(f_2, c_2, \ve),
% \end{split}
% \end{equation}}
% %%%%%------------
\begin{equation}\label{bue.k22}
\begin{split}
	& \ve \dot{f}_2 = 8 e_2 \left[\gamma(1 + 2 f_2) - 2 f_2 e_2\right](1 - e_2)(1 + 2c_2) + O(\ve),\\
	& \ve \dot{c}_2 = -32 c_2 e_2 (1 - e_2)(1 + 2f_2) + O(\ve),\\
    & \dot{e}_2 = -16 e_2 (1-2e_2)(1+2f_2)(1+2c_2) + O(\ve),
\end{split}
\end{equation}
where the derivative is with respect to $t_2$.
Now by setting $\ve=0$ in \eqref{bue.k21} we obtain the corresponding layer problem
\begin{equation}\label{bue.k23}
\begin{split}
	& f_2' = 8 e_2 \left[\gamma(1 + 2 f_2) - 2 f_2 e_2\right](1 - e_2)(1 + 2c_2),\\
	& c_2' = -32 c_2 e_2 (1 - e_2)(1 + 2f_2),\\
    & e_2' = 0,
\end{split}
\end{equation}
which has the associated critical manifold $c_2=0$ and $f_2 = \frac{\gamma}{2(e_2 - \gamma)}$, denoted by $N_2^0$ (see Fig. \ref{chartK2e}).
The Jacobian matrix corresponding to \eqref{bue.k23} along this critical manifold has the eigenvalues
\begin{equation}\label{eigk2e}
	\lambda_{21} = -16e_2(1-e_2)(e_2-\gamma), \qquad\qquad
    \lambda_{22} = \frac{32 e_2^2 (e_2-1)}{(e_2 - \gamma)}.
\end{equation}
{\e{As it is clear form \eqref{eigk2e}, the critical manifold restricted to $e_2\in(\gamma, 1)$ is normally hyperbolic, and specially, is fully attracting since both of the eigenvalues are negative.
As $e_2$ approaches $\gamma$ from above, $f_2$ develops a singularity along $N_2^0$.
Thus, the behavior of $N_2^0$ as $e\to\gamma$ has to be studied in chart $K_3$.}}
Using Fenichel theory and the dynamics in chart $K_2$ for $\ve=0$, one is able to describe the dynamics for $0<\ve\ll 1$ in this chart, i.e., there exists a slow manifold $N_2^\ve$ which is the $\ve$-perturbation of $N_2^0$.
We summarize the properties of the critical manifold of chart $K_2$ in the following lemma.
\begin{lemma}\label{lem:crit2}
The critical manifold
\begin{equation}\label{m2e}
	N_2^0 = \{(f_2, c_2, e_2) \ | \ f_2 = \frac{\gamma}{2(e_2 - \gamma)}, \ c_2 = 0, \ e_2\in I_2^0\},
\end{equation}
is fully attracting, where $I_2^0$ is a compact subset of the interval $(\gamma, 1)$.
Moreover, there exists $\ve_0>0$ such that for any $\ve\in(0,\ve_0)$, there exists a smooth locally invariant attracting one-dimensional slow manifold $N_2^\ve$, which is $O(\ve)$-close to $N_2^0$, with the slow flow
    \begin{equation}\label{eq:bue_slow}
    	\dot{e}_2= -4e_2(\ve+2-2e_2)(\ve+1-\ve f_2)(\ve+2)(1+2f_2).
    \end{equation}
\end{lemma}
Note that $e_2$ is decreasing along $N_2^\ve$, see Fig. \ref{chartK2e}. 
%%%%-----------------------------
Now, we construct the transition map $\Pi_{2}$. For this let us define the sections
% Note that the outgoing section $\Delta_1^{out}$ in chart $K_1$ is mapped to the incoming section $\Delta_2^{in}$ by the transformation $\kappa_{12}$, defined in \eqref{k1tok2}.
% We define $\Delta_2^{in}$ and $\Delta_2^{out}$ in chart $K_2$ as follows:
\begin{align*}
	& \Delta_2^{in} := \{(f_2,c_2,e_2, \ve) \ | \ f_2\in[0, \beta_2], \ c_2=\frac{1}{\alpha_1}, \ e_2\in I_2, \ \ve\in[0, \alpha_2]\},\\
    & \Delta_2^{out} := \{(f_2, c_2, e_2, \ve) \ | \ f_2=\beta_2, \  c_2\in[0, \frac{1}{\alpha_1}], \ e_2\in I_2, \ \ve\in[0, \alpha_2]\}.
\end{align*}
where $\delta_1$ is small, $\beta_2=\frac{\beta_1}{\alpha_1}$, $\alpha_2=\delta_1\alpha_1$, and $I_2$ is a suitable interval. Note that $\Delta_2^{in}=\kappa_{12}(\Delta_1^{out})$. 
Let us define the transition map from $\Delta_2^{in}$ to $\Delta_2^{out}$ as follows:
\begin{equation}\label{mappi2}
	\Pi_{2}: \Delta_2^{in} \to \Delta_2^{out}, \qquad 
        (f_2^{in}, \frac{1}{\alpha_1}, e_2^{in}, \ve) \mapsto (\beta_2, c_2^{out}, e_2^{out}, \ve),
\end{equation}
The map $\Pi_{2}$ is described by the Fenichel theory, i.e., all orbits starting from $(f_2^{in}, \frac{1}{\alpha_1}, e_2^{in}, \ve)$ are attracted by the slow manifold $N_2^\ve$, follow the slow manifold along $N_2^\ve$, and then after some time intersect the section $\Delta_2^{out}$ transversally.
\begin{remark}
In the limit $\ve=0$, the map $\Pi_{2}$ is defined by first projecting $(f_2,e_2)\in\Delta_2^{in}$ onto $N_2^0$ along the stable foliation, and then by following the slow flow \eqref{eq:bue_slow}.
\end{remark}
We summarize the analysis performed in chart $K_2$ in the following lemma.
\begin{lemma}\label{inout.k2e} 
For small $\alpha_1>0$, there exists a sufficiently small $\alpha_2>0$ such that the transition map $\Pi_2$, defined in \eqref{mappi2}, is well-defined.
Moreover, for $\ve=constant$, $\Pi_{2}$ is contracting with the contraction rate $\exp(-K/\ve)$ for some $K>0$.
\end{lemma}
\begin{proof}
The transition map $\Pi_{2}: \Delta_2^{in} \to \Delta_2^{out}$ is described by Fenichel theory, i.e., all orbits starting from $\Delta_2^{in}$ are attracted by the slow manifold $N_2^\ve$, with a contraction rate $\exp(-K/\ve)$ for some $K>0$, and after some time they reach the section $\Delta_2^{out}$.
\end{proof}
{\e{
\begin{remark}
The slow manifold $N_2^\ve$ corresponds to the perturbation of $N_2^0$ when $\ve=constant$.
The family of all such manifolds is denoted by $\mathcal{N}_2$.
\end{remark}
}}
%%%%%%%%%%%%%%%%%%-----------------------------------

\subsubsection{Analysis in chart $K_3$}
Solutions in chart $K_2$ which reach the section $\Delta_2^{out}$ must be continued in chart $K_3$.
For this reason, we continue our analysis in chart $K_3$.
After substituting \eqref{k3e} into \eqref{busys} and dividing out all the equations by the common factor $r_3$, the equations governing the dynamics in chart $K_3$ are given by
\begin{equation}\label{bue.k31}
\begin{split}
    & r_3' = r_3 \Gamma_3 G_{31},\\
    & c_3' = -c_3 \Gamma_3 G_{31} + \left[8r_3(1-r_3c_3)(\ve_3+2c_3) - 4c_3(r_3\ve_3+2-2r_3c_3)\right] G_{32},\\
    & e_3' = r_3  \left[8r_3c_3(1-e_3)(r_3\ve_3+2e_3) - 4e_3(r_3\ve_3+2-2e_3)\right] G_{33},\\
    & \ve_3' = -\ve_3 \Gamma_3 G_{31},\\
\end{split}
\end{equation}
where we denote
%%%%%%%%%%-------------------------------------------
\begin{equation*}
\begin{split}
	& \Gamma_3 := \left[\gamma(1-r_3)(\ve_3+2) - 2e_3(r_3\ve_3+1-r_3)\right], \\
    & G_{31} := (r_3\ve_3+2-2r_3c_3)(\ve_3+2c_3)(r_3\ve_3+2-2e_3)(r_3\ve_3+2e_3), \\
	& G_{32} := (r_3\ve_3+1-r_3)(\ve_3+2)(r_3\ve_3+2-2e_3)(r_3\ve_3+2e_3),\\
    & G_{33} := (r_3\ve_3+1-r_3)(\ve_3+2)(r_3\ve_3+2-2r_3c_3)(\ve_3+2c_3).
\end{split}
\end{equation*}

System \eqref{bue.k31} has three invariant subspaces, namely, $r_3=0$, $\ve_3=0$ and their intersection.
Recall that by definition $e=e_3$ and thus $0<e_3<1$.
\begin{enumerate}
\item $r_3=\ve_3=0$:
 in this case the dynamics is governed by
\begin{equation}\label{bue.k32}
\begin{split}
	c_3'& = -32c_3e_3(1-e_3)[2 + c_3(\gamma-e_3)]\\
    e_3'&=0.
\end{split}
\end{equation}
When $e_3>\gamma$ the equilibria of the system are $p_3^a=(r_3,c_3,e_3,\ve_3)=(0, 0, e_3, 0)$ and $p_3^r=(r_3,c_3,e_3,\ve_3)=(0, \frac{2}{e_3-\gamma}, e_3, 0)$.
Note that the point $p_3^a$ is attracting for the flow in the plane  $(c_3, e_3)$, while the point $p_3^r$ is repelling.
\begin{remark}
Note that when $e_3\to\gamma$, the point $p_3^r\to\infty$ and is not visible any more in the chart $K_3$, see Fig. \eqref{chartK31}.
\end{remark}
%%%%--------------------
\item $\ve_3=0$ and $r_3\geq 0$:
In the invariant plane $\ve_3=0$, the dynamics is governed by
\begin{equation}\label{bue.k33}
\begin{split}
	& r_3' = r_3 c_3 \left[\gamma-e_3\right]V(r_3, c_3, e_3),\\
	& c_3' = c_3\left[(4r_3-2) - c_3(\gamma-e_3)\right]V(r_3, c_3, e_3),\\
    & e_3' = 2r_3 c_3 \left[(2r_3c_3-1)\right]V(r_3, c_3, e_3),
\end{split}
\end{equation}
where $V(r_3, c_3, e_3) = 32e_3(1-r_3)(1-e_3)(1-r_3c_3)$. 
Recall that $c=r_3c_3$ and therefore $V(r_3, c_3, e_3)>0$.
The equilibria of the system are the plane $c_3=0$, denoted by $\hjk{\bar\S_{0,2}}$, and the curve of equilibria given by $c_3=\frac{2}{e_3-\gamma}$, denoted by $\mathcal{M}_3^r$.
The change of stability of the points in $\hjk{\bar\S_{0,2}}$ occurs at $r_3=0.5$, i.e., for $r_3<0.5$ the points are attracting, while for $r_3>0.5$ they are repelling.
We denote the attracting part of $\hjk{\bar\S_{0,2}}$ by $\hjk{\bar\S_{0,2}^\txta}$.
The $e_3$-axis, which we denote by $\ell_{e_3}$, is a boundary of $\hjk{\bar\S_{0,2}^\txta}$, which is a line of equilibria.

\item $r_3=0$ and $\ve_3\geq 0$:
In the invariant plane $r_3=0$, the dynamics is governed by
\begin{equation}\label{bue.k34}
\begin{split}
	& e_3' =0,\\
	& c_3' = -8c_3e_3(1-e_3)\left[(\gamma(\ve_3+2)-2e_3)(\ve_3+2c_3) + 4(\ve_3+2)\right],\\
    & \ve_3' = -8\ve_3e_3(1-e_3)\left[\gamma(\ve_3+2)-2e_3\right](\ve_3+2c_3).
\end{split}
\end{equation}
The equilibria of the system are the planes $c_3=0$, and the line $\ve_3=\frac{2(e_3-\gamma)}{\gamma}$, denoted by $N_3^0$.
The Jacobian of \eqref{bue.k34} along the curve $N_3^0$ has the eigenvalues 
\begin{equation}\label{bue.eigk3}
\lambda_{31}=-64e_3(c_3+1)(1-e_3), \qquad
\lambda_{32}=-8\gamma \ve_3 e_3(1-e_3)(\ve_3+2c_3),
\end{equation}
and hence $N_3^0$ is fully attracting.
In fact, $N_3^0$ is exactly the critical manifold $N_2^0$ that we found in chart $K_2$.
In other words, $N_3^0$ is the image of $N_2^0$ under the transformation $\kappa_{23}$, defined in \eqref{k2tok3}.
\begin{remark}
Note that the attracting manifold $N_2^0$ that is unbounded in chart $K_2$, is now bounded in chart $K_3$.
So the behavior of the critical manifold that is not visible in chart $K_2$ when $e\to\gamma$, is now visible in chart $K_3$.
For $e_3=\gamma$, the critical manifold $N_3^0$ intersects the line $\ell_{e_3}$ at the non-hyperbolic point $q_{e_3}=(e_3,c_3,\ve_3)=(\gamma, 0,0)$.
\end{remark}
% \begin{figure}[H] 
% \centering
% \begin{tikzpicture}
% \pgftext{\includegraphics[scale=1]{figures/chartK32}}
% \node at (3.4, -0.8){$c_3$};
% \node at (-0.6, 0.8){$p_3^a$};
% \node at (0.1, 0){$q_{e_3}$};
% \node at (0.8, 1.1){$p_3^r$};
% \node at (1.4, 1.8){$\mathcal{M}_3^r$};
% \node at (-3.3, -0.3){$N_3^0$};
% \node at (-3.1, -3.1){$\ve_3$};
% \node at (-0.2, 3.4){$e_3$};
% \node at (-0.45, 2.8){1};
% \node at (-1.15, 1.7){$\ell_{e_3}$};
% \end{tikzpicture}
% \caption{Dynamics of system \eqref{bue.k31} in the invariant plane $r_3=0$.
% The slow manifold $N_3^0$ in purple, the exit point $q_{e_3}$ in black, and the line $\ell_{e_3}$ of equilibria in blue.}
% \label{chartK32_old}
% \end{figure}
\end{enumerate}

\hjk{

\begin{figure}[t]
\centering
\begin{subfigure}[tb]{0.45\textwidth}\centering
\begin{tikzpicture}
8\pgftext{\includegraphics[scale=1.2]{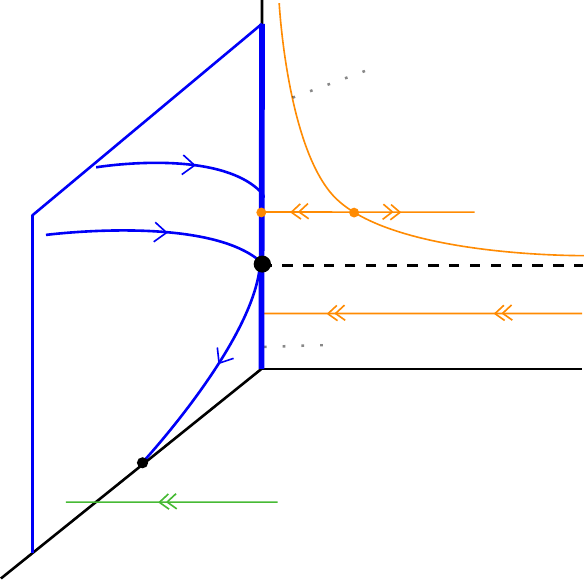}}
\node at (3.5, -1.15){$c_3$};
\node at (-0.7, 0.9){$p_3^a$};
\node at (-0.6, 0.25){$\gamma$};
\node at (0.9, 1.2){$p_3^r$};
\node at (1.3, 2.8){$\mathcal{M}_3^r$};
\node at (0.75, -0.7){$\ell_{e_3}$};
\node at (-3.8, -3.8){$r_3$};
\node at (-0.45, 3.8){$e_3$};
\node at (-2.55, -0.5){$\hjk{\bar\S_{0,2}^\txta}$};
%\node at (-2.4, 0.9){$S_{r,3}^2$};
\node at (-3, -3.55){0.5};
\end{tikzpicture}
\caption[]{Dynamics of \eqref{bue.k31} restricted to $\ve_3=0$. %The attracting plane of equilibria $\hjk{\bar\S_{0,2}^\txta}$ in blue.
}
\label{chartK31}
\end{subfigure}\hfill
\begin{subfigure}[tb]{0.45\textwidth}\centering
\begin{tikzpicture}
\pgftext{\includegraphics[scale=1]{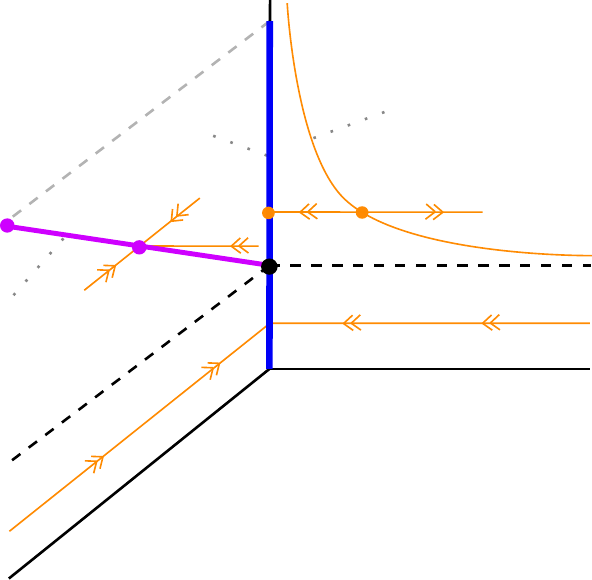}}
\node at (3., -1){$c_3$};
\node at (-0.6, 0.8){$p_3^a$};
\node at (0.1, 0){$q_{e_3}$};
\node at (0.8, 1.1){$p_3^r$};
\node at (1.4, 1.8){$\mathcal{M}_3^r$};
\node at (-3.3, -0.3){$N_3^0$};
\node at (-3.1, -3.1){$\ve_3$};
\node at (-0.2, 3.4){$e_3$};
\node at (-0.45, 2.8){1};
\node at (-1.15, 1.7){$\ell_{e_3}$};
%\draw[white] (-3,-3.5) rectangle (0,-4.4);
\end{tikzpicture}
\caption[]{Dynamics of \eqref{bue.k31} restricted to $r_3=0$.
%The slow manifold $N_3^0$ in purple, the exit point $q_{e_3}$ in black, and the line $\ell_{e_3}$ of equilibria in blue.
}
\label{chartK32}
\end{subfigure}
\caption[]{Dynamics of system \eqref{bue.k31} restricted to invariant subspaces.}
\label{figsK3}
\end{figure}
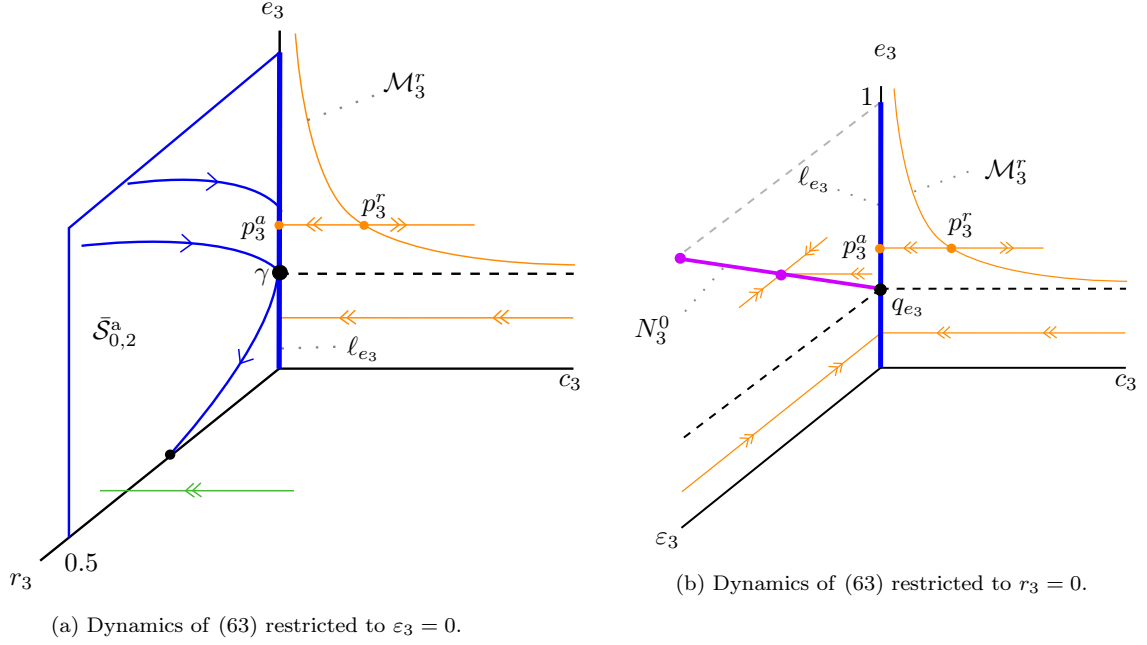
}

%%%%%------------------
% In order to connect the flow in chart $K_2$ with the flow in chart $K_3$, let us define the sections $\Delta_3^{in}$ and $\Delta_3^{out}$ as follows 

% %%%%%%%------------------------------
% \vspace{2cm}
% Let us define the sections $\Delta_3^{in}$ and $\Delta_3^{out}$ as follows 
% \begin{align*}
% 	& \Delta_3^{in} = \{(r_3, c_3, e_3, \ve_3) \ | \ c\in\left(0, \frac{\delta_2}{r_3}\right), \ e\in(\gamma-\delta_2, \gamma+\delta_2), \ \ve_3=\frac{\ve}{r_3}\},\\
% 	& \Delta_3^{out} =  \{(r_3, c_3, e_3, \ve_3) \ | \ 0<f< p_f , c\in(0, \delta_3), \ e\in(0, \delta_4)      \},
% \end{align*}
% where $\delta_3, \delta_4>0$ is sufficiently small.
% %%%%%%----------------------------------
%
We summarize the analysis of the invariant planes, performed in this subsection, in the following Lemma.
%%%%%--------------------------
\begin{lemma}
The following properties hold for system \eqref{bue.k31}:
\begin{enumerate}
\item The equilibria are the plane $\hjk{\bar\S_{0,2}}$ which intersects the line $\ell_{e_3}$, 
and the following two one-dimensional manifolds
\begin{align*}
& \mathcal{M}_3^r=\{(r_3.c_3,e_3,\ve_3) \ | \  r_3=\ve_3=0, \ e_3\in(\gamma, 1), \ c_3=\frac{2}{e_3-\gamma}\},\\
& N_3^0=\{(r_3.c_3,e_3,\ve_3) \ | \  r_3=c_3=0, \ e_3\in[\gamma, 1), \ \ve_3=\frac{2(e_3-\gamma)}{\gamma}\}.
\end{align*}

\item For $e_3>\gamma$, the equilibria of system \eqref{bue.k31} along $N_3^0$ have
\begin{enumerate}
\item a two-dimensional stable manifold corresponding to the negative eigenvalues given in \eqref{bue.eigk3}.
\item a two-dimensional center manifold corresponding to a double zero eigenvalue.
\end{enumerate}
\item The linearization of the system in $\hjk{\bar\S_{0,2}}$ has a triple zero eigenvalue, and the eigenvalue $\lambda=64e_3(e_3-1)(r_3-1)(r_3-0.5)$ changes its stability at $r_3=0.5$.
\item The linearization of system \eqref{bue.k31} at the steady states in the line $\ell_{e_3}$ has a stable eigenvalue $\lambda = 64e_3(e_3-1)$, and a triple zero eigenvalue.
Moreover, there exists a three-dimensional center manifold $\mathcal{W}_{a,\ve}^c$ at the point $(r_3, c_3, e_3, \ve_3)=(0, 0, e_3, 0)\in\ell_{e_3}$.
In chart $K_3$ close to the point $e_3=\gamma$, the center manifold $\mathcal{W}_{a,\ve}^c$ is given as the graph
\begin{equation}\label{1d.center}
c_3= r_3\ve_3(1+O(r_3\ve_3)).
\end{equation}
\end{enumerate}
\end{lemma}
%%%-----------------
The main goal of chart $K_3$ is to analyze the behavior of the solutions of \eqref{bue.k31} close to the exit point $q_{e_3}\in\ell_{e_3}$.
{\e{From our analysis in chart $K_2$, we know that there exists the family of attracting slow manifolds $\mathcal{N}_2$.
This in chart $K_3$ is denoted by $\mathcal{N}_3$ which is the image of $\mathcal{N}_2$ under the transformation $\kappa_{23}$, i.e. $\mathcal{N}_3=\kappa_{23}(\mathcal{N}_2)$.
In order to know how $\mathcal{N}_3$ is continued close to the point $q_{e_3}$, we restrict the dynamics to the sets}}
\begin{align*}
	& D_3^{in}:=\{(r_3, c_3, e_3, \ve_3) \ | \ r_3\in[0,\alpha_3], e_3\in(\gamma, 1], \ve_3\in[0, \beta_3] \},\\
    & D_3^{out}:=\{(r_3, c_3, e_3, \ve_3) \ | \ r_3\in[0,\alpha_3], e_3\in[0, \gamma), \ve_3\in[0, \beta_3] \},
\end{align*}
where $\alpha_3=\alpha_2\beta_2$ and $\beta_3=\frac{1}{\beta_2}$, due to the transformation $\kappa_{23}$ defined in \eqref{k2tok3}.
Now we define the sections as follows
\begin{align*}
& \Delta_3^{in}:=\{(r_3, c_3, e_3, \ve_3)\in D_3^{in} \ | \ \ve_3=\beta_3\},\\
& \Delta_3^{out}:=\{(r_3, c_3, e_3, \ve_3)\in D_3^{out} \ | \ r_3=\alpha_3\}.
\end{align*}
Let us denote $\Pi_{3}$ as the transition map from $\Delta_3^{in}$ to $\Delta_3^{out}$, induced by the flow of \eqref{bue.k31}.
In order to construct the map $\Pi_{3}$, we reduce system \eqref{bue.k31} to its center manifold, namely, $\mathcal{W}_{a,3}^c$ and analyze the system based on the the dynamics on $\mathcal{W}_{a,3}^c$.
This is done by substituting \eqref{1d.center} into system \eqref{bue.k31}, and rescaling time by dividing out the common factor 
\begin{equation}
\left[r_3\ve_3 + 2 − 2r_3^2\ve_3(1+O(r_3\ve_3)\right]\left[\ve_3 + 2r_3\ve_3(1+O(r_3\ve_3)\right].
\end{equation}
The flow of the center manifold is given by 
\begin{equation}\label{center.manifold31}
\begin{split}
& r_3' = r_3 G_{34},\\%(r_3,e_3,\ve_3),\\
& e_3' = r_3 G_{35},\\%(r_3,e_3,\ve_3),\\
& \ve_3' = -\ve_3 G_{34},%(r_3,e_3,\ve_3),
\end{split}
\end{equation}
where we denote
\begin{equation*}
\begin{split}
    & G_{34}:= \left[\gamma(1-r_3)(\ve_3+2) - 2e_3(r_3\ve_3+1-r_3)\right] (r_3\ve_3+2-2e_3)(r_3\ve_3+2e_3),\\
   & G_{35}:= \left[8r_3^2\ve_3(1+O(r_3\ve_3))(1-e_3)(r_3\ve_3+2e_3) - 4e_3(r_3\ve_3+2-2e_3)\right](r_3\ve_3+1-r_3)(\ve_3+2).
\end{split}
\end{equation*}
It is clear from \eqref{center.manifold31} that the planes $r_3=0$ and $\ve_3=0$ are invariant.
Setting $r_3=0$ in \eqref{center.manifold31}, one obtains 
\begin{equation}\label{center.manifold32}
\begin{split}
& e_3' = 0,\\
& \ve_3' = -4\ve_3 e_3(1 - e_3)[\gamma(\ve_3+2) - 2e_3].
\end{split}
\end{equation}
The equilibria of \eqref{center.manifold31} are again the line $\ell_{e_3}$ and the manifold $N_3^0$.
The Jacobian of \eqref{center.manifold32} evaluated at the line $\ell_{e_3}$ has the eigenvalue $\lambda=8e_3(1-e_3)(e_3 - \gamma)$, showing that $\ell_{e_3}$  is repelling for $e_3>\gamma$, while attracting for $e_3<\gamma$.
Moreover, the manifold $N_3^0$ is attracting for the flow in the plane $r_3=0$.
The eigenvalue at the point $(r_3,e_3,\ve_3)=(0, \gamma, 0)\in\ell_{e_3}$ is zero and hence this point is degenerate.

Setting $\ve_3=0$ in system \eqref{center.manifold31} results in
\begin{equation}
\begin{split}
& r_3' = 8 r_3 e_3(1-e_3)(1-r_3)[\gamma-e_3],\\
& e_3' = -16r_3e_3(1-e_3)(1-r_3).
\end{split}
\end{equation}
In the plane $\ve_3=0$, the line $\ell_{e_3}$ is attracting for $e_3>\gamma$ while repelling for $e_3<\gamma$.
\begin{remark}
Note that the dynamics in the invariant plane $\ve=0$ corresponds to the reduced flow on $S^2_{a}$ in the original system.
\end{remark}
Summarizing the analysis, we have the following lemma.
\begin{lemma}\label{nilpotent}
The following properties hold for system \eqref{center.manifold31}:
\begin{enumerate}
\item The curve $N_3^0$ has a one dimensional stable manifold, and a two dimensional center manifold away from the point $q_{e_3}$.
\item The linearization of \eqref{center.manifold31} at the points in $\ell_{e_3}$ is given by 
\begin{equation*}
\begin{pmatrix}
8e_3(e_3-1)(e_3-\gamma) & 0 & 0\\
16e_3(e_3-1) & 0 & 0\\
0 & 0 & -8e_3(e_3-1)(e_3-\gamma)
\end{pmatrix},
\end{equation*}
\item The point $q_{e_3}$ is nilpotent.
\end{enumerate}
\end{lemma}
\begin{figure}[h]
\centering
\begin{tikzpicture}
\pgftext{\includegraphics[scale=1]{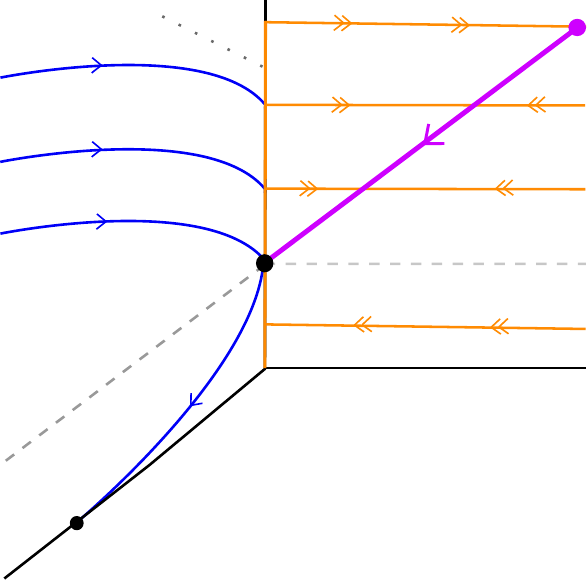}}
\node at (3.3, -0.9){$\ve_3$};
\node at (0.0, 0.0){$q_{e_3}$};
\node at (-3.2, -3.1){$r_3$};
\node at (-0.3, 3.2){$e_3$};
\node at (-0.45, 2.7){1};
\node at (-1.65, 2.95){$\ell_{e_3}$};
\node at (1.8, 2.25){$N_3^0$};
\end{tikzpicture}
\caption{Dynamics of the system \eqref{center.manifold31}; the attracting manifold $N_3^0$ in purple, and the nilpotent point $q_{e_3}$ in black.}
\label{chartK33}
\end{figure}
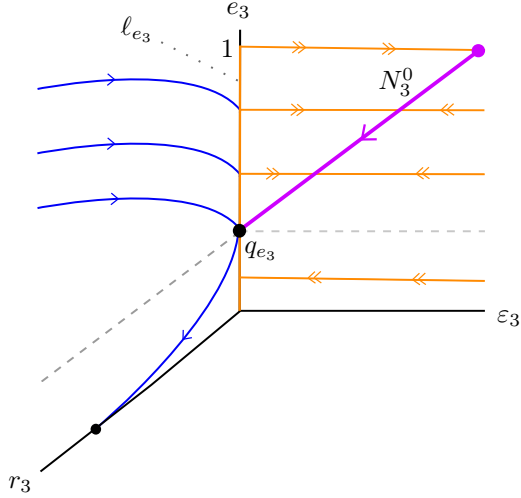
%%%%%%-----------
As we already mentioned, our goal in chart $K_3$ is to describe the dynamics \eqref{bue.k31} close to the line $\ell_{e_3}$, and especially at the point $q_{e_3}$.
To this end, we defined the map $\Pi_{3}:\Delta_3^{in}\to \Delta_3^{out}$ where $\Delta_3^{in}$ is transversal to $N_3^0$ for $e>\gamma$, while $\Delta_3^{out}$ is transversal to the slow manifold in the plane $\ve_3=0$ for $e<\gamma$.
From Lemma \ref{nilpotent} we know that the point $q_{e_3}$ is nilpotent.
Thus, in order to describe the transition map $\Pi_3$ we need to blow-up the point $q_{e_3}$. 
{\e{For such a point, a similar analysis has been carried out in \cite[Theorem 5.8]{kosiuk2012relaxation}, in view of which we have the following theorem.
\begin{theorem}\label{transition.K3}%%%\label{second.blowup}
Assume that $R_3\subset\Delta_3^{in}$ is a small rectangle centered at the intersection point $N_3^0\cap\Delta_3^{in}$.
For sufficiently small $\alpha_3$, the transition map $\Pi_3:R_3\to\Delta_3^{out}$ induced by the flow of \eqref{center.manifold31} is well-defined and satisfies the following properties:
\begin{enumerate}
\item The continuation of $\mathcal{N}_3$ by the flow intersects the section $\Delta_3^{out}$ in a curve denoted by $\sigma_3^{out}$.
\item Restricted to the lines $r_3=constant$ in $R_3$, the map is contracting with the rate $\exp(-K/r_3)$ for some $K>0$.
\item The image $\Pi_3(R_3)$ is an exponentially thin wedge-like containing the curve $\sigma_3^{out}$. 
\end{enumerate}
\end{theorem}}}
% In view of Theorem \ref{second.blowup}, we have the following lemma that describes the the transition map $\Pi_{3}$.
% \begin{lemma}\label{transition.K3}
% The transition map $\Pi_{3}:\Delta_3^{in}\to\Delta_3^{out}$ is well-defined by the flow of \eqref{center.manifold31}, and is given by
% \begin{equation*}
% \Pi_{3}
% \begin{pmatrix}
% r_3\\
% c_3\\
% e_3\\
% \beta_3
% \end{pmatrix}
% =
% \begin{pmatrix}
% \alpha_3\\
% h_{a,3}^{out}(\alpha_3, e_3^{out}, \ve_3^{out}) + \Phi(r_3,c_3,e_3)\\
% e_3^{out}(\ve_3^{out})\\
% \ve_3^{out}
% \end{pmatrix},
% \end{equation*}
% where $h_{a,3}^{out}(\cdot)$ describes the three dimensional center manifold $\mathcal{W}_{a,3}^c$, and $e_3^{out}(\ve_3^{out})$ describes the curve $\sigma_3^{out}\subset\Delta_3^{out}$, which is formed by the intersection of $\mathcal{N}_3^0$ with $\Delta_3^{out}$.  
% \end{lemma}

{\e{Finally, if we set $\alpha_3=\delta_2$ (recall the definition of $\Sigma_2$) we actually have that $\Delta_3^{out}=\bar{\Sigma}_2:=\Phi^{-1}(\Sigma_2\times\{[0,\rho_2]\})$ for some $\rho_2>0$.}}
%%%%%%------------------------------
\begin{figure}[t]
\centering
\begin{tikzpicture}
\pgftext{\includegraphics[scale=1.2]{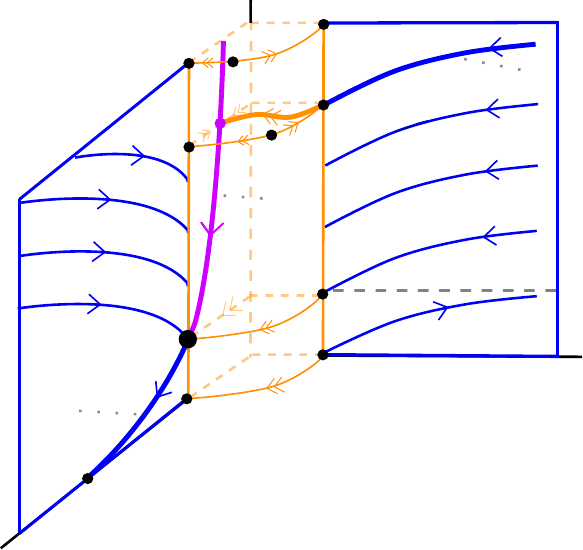}}
\node at (3.8, -1){$\bar{c}$};
\node at (-1.6, -0.85){$\bar{q}_e$};
\node at (2.95, 2.5){$\mathcal{\bar{\omega}}_5$};
\node at (1.7, 3.35){$\hjk{\bar\S_{0,1}^\txta}$};
\node at (-2.7, 2){$\hjk{\bar\S_{0,2}^\txta}$};
\node at (-3.8, -3.4){$\bar{f}$};
\node at (-0.5, 3.7){$\bar{e}$};
\node at (0.6, 1.9){$\bar{q}^e$};
\node at (-2.9, -1.6){$\bar{\omega}_7$};
\node at (0, 0.95){$\bar{\mathcal{N}}^0$};
\end{tikzpicture}
\caption{Geometry of the blown-up space and the singular cycle close to the non-hyperbolic line $\ell_1$, which is blown-up to the orange cylinder.
The reduced flows in $\hjk{\bar\S_{0,1}^\txta}$ and $\hjk{\bar\S_{0,1}^\txta}$ are illustrated in blue. 
The thick orange manifold inside the cylinder corresponds to the three-dimensional center manifold in chart $K_1$.
The attracting critical manifold  in chart $K_2$ is shown in purple.}
\label{CompleteblowupE}
\end{figure}
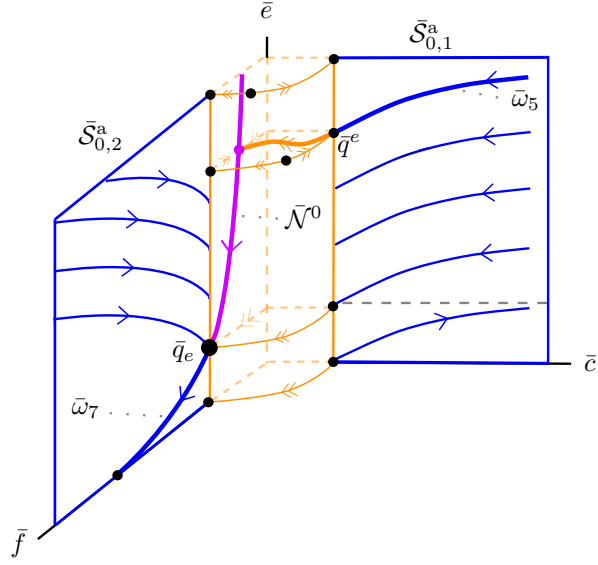

\subsubsection{Properties of the blow-up of the non-hyperbolic line $\ell_1\times\{0\}$ and proof of Lemma \ref{entry1}}
In the above subsections, we have presented the detailed analysis of the blow-up of the non-hyperbolic line $\ell_1\times\{0\}$ in charts $K_1, K_2$ and $K_3$, which has been summarized in Fig. \ref{CompleteblowupE}.
{\color{black}{A summary of the analysis, carried out in such charts, is as follows.
First of all, the critical manifolds \hjk{$\S_{0,1}$} (i.e., $f=0$) and \hjk{$\S_{0,2}$} (i.e., $c=0$) intersect in the non-hyperbolic line $\ell_1$, which is replaced by the orange cylinder, see Figs. \ref{blowup.space} and \ref{CompleteblowupE}.
Note that in Fig. \ref{CompleteblowupE}, the orbits $\bar{\omega}_5$ and $\bar{\omega}_7$ in the blown-up space correspond to the orbits $\omega_5$ and $\omega_7$, respectively.
The point at which $\bar{\omega}_5$ reaches the cylinder is denoted by $\bar{q}^e$, and the point at which $\bar{\omega}_7$ starts is denoted by $\bar{q}_e$.  
Starting from the section $\bar{\Sigma}_1$, the trajectory follows the orbit $\bar{\omega}_5$ on $\bar{f}=0$ until it reaches the point $\bar{q}^e$.
Our analysis in chart~$K_1$ (Lemma~\ref{centermanifold.k1e}) shows that there exists a three-dimensional attracting center manifold which is the continuation of the family of orbits (indexed by $\ve$) of the attracting slow manifold $\hjk{\S_{\ve,1}^\txta}$.
This allows us to connect the family $\hjk{\S_{\ve,1}^\txta}$ into the chart $K_2$ which is inside the cylinder (see the thick orange manifold from $\bar{q}^e$ to $\bar{\mathcal{N}}^0$ in  Fig. \ref{CompleteblowupE}).
Our analysis in chart $K_2$ (Lemma~\ref{lem:crit2}) shows that the slow manifold $N_2^0$ is normally hyperbolic and stable.
Therefore, the family $\hjk{\S_{\ve,1}^\txta}$ is exponentially attracted by the slow manifold $N_2^\ve$.
Next, our analysis in chart $K_3$ shows that the unbounded critical manifold $N_2^0$ (see Figs. \ref{chartK2e}, \ref{chartK33}) limits in the point $q_{e_3}$, which is exactly the point $\bar{q}_e$ in Fig. \ref{CompleteblowupE}.
Next, our analysis in chart $K_3$ (see Lemma~\ref{nilpotent} and Fig.~\ref{chartK33}) demonstrates that the unbounded critical manifold $N_2^\ve$ (see Figs.~\ref{chartK2e} and~\ref{chartK33}) limits at the point~$q_{e_3}$, which is exactly the point $\bar{q}_e$ in Fig.~\ref{CompleteblowupE}.
In addition, we have proven that the point~$q_{e_3}$ is degenerate, i.e., the linearization of the dynamics at $q_{e_3}$ has a nonzero (stable) eigenvalue and a triple zero eigenvalue (see Lemma~\ref{nilpotent}), which allows us to construct a three-dimensional center manifold at the point $q_{e_3}$.
Now, by following the family $\mathcal{N}_2$ along such a center manifold, we conclude (Lemma \ref{transition.K3}) that the continuation of $N_2^\ve$ for a sufficiently small $\ve>0$ intersects the section $\bar{\Sigma}_2$ in a point, namely, $(\alpha_3, c_3(\ve_3), e_3(\ve_3), \ve_3)\in\bar{\Sigma}_2$, for some $\ve_3\in[0,\beta_3]$, which is exponentially close to the slow manifold $\hjk{\S_{\ve,2}^\txta}$.
Note that the point $(\alpha_3, c_3(\ve_3), e_3(\ve_3), \ve_3)$ converges to the point $q_2:=\Sigma_2\cap\omega_7$ as $\ve_3\to 0$.
All these analyses in charts $K_1, K_2$, and $K_3$ show that the transition map $\bar{\pi}_{1}:\bar{\Sigma}_1\to\bar{\Sigma}_2$ is well-defined for $\ve\in[0,\ve_0]$ and is smooth for $\ve\in(0,\ve_0]$, for some $\ve_0>0$.  
}}

%Note that due to the blow-up, the non-hyperbolic line $\ell_1$ has been replaced by the cylinder, depicted in orange color, see Fig. \ref{blowup.space}.
We are now ready to prove Lemma \ref{entry1}.

\textit{Proof of Lemma \ref{entry1}.}
The proof is carried out by constructing the map $\pi_{1}:\Sigma_1\to\Sigma_2$ for $\ve>0$ as
\begin{equation}
	\pi_{1}=\Phi\circ\bar{\pi}_1\circ\Phi^{-1},
\end{equation}
where $\Phi$ is given by \eqref{trans.e}, $\Phi^{-1}$ is the corresponding blown-up transformation, and $\bar{\pi}_{1}:\bar{\Sigma}_1\to\bar{\Sigma}_2$ is a transition map which can equivalently be regarded as
\begin{equation}\label{eq:pibar}
\bar\pi_1(\bar{\Sigma}_1) = \Pi_3\circ\kappa_{23}\circ\Pi_2\circ\kappa_{12}\circ\Pi_1(\bar{\Sigma}_1)\subset\bar{\Sigma}_2
\end{equation}

The proof is based on the corresponding transition map $\bar{\pi}_{1}:\bar{\Sigma}_1\to\bar{\Sigma}_2$ in the blown-up space and interpreting the result for fixed $\ve\in[0, \ve_0]$ with $\ve_0>0$. 
Recall that the transition $\bar{\pi}_{1}:\bar{\Sigma}_1\to\bar{\Sigma}_2$ is equivalent to the transition map $\pi_{1}:\Sigma_1\to\Sigma_2$ in the sense that it has the same properties. 
Furthermore, via the matching maps $\kappa_{ij}$ defined in Lemma \ref{lem:kappas}, we have appropriately identified the relevant sections in each of the charts, allowing us to follow the flow of the blown-up vector field along the three charts. 
%Thus, we shall continue this proof by just detailing the transition \eqref{eq:pibar}. 
%$\bar{\pi}_{1}:\Delta_1^{in}\to\Delta_3^{out}$ which can equivalently be regarded as
% \begin{equation}
% \bar\pi_1(\Delta_1^{in}) = \Pi_3\circ\kappa_{23}\circ\Pi_2\circ\kappa_{12}\circ\Pi_1(\Delta_1^{in})\subset\Delta_3^{out}
% \end{equation}
%
%We start by explaining the properties of the blown-up vector fields, summarized in Fig. \ref{CompleteblowupE}.

{\e{As summarized above, the transition map $\bar{\pi}_{1}:\bar{\Sigma}_1\to\bar{\Sigma}_2$ is well-defined for $\ve\in[0,\ve_0]$ and smooth for $\ve\in(0,\ve_0]$ for some $\ve_0>0$.}}
It remains to prove that $\bar{\pi}_{1}$ is a contraction.
From Lemma \ref{contract.s1} we know that the solutions started in $\bar{\Sigma}_1$ are contracting, see (Fig. \ref{redVF1}).
This family of orbits is continued to chart $K_2$ by spending an $O(1)$-time on the time scale $t_2$ of system \eqref{bue.k22}.\
This continuation persists (Theorem \ref{transition.K3}) during the passage near the point $q_{e_3}$ in chart $K_3$ until it reaches the section $\bar{\Sigma}_2$.
As the contraction persists from $\bar{\Sigma}_1$ to $\bar{\Sigma}_2$, one concludes that $\bar{\pi}_{1}$ is a contraction.
This completes the proof. \qed

%% file: subfiles/BlowUpFaxis.tex
%%% blow-up the f-axis
\renewcommand{\bf}{F}
\newcommand{\be}{E}
\newcommand{\bc}{C}
\newcommand{\br}{R}
\newcommand{\beps}{\varepsilon}

In this subsection, for the sake of brevity, we summarize the blow-up of the non-hyperbolic line $\ell_2\times\{0\}$ and give a sketch of the proof of Lemma \ref{entry2}.
To this end, we transform the non-hyperbolic line of steady states $\ell_2\times\{0\}$ by
\begin{equation}\label{trans.f}
	f = \tilde{f}, \qquad\qquad c=r\tilde{c}, \qquad\qquad \ve=r\tilde{\ve},\qquad\qquad e=r\tilde{e},
\end{equation}
where {\e{$\tilde{c}^2 + \tilde{e}^2 + \tilde{\ve}^2 = 1$, $\bar{f}\in[0,1]$}} and $r\geq 0$, and define the charts $\tilde{K}_1, \tilde{K}_2$ and $\tilde{K}_3$ as follows
\begin{align*}
& \tilde{K}_1: \qquad f = \tilde{f}_1, \qquad\qquad  c = \tilde{r}_1 \tilde{c}_1, \qquad\qquad \ve = \tilde{r}_1 \tilde{\ve}_1,\qquad\qquad e = \tilde{r}_1,\\ %\label{k1f}\\
& \tilde{K}_2: \qquad f = \tilde{f}_2, \qquad\qquad c = \tilde{r}_2 \tilde{c}_2, \qquad\qquad \ve = \tilde{r}_2,\qquad\qquad\quad  e = \tilde{r}_2 \tilde{e}_2,\\ %\label{k2f}\\
& \tilde{K}_3: \qquad f = \tilde{f}_3, \qquad\qquad c = \tilde{r}_3, \qquad\qquad\quad \ve = \tilde{r}_3 \tilde{\ve}_3,\qquad\qquad e = \tilde{r}_3 \tilde{e}_3. %\label{k3f}
\end{align*}

$\tilde K_1=\left\{ \tilde e=1\right\}$

%\begin{lemma}
% The changes of coordinates $\tilde{\kappa}_{ij}$ from the chart $\tilde{K}_i$ to $\tilde{K}_j$ $(i,j=1,2,3,  \ i\neq j)$ are given by
% \begin{align}
% 	& \tilde{\kappa}_{12}:\qquad
%     \tilde{f}_2 =  \tilde{f}_1, \qquad\qquad 
%     \tilde{c}_2 =  \frac{\tilde{c}_1}{\tilde{\ve}_1},\qquad\qquad 
%     \tilde{\ve}_2 = \tilde{r}_1 \tilde{\ve}_1, \qquad\qquad 
%     \tilde{e}_2 = \frac{1}{\tilde{\ve}_1}, \qquad \tilde{\ve}_1>0,\label{tk1totk2}\\
%     & \tilde{\kappa}_{23}:\qquad 
%     \tilde{f}_3 = \tilde{f}_2, \qquad\qquad 
%     \tilde{c}_3= \frac{1}{\tilde{c}_2},\qquad\qquad 
%     \tilde{\ve}_3 = \frac{1}{\tilde{c}_2}, \qquad\qquad\ \  
%     \tilde{e}_3 = \frac{\tilde{e}_2}{\tilde{c}_2}, \quad\  
%     \tilde{c}_2>0.\label{tk2totk3}
% \end{align}
% \end{lemma}
%%%%%-----------------------
%%%%%%%----------------------------
%\subsubsection{Properties of the blow-up of the non-hyperbolic line $\ell_2\times\{0\}$ and proof of Lemma \ref{entry2}}
%In this subsection, we present a sketch of the proof of Lemma \ref{entry2}.
%%%%%%%%%%%%%%-----------------------------------------------
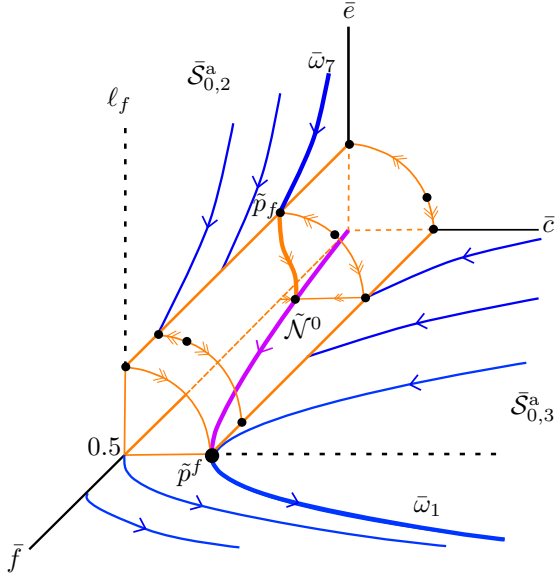
\begin{figure}[t]
\centering
\begin{tikzpicture}
\pgftext{\includegraphics[scale=1.15]{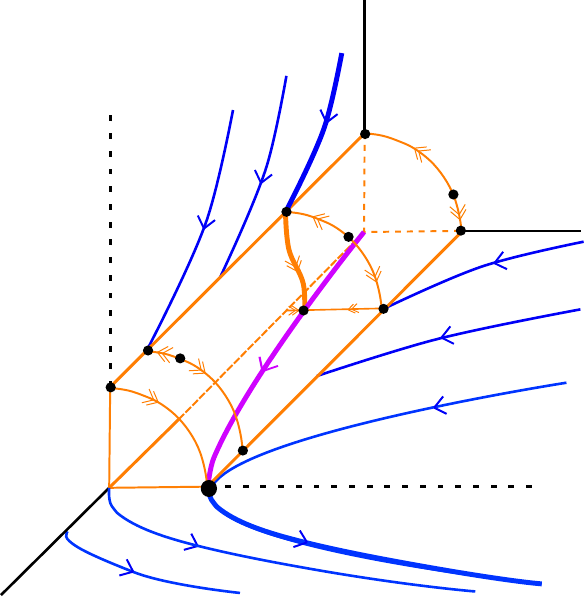}}
\node at (-0.25, 1.1){$\tilde{p}_f$};
\node at (-1.25, -2.45){$\tilde{p}^f$};
\node at (-3.6, -3.6){$\bar{f}$};
\node at (3.5, .85){$\bar{c}$};
\node at (.85, 3.7){$\bar{e}$};
\node at (0.5, 3){$\bar{\omega}_7$};
\node at (1.9, -2.9){$\bar{\omega}_1$};
\node at (-2.4, -2.1){$0.5$};
\node at (-2.2, 2.5){$\ell_f$};
\node at (0.25, -0.5){$\tilde{\mathcal{N}}^0$};
\node at (-1, 2.8){$\hjk{\bar\S_{0,2}^\txta}$};
\node at (3.3, -1.6){$\hjk{\bar\S_{0,3}^\txta}$};
\end{tikzpicture}
\caption{Geometry of the blown-up space and the singular cycle close to the non-hyperbolic line $\ell_2$, which is blown-up to the orange cylinder.
The reduced flows in $\hjk{\bar\S_{0,2}^\txta}$ and $\hjk{\bar\S_{0,3}^\txta}$ are shown in blue. 
The thick orange manifold inside the cylinder corresponds to the three-dimensional center manifold in chart $\tilde{K}_1$.
The attracting critical manifold  in chart $\tilde{K}_2$ is illustrated in purple.}
\label{CompleteblowupF}
\end{figure}

Recall that the goal of Lemma \ref{entry2} is to describe the map $\pi_{2}:\Sigma_2\to\Sigma_3$ in the original space.
In this subsection, we present a sketch of the proof of Lemma \ref{entry2} by constructing the corresponding map $\bar{\pi}_{2}:\bar{\Sigma}_2\to \bar{\Sigma}_3$ in the blown-up space, and interpreting the results for fixed $\ve\in[0,\ve_0]$ for some $\ve_0>0$. 
For the sake of brevity, we have summarized the analysis of the blow-up of the non-hyperbolic line $\ell_2\times\{0\}$ in Fig.~\ref{CompleteblowupF}.

First of all, note that the non-hyperbolic line $\ell_2$, which is the intersection of the critical manifolds $c=0$ and $e=0$, has been blown-up to the orange cylinder (see Fig.~\ref{CompleteblowupF}).
We have illustrated the slow flows in the planes $c=0$ and $e=0$ in blue.
Note that the orbits $\bar{\omega}_7$ and $\bar{\omega}_1$ which are in the blown-up space correspond, respectively, to the orbits $\omega_7$ and $\omega_1$ in the original space (see Figs.~\ref{SinCyc} and~\ref{CompleteblowupF}).
As it is shown in Fig.~\ref{CompleteblowupF}, the intersection of the cylinder with $\bar{\omega}_7$ and $\bar{\omega}_1$ is denoted by $\tilde{p}_f$ and $\tilde{p}^f$, respectively.

Our analysis in chart $\tilde{K}_1$ proves that there exists a three-dimensional attracting center manifold at the point $\tilde{p}_f$, which is the continuation of the family indexed by $\ve$ of the attracting slow manifold $\hjk{\bar\S_{0,2}^\txta}$.
In view of such a center manifold, the family of the slow manifold $\hjk{\bar\S_{0,2}^\txta}$ enters the chart $\tilde{K}_2$.
Our analysis in chart $\tilde{K}_2$ proves that there exists a hyperbolic attracting one-dimensional slow manifold $\tilde{\mathcal{N}}^0$, which attracts the interior of the cylinder.
Our analysis in chart $\tilde{K}_3$ shows that the critical manifold $\tilde{\mathcal{N}}^0$ limits at the point $\tilde{p}^f$ (see Fig.~\ref{CompleteblowupF}).
Note that $\tilde{p}^f$ is a \textit{degenerate} point, i.e., the linearization of the blown-up dynamics in chart $\tilde{K}_3$ at $\tilde{p}^f$ has a stable eigenvalue, and a triple zero eigenvalues which allows us to construct a three-dimensional attracting center manifold.
Therefore the family of flows follows such a center manifold and then intersects the section $\Sigma_3$ in a point $(f(\ve), \delta_3, e(\ve))$, for some $\delta_3>0$, which is exponentially close to the slow manifold $\hjk{\S_{\ve,3}^\txta}$ and converges to the point $q_3:=\Sigma_3\cap\omega_1$ as $\ve\to 0$. 
This proves that the transition map $\bar{\pi}_{2}:\bar{\Sigma}_2\to \bar{\Sigma}_3$ and hence $\pi_{2}:\Sigma_2\to\Sigma_3$ are well-defined for $\ve\in[0,\ve_0]$ and also are smooth for $\ve\in(0,\ve_0]$, for some $\ve_0>0$.
%Using Fenichel theory, all orbits starting in $\Sigma_2$ are exponentially attracted onto the continuation of the slow manifold $\hjk{\bar\S_{0,2}^\txta}$, and then enter chart $\tilde{K}_2$ where all orbits are exponentially attracted onto the 1-dimensional slow manifold $\tilde{\mathcal{N}}_2^\ve$.
%Next, they follow $\tilde{\mathcal{N}}_2^\ve$ through $\tilde{K}_2$ to $\tilde{K}_3$, pass the point $\tilde{p}^f$ and hence intersect $\Sigma_3$ exponentially close to the point $(f(\ve), \delta, e(\ve))$.
The proof of contraction of the transition map $\pi_{2}$ follows the same line of reasoning as that of the map $\pi_1$, and hence is omitted for brevity.

%% file: subfiles/discuss.tex
In this work we have studied a model that describes several important properties of myxobacteria during development \cite{igoshin2004biochemical}.
This model, which is in line with observation from experiments \cite{igoshin2004biochemical}, acts as an internal clock to control the gliding motions in myxobacteria.
When two cells collide with each other, the speed of the clock in both cells is affected, some spatial wave patterns are created, and hence leads to synchronization of cells, i.e., fruiting body.
The model presented in \cite{igoshin2004biochemical} can reproduce observed spatial patterns in experiments, and furthermore, it can explain both the cellular oscillations and the developmental stage of myxobacteria from vegetative swarming to the rippling phase and hence to the formation of the fruiting body.

The model, described by a system of three ordinary differential equations, has oscillatory behavior for certain parameter values, and sufficiently small Michaelis-Menten constants which we have unified them by a parameter $\ve$.
We have analyzed the dynamics of this oscillator in the limits of $\ve$, and proven that for sufficiently small $\ve$, there exists a strongly attracting limit cycle.
%Note that from biological point of view, the uniqueness of a periodic solution admits that the oscillations are more robust than if multiple periodic solutions coexist, because in the latter conditions the solutions may switch from one to another mode of oscillations in the same experimental conditions.
%So the existence of a unique solution biologically means that the domain of parameters on which the solutions are robust is much larger than the case that multiple solutions coexist.
{\e{
The geometric method could be pushed to analyze the global uniqueness of the limit cycle 
which is clearly of great interest from both the mathematical and biological point of view.
This requires a more global analysis of the singular flows, and in particular, connecting orbits between the critical manifolds $\hjk{\S_{0,i}}$ by orbits of the layer problem. 
As the layer problem is linear, this is possible.}}
Our approach has been based on the geometric perturbation analysis and blow-up method.
The geometric perturbation theory and geometric desingularization by several blow-ups allow us to fully understand the structure of the limit cycle.
We emphasize that the approach and tools presented in this paper, i.e. geometric singular perturbation theory and the blow-up method, are not limited to the analysis of the system \eqref{syseps};
these tools can be applied to similar systems such as \cite{kut2009analytical} whose parameters have the property of zero-order ultrasensitivity.

%% file: subfiles/RelaxGamma.tex
%%%%%---------------------------------------
%{\color{blue}
%As discussed in two-parameter bifurcation analysis in Subsection \ref{2para.bif}, system \eqref{syseps} has oscillatory behavior when $\gamma\in(0,1)$ and $\ve$ is sufficiently small, see Fig. \ref{bif.2para}.
Although in our analysis we have fixed the parameter $\gamma=0.08$, in this appendix we show that the behavior of the singular cycle $\Gamma_0$, illustrated in Fig. \ref{SinCyc}, will remain qualitatively the same under a sufficiently small perturbation of $\gamma$, and hence Theorem 3.22 holds for these values as well.

As it is shown in Fig. \ref{SinCyc}, $\omega_1$ and $\omega_3$ are described by the layer problem \eqref{identical}, whose behavior highly depends on parameter $\gamma$.  
For the former, observations from numerical simulations show that if the layer problem starts from a point in $\hjk{\S_{0,3}^\txta}$, namely $p_\gamma:=(f_\gamma, \frac{1}{2}, 0)\in\ell_e$ when $\frac{1}{2}<f_\gamma<1$, the parameter $\gamma$ can influence the fast dynamics to move from $p_\gamma$ to a point either in $\hjk{\S_{0,1}^\txta}$, $\hjk{\S_{0,2}^\txta}$, $\hjk{\S_{0,3}^\txta}$, $\hjk{\S_{0,4}^\txta}$, $\hjk{\S_{0,5}^\txta}$ or $\hjk{\S_{0,6}^\txta}$. 
For the latter, numerical simulations show that if the layer problem starts from a point in $\hjk{\S_{0,6}^\txta}$, namely $p^\gamma:=(f^\gamma, \frac{1}{2}, 1)\in\ell^e$ when $0<f^\gamma<\frac{1}{2}$, the parameter $\gamma$ can influence the fast dynamics to move from $p^\gamma$ to a point either in $\hjk{\S_{0,1}^\txta}$ or $\hjk{\S_{0,2}^\txta}$.
In this regard, it is interesting to find a certain range for $\gamma$ such that the qualitative behavior of the fast dynamics remains the same as $\omega_1$ and $\omega_3$, shown in Fig. \ref{SinCyc}, i.e., the fast dynamics moves \emph{directly} from $p_\gamma$ to a point in $\hjk{\S_{0,6}^\txta}$, and from $p^\gamma$ to a point in $\hjk{\S_{0,1}^\txta}$, while does not intersect with the other planes.

In view of equations \eqref{Red.S33}, one can show that the slow flow, started from the point $p^f=(\frac{1}{2}, 0, 0)$, will reach the point $p_1=(\frac{1+\sqrt{\gamma}}{2}, \frac{1}{2}, 0)$.
In order to find a certain range for $\gamma$, as the layer problem \eqref{identical} is linear, one can find the closed form of solutions.
In view of the boundary conditions in $\hjk{\S_{0,3}^\txta}$ and $\hjk{\S_{0,6}^\txta}$, we will get a system of transcendental equations, whose solution determines a point at which $\omega_1$ intersects with $\hjk{\S_{0,6}^\txta}$.
However, due to the fact that it is impossible to solve such a system of equations analytically, we have used numerical methods to calculate the solution of transcendental equations.
Computed numerically, for any $\gamma\in\mathcal{R}_1:=(0.0561, \, 0.1177)$, the qualitative behavior of the fast dynamics is the same as $\omega_1$, illustrated in Fig. \ref{SinCyc}.
Moreover, one can check that for such a range, the qualitative behavior of the fast dynamics is the same as $\omega_3$ as well, illustrated in Fig. \ref{SinCyc}.
Therefore, one concludes that Theorem \ref{main.result} is valid for all $\gamma\in\mathcal{R}_1$.
Analogously, one can find a range for the case that $\gamma$ is close to 1, see Remark \ref{gamma.relax}. 

%% file: main.bbl
\begin{thebibliography}{10}

\bibitem{brons2015mixed}
Morten Br{\o}ns, Mathieu Desroches, and Martin Krupa.
\newblock Mixed-mode oscillations due to a singular hopf bifurcation in a
  forest pest model.
\newblock {\em Mathematical Population Studies}, 22(2):71--79, 2015.

\bibitem{chicone2006}
Carmen Chicone.
\newblock {\em Ordinary differential equations with applications}.
\newblock Springer Science \& Business Media, 2006.

\bibitem{chow1994normal}
Shui-Nee Chow, Chengzhi Li, and Duo Wang.
\newblock {\em Normal forms and bifurcation of planar vector fields}.
\newblock Cambridge University Press, 1994.

\bibitem{decroly1982birhythmicity}
Olivier Decroly and Albert Goldbeter.
\newblock Birhythmicity, chaos, and other patterns of temporal
  self-organization in a multiply regulated biochemical system.
\newblock {\em Proceedings of the National Academy of Sciences},
  79(22):6917--6921, 1982.

\bibitem{dhooge2003matcont}
Annick Dhooge, Willy Govaerts, and Yu~A Kuznetsov.
\newblock Matcont: a matlab package for numerical bifurcation analysis of odes.
\newblock {\em ACM Transactions on Mathematical Software (TOMS)},
  29(2):141--164, 2003.

\bibitem{dumortier1996canard}
Freddy Dumortier and Robert~H Roussarie.
\newblock {\em Canard cycles and center manifolds}, volume 577.
\newblock American Mathematical Soc., 1996.

\bibitem{fenichel1979geometric}
Neil Fenichel.
\newblock Geometric singular perturbation theory for ordinary differential
  equations.
\newblock {\em Journal of Differential Equations}, 31(1):53--98, 1979.

\bibitem{goldbeter1991minimal}
Albert Goldbeter.
\newblock A minimal cascade model for the mitotic oscillator involving cyclin
  and cdc2 kinase.
\newblock {\em Proceedings of the National Academy of Sciences},
  88(20):9107--9111, 1991.

\bibitem{goldbeter2017dissipative}
Albert Goldbeter.
\newblock Dissipative structures and biological rhythms.
\newblock {\em Chaos: An Interdisciplinary Journal of Nonlinear Science},
  27(10):104612, 2017.

\bibitem{goldbeter1984ultrasensitivity}
Albert Goldbeter and Daniel~E Koshland.
\newblock Ultrasensitivity in biochemical systems controlled by covalent
  modification. interplay between zero-order and multistep effects.
\newblock {\em Journal of Biological Chemistry}, 259(23):14441--14447, 1984.

\bibitem{goldbeter1996biochemical}
Albert Goldbeter and John~J Tyson.
\newblock Biochemical oscillations and cellular rhythms: The molecular bases of
  periodic and chaotic behaviour.
\newblock {\em Nature}, 380(6571):213--213, 1996.

\bibitem{goodwin1965oscillatory}
Brian~C Goodwin.
\newblock Oscillatory behavior in enzymatic control processes.
\newblock {\em Advances in enzyme regulation},
  3:425IN1429IN3431--428IN2430IN6437, 1965.

\bibitem{guckenheimer2013nonlinear}
John Guckenheimer and Philip~J Holmes.
\newblock {\em Nonlinear oscillations, dynamical systems, and bifurcations of
  vector fields}, volume~42.
\newblock Springer Science \& Business Media, 2013.

\bibitem{gucwa2009geometric}
Ilona Gucwa and Peter Szmolyan.
\newblock Geometric singular perturbation analysis of an autocatalator model.
\newblock {\em Discrete and Continuous Dynamical Systems / S, v.2, 783-806
  (2009)}, 2009.

\bibitem{hodgkin1952quantitative}
Alan~L Hodgkin and Andrew~F Huxley.
\newblock A quantitative description of membrane current and its application to
  conduction and excitation in nerve.
\newblock {\em The Journal of physiology}, 117(4):500--544, 1952.

\bibitem{igoshin2004biochemical}
Oleg~A Igoshin, Albert Goldbeter, Dale Kaiser, and George Oster.
\newblock A biochemical oscillator explains several aspects of myxococcus
  xanthus behavior during development.
\newblock {\em Proceedings of the National Academy of Sciences of the United
  States of America}, 101(44):15760--15765, 2004.

\bibitem{jardon2019survey}
Hildeberto Jard{\'o}n-Kojakhmetov and Christian Kuehn.
\newblock A survey on the blow-up method for fast-slow systems.
\newblock {\em arXiv preprint arXiv:1901.01402}, 2019.

\bibitem{kaiser1989multicellular}
Dale Kaiser.
\newblock Multicellular development in myxobacteria.
\newblock In {\em Genetics of bacterial diversity}, pages 243--263. Elsevier,
  1989.

\bibitem{kosiuk2012relaxation}
Ilona Kosiuk.
\newblock {\em Relaxation oscillations in slow-fast systems beyond the standard
  form}.
\newblock PhD thesis, Univ. Leipzig, 2012.

\bibitem{kosiuk2016geometric}
Ilona Kosiuk and Peter Szmolyan.
\newblock Geometric analysis of the goldbeter minimal model for the embryonic
  cell cycle.
\newblock {\em Journal of mathematical biology}, 72(5):1337--1368, 2016.

\bibitem{krupa2001extending}
Martin Krupa and Peter Szmolyan.
\newblock Extending geometric singular perturbation theory to nonhyperbolic
  points---fold and canard points in two dimensions.
\newblock {\em SIAM journal on mathematical analysis}, 33(2):286--314, 2001.

\bibitem{kuehn2015multiple}
Christian Kuehn.
\newblock {\em Multiple time scale dynamics}, volume 191.
\newblock Springer, 2015.

\bibitem{kuehn2015multiscale}
Christian Kuehn and Peter Szmolyan.
\newblock Multiscale geometry of the olsen model and non-classical relaxation
  oscillations.
\newblock {\em Journal of Nonlinear Science}, 25(3):583--629, 2015.

\bibitem{kut2009analytical}
Carmen Kut, Vahid Golkhou, and Joel~S Bader.
\newblock Analytical approximations for the amplitude and period of a
  relaxation oscillator.
\newblock {\em BMC systems biology}, 3(1):6, 2009.

\bibitem{kuznetsov2013elements}
Yuri~A Kuznetsov.
\newblock {\em Elements of applied bifurcation theory}, volume 112.
\newblock Springer Science \& Business Media, 2013.

\bibitem{novak2008design}
B{\'e}la Nov{\'a}k and John~J Tyson.
\newblock Design principles of biochemical oscillators.
\newblock {\em Nature reviews Molecular cell biology}, 9(12):981--991, 2008.

\bibitem{roberts2016mixed}
Andrew Roberts, John Guckenheimer, Esther Widiasih, Axel Timmermann, and
  Christopher~KRT Jones.
\newblock Mixed-mode oscillations of el nino--southern oscillation.
\newblock {\em Journal of the Atmospheric Sciences}, 73(4):1755--1766, 2016.

\bibitem{szmolyan2004relaxation}
Peter Szmolyan and Martin Wechselberger.
\newblock Relaxation oscillations in $\mathbb{R}^3$.
\newblock {\em Journal of Differential Equations}, 200(1):69--104, 2004.

\bibitem{taghvafard2018modeling}
Hadi Taghvafard.
\newblock {\em Modeling, Analysis, and Control of Biological Oscillators}.
\newblock University of Groningen, 2018.

\bibitem{taghvafard2017parameter}
Hadi Taghvafard, Hildeberto Jard{\'o}n-Kojakhmetov, and Ming Cao.
\newblock Parameter-robustness analysis for a biochemical oscillator model
  describing the social-behaviour transition phase of myxobacteria.
\newblock {\em Proc. R. Soc. A}, 474(2209):20170499, 2018.

\bibitem{takens1976constrained}
Floris Takens.
\newblock Constrained equations: a study of implicit differential equations and
  their discontinuous solutions.
\newblock {\em Lecture Notes in Mathematics}, 525:143--235, 1976.

\bibitem{van1928lxxii}
Balth van~der Pol and Jan van~der Mark.
\newblock The heartbeat considered as a relaxation oscillation, and an
  electrical model of the heart.
\newblock {\em The London, Edinburgh, and Dublin Philosophical Magazine and
  Journal of Science}, 6(38):763--775, 1928.

\bibitem{winfree1967biological}
Arthur~T Winfree.
\newblock Biological rhythms and the behavior of populations of coupled
  oscillators.
\newblock {\em Journal of theoretical biology}, 16(1):15--42, 1967.

\end{thebibliography}
